\documentclass[11pt]{article}
\usepackage{amsmath,amssymb,amscd,verbatim,color}
\usepackage{bbm}
\usepackage{stmaryrd}
\usepackage{bbding}
\usepackage{mathrsfs}
\usepackage{amsfonts}
\usepackage{cases}
\usepackage{amsthm}
\usepackage[square, comma, sort&compress, numbers]{natbib}

\usepackage[top=1in,bottom=1in,left=1.20in,right=1.20in]{geometry}

\def\N{{\mathbb N}} \def\Z{{\mathbb Z}}
\def\R{{\mathbb R}} \def\T{{\mathbb T}}

\newtheorem{Theorem}{Theorem}
\newtheorem{Lemma}{Lemma}[section]

\newtheorem{Remark}{Remark}[section]

\newtheorem{remark}{Remark}[section]

\newcommand{\la}{\langle } \newcommand{\ra}{\rangle }

\numberwithin{equation}{section}
\theoremstyle{plain}

\newcommand{\Cal}{\mathcal}
\newcommand{\sss}{\smallskip}
\newcommand{\ms}{\medskip}
\newcommand{\bs}{\bigskip}

\newcommand{\kth}{e^{{\rm i}\langle k, \theta\rangle}\; }

\makeatletter

\newcommand{\Rmnum}[1]{\expandafter\@slowromancap\romannumeral #1@}
\makeatother
\newcommand{\beq}{\begin{equation} } \newcommand{\eeq}{\end{equation} }

\begin{document}
	\title{{Linearly Stable KAM Tori for One Dimensional Forced Kirchhoff Equations under Periodic Boundary Conditions}
}

\author{\\ Yin Chen$^1$\thanks{The first author is partially supported
by NSFC grant 12301236.},  Jiansheng Geng$^2$\thanks{The corresponding author. The second author is partially supported
by NSFC grant 12571199.}, Guangzhao Zhou$^2$\\
{\footnotesize $^1$ School of Mathematical Sciences, Yangzhou University, Yangzhou 225002, P.R.China}\\
{\footnotesize$^2$ School of Mathematics, Nanjing University, Nanjing 210093, P.R.China}\\
{\footnotesize Email: yinchen@yzu.edu.cn(Chen); jgeng@nju.edu.cn(Geng); 602023210023@smail.nju.edu.cn(Zhou)}}

%\date{put in custom date; delete this line for today's date}
%\maketitle

\date{}
\maketitle

\textbf{Abstract}:
 We prove an abstract infinite  dimensional KAM theorem, which could be applied to prove the existence and linear stability of small-amplitude  quasi-periodic solutions for one dimensional forced Kirchhoff equations with periodic boundary conditions
 \[ u_{tt}-(1+\int_{0}^{2\pi} |u_x|^2 dx)u_{xx}+
 M_\xi u+\epsilon g(\bar{\omega}t,x) =0,\quad  u(t,x+2\pi)=u(t,x),\]
 where $M_\xi$ is a real Fourier multiplier, $g(\bar{\omega}t,x)$ is real analytic  with forced Diophantine frequencies $\bar\omega$, $\epsilon$ is a small parameter.
 The paper generalizes the previous results from the simple eigenvalue to the double eigenvalues under the quasi-linear perturbation.

\textbf{Keywords}: Kirchhoff equation;  KAMPDE;  T\"{o}plitz-Lipschitz; double eigenvalues.\\

\section{Introduction and Main Results }

\ \ \ \   Kirchhoff equation has been introduced for the first time in \cite{Ki} in 1876 in one space dimension, without forcing term and with Dirichlet boundary conditions
     which describes the transversal free vibrations of a clamped string with the tension on the deformation. It is a quasi-linear  wave-type PDE (partial differential equation) with unbounded nonlinearity, namely the nonlinear part of the equation contains as many derivatives as the linear part. We distinguish the quasi-periodic solutions according to the following two cases: the corresponding quasi-periodic solutions are called {\it response}  solutions if one only excites the forced frequencies; the corresponding quasi-periodic solutions are called {\it non-response} quasi-periodic solutions (quasi-periodic solutions for short) if one excites the internal frequencies. For PDEs with unbounded nonlinearities, Kuksin firstly proved
the existence of quasi-periodic solutions for KdV in \cite{K2} (see also Kappeler-P\"oschel \cite{KP1}).
This approach has been improved by Liu-Yuan \cite{LY2}
to deal with DNLS (Derivative Nonlinear Schr\"odinger) (see also \cite{Gao}). We mention that
Corsi-Feola-Procesi \cite{CFP} establish a general abstract KAM method to
prove the existence of analytic solutions of quasi-linear PDEs.
Besides, the response solutions for quasi-linear (either
fully nonlinear) PDEs have been proved by Baldi-Berti-Montalto  \cite{Ba3} for perturbations of Airy equations,
by Feola-Procesi \cite{FP} for fully nonlinear reversible
Schr\"odinger equation. The quasi-periodic solutions for quasi-linear PDEs have been proved by Baldi-Berti-Montalto \cite{Ba4,Ba5} for perturbations of KdV and mKdV equations. Berti-Montalto \cite{Be7} have proved the existence of quasi-periodic standing wave solutions of the water
waves equations with surface tension and Baldi-Berti-Haus-Montalto \cite{BBHM} have proved the similar results without surface tension. Feola-Giuliani \cite{FG}  has established  the small amplitude, quasi-periodic traveling waves for the pure gravity water waves system in infinite depth.
 Such results are all obtained by imposing the second Melnikov conditions and provide the linear stability of the solutions. See also Baldi-Montalto \cite{BM} and Berti-Hassainia-Masmoudi\cite{BHM} for Euler equation case.

Besides, by imposing only the first Melnikov conditions, the existence of response solutions and quasi-periodic solutions can be also proved
with the multi-scale approach. This method called CWB method comes from Nash-Moser iteration scheme developed by Craig-Wayne \cite{CW},
Bourgain \cite{B1,B2,B4} for analytic NLS (Nonlinear Schr\"odinger) and NLW (Nonlinear Wave). This approach is based on the
 multi-scale analysis of the linearized operators around the quasi-periodic solutions and it has been recently improved by Berti-Bolle \cite{Be2, Be3, Be4} for NLW, NLS with smooth nonlinearity, by Berti-Corsi-Procesi \cite{Be6} on compact Lie-groups and recently by Wang \cite{Wang} for energy supercritical nonlinear Schr\"odinger equations. This
method does not provide any information about the linear stability of the quasi-periodic solutions since the
linearized equations have variable coefficients. Comparing \cite{Be2} with \cite{Be3}, we should realize that there is a big difference between response solution case and quasi-periodic solution case.

Indeed the second  Melnikov
conditions are seriously violated in the case of multiple eigenvalues for one space dimension and higher space dimension. There are very few results about linear stability of quasi-periodic solutions, for example, Chierchia-You
\cite{Ch}, for analytic one dimensional NLW equation with periodic boundary conditions and Geng-Yi \cite{GY5}, Geng-You \cite{GY3} for analytic one dimensional Schr\"odinger equation with periodic boundary conditions (double eigenvalues), people can refer to Kuksin \cite{K1}, Kuksin-P\"oschel \cite{KP}
 and P\"oschel \cite{P2} for simple eigenvalue case. Geng-You \cite{GY1, GY2} proved that the higher dimensional nonlinear beam equations and nonlocal smooth Schr\"odinger equations
admit small-amplitude linearly-stable quasi-periodic solutions. Chen-Geng-Xue \cite{CGX} proved that the higher dimensional nonlinear wave equations under nonlocal and forced perturbation admit small-amplitude linearly-stable quasi-periodic solutions (see also \cite{CG}). The breakthrough of
constructing quasi-periodic solutions for higher dimensional Schr\"odinger
equation by modified KAM method was made recently by Eliasson-Kuksin \cite{E2}. They
proved that the higher dimensional nonlinear Schr\"odinger equations admit small-amplitude
linearly-stable quasi-periodic solutions. Eliasson and Kuksin introduced the conception of
the Lipschitz domain, in the Lipschitz domain, the corresponding normal frequencies satisfy
T{\"o}plitz-Lipschitz property, thus the measure estimates are feasible (see also \cite{E-G-K,Ge,GXY,GY6,P4}).

The existence of response periodic solutions for the forced Kirchhoff equation in any space dimension has been proved
by Baldi \cite{Ba1}, both for Dirichlet boundary conditions and for periodic boundary conditions.
This approach does not imply the linear stability and it does not work in the
quasi-periodic case since the small divisor problem is more difficult. More recently the existence and linear stability of response solutions in one space dimension under the periodic boundary conditions has been proved by Montalto \cite{M}, and the existence of response solutions for the forced Kirchhoff equation in higher space dimension has been proved by Corsi-Montalto \cite{CM}, but they did not provide the linear stability. Moreover, they didn't excite the internal frequencies, i.e., they only handled the forced frequency as the exciting frequency. In \cite{CG1}, Chen--Geng proved that the higher dimensional Kirchhoff equations without the forced perturbation admit small--amplitude linearly--stable quasi--periodic solutions, where the pair--property of the normal coordinates $w_n, \bar w_n$ is crucial in that paper. In \cite{CG2}, Chen-Geng excited the internal frequencies and
prove the existence and linear stability of quasi-periodic  solutions for one dimensional forced Kirchhoff equation under Dirichlet boundary conditions.
Compared to Montalto\cite{M} and Corsi-Montalto\cite{CM}, Chen-Geng\cite{CG2}  is based on an improved Kuksin lemma together with the refined T\"oplitz-Lipschitz property, while  \cite{M} and \cite{CM} are based on KAM methods together with pseudo-differential calculus. In addition, the obtained solutions in \cite{M} and \cite{CM} are $C^k$ ($k$ finite), while the obtained solutions in \cite{CG2} are at least $C^\infty$ even Gevrey smooth. Compared to \cite{CG1}, the pair--property of the normal coordinates $w_n, \bar w_n$ is seriously violated in the forced perturbation, hence, Chen-Geng\cite{CG2}  developed off--diagonal decay property of the forced perturbation together with the refined T\"oplitz-Lipschitz property. In this paper, we generalizes Chen-Geng\cite{CG2} from Dirichlet boundary conditions to periodic boundary conditions, which will bring the essential difficulties. As is well known, the eigenvalues associated with Dirichlet boundary conditions are simple, while the eigenvalues associated with periodic boundary conditions are double, together with quasi-linear perturbation, KAM theory for this kind of partial differential equations is more difficult. In fact, we make use of the pair--property of the normal coordinates $w_n, \bar w_n$ along each KAM iteration, i.e., the pair--property of the normal coordinates $w_n, \bar w_n$ along each KAM iteration is preserved (which need to be clarified), the contribution of the finite-rank perturbation to the normal form $N$ is constant-coefficient  non-diagonal $2\times 2$ block, i.e., the different normal coordinates $w_n, \bar w_{-n}(|n|\leq EK)$ is coupled, we can handle them with the help of the finite-rank perturbation.

Considering back the forced  Kirchhoff equation under periodic boundary conditions
\begin{equation}\label{1.2}
  u_{tt}-(1+\int_{0}^{2\pi} |u_x|^2 dx)u_{xx}+M_\xi u+\epsilon g(\bar{\omega}t,x) =0,\quad  u(t,x+2\pi)=u(t,x),
\end{equation}
 it is a quasi-linear PDE so we could not directly apply the so-called Kuksin lemma in \cite{KP1, LY2} to obtain an abstract KAM theory. A critical strategy for proving the existence and linear stability of small-amplitude quasi-periodic solutions of $(\ref{1.2})$ is to keep the pair--property of the normal coordinates $w_n, \bar w_n(|n|>EK)$ along each KAM iteration and decay property of the nonlinear term (\ref{G})(see also $(A5)$), which will always be preserved throughout the KAM iteration. Hence it is feasible for us to further develop  and establish an abstract KAM theory to prove our results. Moreover, the refined T\"oplitz-Lipschitz property $(A6)$ will also be verified at each KAM step, which is critical to solve the homological equations and estimate the measure of the parameter set. Once the assumption $(A6)$ has been satisfied, we can consequently obtain the form of each normal frequency $\Omega_n$ satisfying $(\ref{asymp1})$, where the function $f$ only depends on the angle variable $\theta$ and parameter $\sigma$, namely $f$ is uniform in each space index $n$.

 In fact, $\tilde \Omega_n$ in $( \ref{asymp1})$ comes from the coefficients of the second-order terms $w_n\bar w_n$ which can not be eliminated in the KAM iteration. Specifically, in the subsection 4.1, after the initial KAM iteration, we observe that all these second-order terms $w_n\bar w_n$  originate from the two aspects. One is directly from the second-order term $w_n\bar w_n$ in $P^1$ (see $(A5)$) which can not be eliminated. In this case, the second Melnikov conditions are like
 $$ |\la k,\omega\ra\pm 2\bar \Omega_n|\geq \frac{\gamma_0\cdot |n|}{K^\tau}, $$
 coming from the special form of the Kirchhoff equation and  we can obtain one more regularity from these denominators such that the unbounded terms can be controlled when solving the homological equations. Furthermore, due to the (\ref{G}), the coefficients of $w_n\bar w_n$ obviously have the same order as $|n|$. The other is from $P^3$ (see $(A5)$) which is the result of the Poisson brackets
  $$\{ P-R, F\}, \{\{ P-R, F\}, F\}, \cdots,$$
  where $R$ in $(\ref{R0})$ and $F$ defined in $(\ref{F0})$. Among these Poisson brackets, the terms $w_nw_m$, $w_n\bar w_m$,$\bar w_n\bar w_m$, $|n|,|m|\leq EK$  appear in $P^3$  and can be eliminated in each KAM iteration except for $w_n\bar w_m, |n|=|m|\leq EK$ since their coefficients are always bounded thanks to the exponential decay property in Lemma 3.2, which is related to coefficients of the first-order term $w_n,\bar w_n$ in (\ref{decay}). In this case, when solving the homological equations, the second Melnikov conditions are like
  $$|\la k,\omega\ra \pm(\bar \Omega_n+\bar \Omega_m)|\geq \frac{\gamma_0}{K^\tau},\ |n|,|m|\le EK,$$
$$|\la k,\omega\ra \pm(\bar \Omega_n-\bar \Omega_m)|\geq \frac{\gamma_0}{K^\tau},\ |n|,|m|\le EK,\ |k|+||n|-|m||\neq 0.$$
  Besides, all the coefficients of $w_n\bar w_m, |n|=|m|\leq EK$ come from the coefficients of
  $(w_n+\bar {w}_n)^2(w_m+\bar {w}_m)^2$ multiplied by the coefficients of
  $w_n,\bar w_n$ in $F$. Due to (\ref{G}),  the coefficients of the fourth-order terms $(w_n+\bar {w}_n)^2(w_m+\bar {w}_m)^2$ in $P-R$ have the same order as $|n||m|$.  By Lemma 3.2 and the construction of the Hamiltonian function $F$ in (\ref{F0}), the coefficients of $w_n,\bar w_n$ in $F$ inherit the exponential decay $ e^{-|n|\bar \rho}$ of the coefficients of the term $w_n,\bar w_n$ in (\ref{decay}), then $ e^{-|n|\bar \rho}$ can be used to control $|n|$. So it is natural for us to compute
\begin{eqnarray}\label{coeff}
\frac{\partial^2P}{\partial w_n\partial w_n}+\frac{\partial^2P}{\partial w_n\partial \bar w_n}+\frac{\partial^2P}{\partial \bar w_n\partial \bar w_n}, n\in\Bbb Z,
\end{eqnarray}
which include all the possible coefficients of $w_n\bar w_n$, namely $\tilde \Omega_n$,$n\in\Bbb Z$.  Due to the above discussion, (\ref{coeff}) have the same order as $|n|$ so the factor $\frac{1}{|n|}$ is used to eliminate  the number $|n| $ --the effect of the quasi-linear perturbation. Therefore it is necessary to prove the T\"{o}plitz-Lipschitz property, namely
\[\|\lim_{n\rightarrow \infty}\frac{1}{|n|}\sum_{\upsilon=\pm}\frac{\partial^2 P}{\partial w^\upsilon_n\partial w^\upsilon_n}\|_{\!{}_{D(r,s),\Cal O}}\leq \varepsilon,\]
\[\|\frac{1}{|n|}\sum_{\upsilon=\pm}\frac{\partial^2 P}{\partial w^\upsilon_n\partial w^\upsilon_n}-\lim_{n\rightarrow \infty}\frac{1}{|n|}\sum_{\upsilon=\pm}\frac{\partial^2 P}{\partial w^\upsilon_n\partial w^\upsilon_n}\|_{\!{}_{D(r,s),\Cal O}}\leq \frac{\varepsilon}{|n|},n\in\Bbb Z\]\\
the second inequality indicates the uniform decay of the drift of the normal frequencies.
 According to the above discussion, it is sufficient for us to impose the non-resonance conditions for the
 difference between two normal frequencies and the second Melnikov non-resonance
 conditions defined in the assumption $(A3)$ have two kinds of formulas according to the size of $n$.
  Moreover, the perturbation $P$ can be divided into three parts with the special form defined in the assumption $(A5)$.
 In this paper using only the KAM scheme is  more convenient than that in \cite{M},
 where the authors made use of pseudo-differential calculus together with quadratic KAM reduction.

Specifically, here we assume that the operator $A:=-\partial_{xx} +M_\xi$ with periodic boundary conditions has eigenvalues $\{\mu_n\}$
satisfying
\[\tilde{\omega}_j=\lambda_j=\sqrt{\mu_{i_j}}=\sqrt{i_j^2+\xi_j},\ 1\leq j\leq b;\
\Omega_n=\lambda_n=\sqrt{\mu_n}=|n|,\ n\neq i_1,\ldots,i_b,\]
 and the corresponding orthonormal basis of eigenfunctions $\{\phi_n(x)\}\in L^2(\T),n\in \Z.$ For the sake of convenience, we choose real eigenfunctions $\phi_n(x)$ as follows:
 \beq\label{eigenfunction}\phi_n(x)=\left\{\begin{array}{l}\sqrt{\frac{1}{2\pi}}, n=0\\
 \sqrt{\frac{1}{\pi}}\sin(nx), n>0 \\
\sqrt{\frac{1}{\pi}}\cos(nx), n<0
\end{array}
\right.\eeq
We assume $0\in\{i_1,\ldots,i_b\}$
in order to take care of $(\mu_n,k)=(0,0)$,
and we assume the parameter $\sigma=(\bar\omega,\xi)\in\Cal{O}\subset\R^{\nu+b}$, where $\xi=(\xi_1,\ldots,\xi_b)\in (0,1)^b\subset \R^b$, $\Cal{O}$ is a compact subset.\\

Now we state the main theorem as follows.

 \begin{Theorem}\label{maintheorem}
For any $0< \gamma \ll 1$, there is a Cantor subset $\Cal
O_\gamma\subset\Cal O$ with ${\rm meas}(\Cal O\setminus \Cal
O_\gamma)=O(\gamma)$, such that for any $(\bar{\omega},\xi)\in \Cal O_\gamma$, equation
$(1.1)$ with the analytic forced term $g(\bar\omega t, x)$, admits a $C^{\infty}$-smooth
small-amplitude, linearly stable
 quasi-periodic solution of the form
\[ u(t, x)=\sum_{n\in\Bbb Z}u_n(\bar{\omega}t,\tilde{\omega}_{1}^*t, \cdots,\tilde{ \omega}_{b}^* t)\phi_n(x),
\] where $u_n: \T^{\nu+b}\to \R $ and
$\tilde{\omega}_{1}^*,\cdots, \tilde{\omega}_{b}^*$ are close to the unperturbed
frequencies $\tilde{\omega}_{1},\cdots, \tilde{\omega}_{b}$.
\end{Theorem}

This paper is organized as follows: In Section 2 we give an infinite dimensional KAM theorem;
  in Section 3, we give its applications
 to the forced  Kirchhoff equations under periodic boundary conditions.
  The proof of the KAM theorem is
  given in the Section
 4, 5, 6. Some technical lemmata are put into the Appendix.

\section {An Infinite Dimensional KAM Theorem for One Dimensional Forced Kirchhoff Equations under Periodic Boundary Conditions}

We start by introducing some notations. For given $b$ vectors $0\in\{i_1, \cdots, i_b\}$ in
$\Z$, denote its complementary set
$\Z_1=\Z\setminus \{ i_1, \cdots, i_b\}$. Let $w=(\cdots,
w_n,\cdots)_{n\in \Z_1}$, and its complex conjugate $\bar
w=(\cdots,\bar w_n,\cdots)_{n\in \Z_1}$. We introduce a Banach space $l^{a,\rho}_1$
with weighted norm
$$\|w\|_{a,\rho} =\sum_{{n\in\Z_1}}|w_n||n|^{a}e^{|n|\rho},$$ where
$a> 0, \rho>0$. Denote a complex neighborhood of $\T^{\nu+b}\times\{I=0\}\times\{w=0\}\times\{\bar w=0\}$ by
$$D(r,s)=\{(\theta,I,w,\bar w):|{\rm Im} \theta|<r,|I|<s^2,{\|w\|}_{a,\rho}<s,
{\|\bar w\|}_{a,\rho}<s\},$$
 where $|\cdot|$ denotes the sup-norm of complex vectors.
Moreover, we denote by $\Cal O$ a positive measure parameter set
in $\R^{\nu+b}:=\R^{\tilde b}$.\\
 A function $F(\theta,\sigma)$ is $C^1_W$ of parameter $\sigma \in \Cal O$  in the sense of whitney and we denote $D(r)=\{\theta: |\textmd{Im}\theta|< r\}$,
$$\|F\|_{D(r),\Cal O}=\sup_{\sigma \in \Cal O}\sup_{\theta\in D(r)}(|F(\theta,\sigma)|+|\frac{\partial}{\partial\sigma}F(\theta,\sigma)|),\quad [F]=\frac{1}{(2\pi)^{\nu+b}}\int_{\T^{\nu+b}}F(\theta,\sigma)d\theta,$$
where if $F$ is independent of $\theta$, then we denote the norm $\|\cdot\|_{D(r),\Cal O}=|\cdot|_{\Cal O}$ for simplicity. For any finite dimensional parameter dependent matrix $A(\sigma)=(a_{ij}(\sigma))$, the matrix norm $\|A\|_{\Cal O}$ is defined by
$$\|A\|_{\Cal O}=\sup_{\sigma\in \Cal O}\max_{i}(\sum_{j}|a_{ij}|+|\frac{\partial}{\partial\sigma}a_{ij}|).$$
 Besides, we introduce a truncation operator $\Gamma_K$ as follows
$$(\Gamma_K F)(\theta):=\sum_{|k|\leq K}\hat{F}_k\kth,\ \ (1-\Gamma_K)F(\theta)=\sum_{|k|> K}\hat{F}_k\kth,$$
where $\hat{F}_k$ is the $k$-Fourier coefficient of $F$.\\
For $F=F(\theta,I,w,\bar w,\sigma)$, we expand $F$ into Taylor series
$$ F(\theta,I,w,\bar w)=\sum_{l\in \Z^b,\alpha,\beta\in \N^{\Z_1}}F_{l,\alpha,\beta}(\theta,\sigma)I^lw^{\alpha}\bar w^{\beta},$$
where $F_{l,\alpha,\beta}$ are $C^1_W$ functions of parameter $\sigma$ in the sense of whitney, $w^\alpha=\Pi_{n\in Z_1}w_n^{\alpha_n}, \bar w^\beta= \Pi_{n\in Z_1}\bar w_n^{\beta_n}$,
$w=(w_n)_{n\in \Z_1}$, $\bar w=(\bar w_n)_{n\in\Z_1}$, $ \alpha, \beta\in \N^{\Z_1}$, $\alpha=(\alpha_n)_{n\in \Z_1}$, $ \beta=(\beta_n)_{n\in \Z_1}$, $\alpha_n\in \N, \beta_n\in \N$. We define the weighted form of function $F$ by
$$ \|F\|_{D(r,s),\Cal O}=\sup_{\|w\|_{a,\rho}< s\atop \|\bar w\|_{a,\rho}< s}\sum_{l,\alpha,\beta}\|F_{l,\alpha,\beta}\|_{D(r),\Cal O}s^{2|l|}|w^\alpha||\bar{w}^\beta|,$$
and the vector $X_{F}=(F_I,-F_\theta,-{\rm i}F_{\bar w},{\rm i}F_w)$ with weighted norm
\begin{eqnarray*}
\|X_F\|_{s,\bar a,\rho,D(r,s),\Cal O}&=&\|F_I\|_{D(r,s),\Cal O}+\frac{1}{s^2}\|F_{\theta}\|_{D(r,s),\Cal O}+\frac{1}{s}\sum_{n\in \Z_1}\|F_{w_n}\|_{D(r,s),\Cal O}|n|^{\bar a}e^{|n|\rho}\\
&+&\frac{1}{s}\sum_{n\in \Z_1}\|F_{{\bar w}_n}\|_{D(r,s),\Cal O}|n|^{\bar a}e^{|n|\rho}.
\end{eqnarray*}
In the analogous way, the norm of the frequencies $\omega=(\omega_j)_{1\leq j\leq \nu+b}$ and semi-norm of $\Omega=(\Omega_n)_{n\in \Z_1}$ are defined as
$$|\omega|_{\Cal O}=\sup_{\sigma\in \Cal O}\sup_{1\leq j\leq \nu+b}(|w_j|+|\frac{\partial \omega_j}{\partial \sigma}|),\quad |\Omega|_{-1, D(r),\Cal O}=\sup_{\sigma \in \Cal O \atop \theta\in D(r)}\sup_{n\in \Z_1}\frac{1}{|n|}|\frac{\partial \Omega_n}{\partial \sigma}|.$$
\begin{Remark}
In this paper, we require that $\bar a=a-1\geq 0$, namely the weight of the vector fields is  weaker than that of $w,\bar w$. This is due to Lemma 3.3.
\end{Remark}
 In this paper,the generalized normal form $N$ depending  on the angle variable $\theta$ is
\begin{equation}\label{hamN}
N =\la\bar{\omega},\bar{I}\ra+\la\tilde{\omega},I\ra+\sum_{n\in \Z_1}\Omega_n(\theta,\sigma)w_n \bar
w_n+\sum_{|n|\leq EK}\langle A_{|n|}z_{|n|},\bar z_{|n|}\rangle,
\end{equation}
where  $\omega=(\bar{\omega},\tilde{\omega}),\theta=(\bar\theta,\tilde\theta)$, $\sigma\in \Cal O$ is a parameter, the phase space is endowed
with the symplectic structure $\displaystyle d\bar{I}\wedge d\bar{\theta}+ dI\wedge d\tilde{\theta}
+{\rm i} \sum_{n\in \Z_1} dw_n \wedge d \bar w_n$. And
$$A_{|n|}= (\overline{A_{|n|}})^T=\left(\begin{array}{ccc}a_{nn}(\sigma)&a_{n(-n)}(\sigma)\\
a_{(-n)n}(\sigma)&a_{(-n)(-n)}(\sigma)
\end{array}\right), z_{|n|}=\left(\begin{array}{c}w_n\\
w_{(-n)}
\end{array}\right),\bar z_{|n|}=\left(\begin{array}{c}\bar w_n\\
\bar w_{(-n)}
\end{array}\right).$$

\sss Now we consider the perturbed Hamiltonian
\begin{equation}\label{hamH}
H=N+P =\la\bar{\omega},\bar{I}\ra+\la\tilde{\omega},I\ra+\sum_{n\in \Z_1}\Omega_n(\theta,\sigma)w_n \bar
w_n+\sum_{|n|\leq EK}\langle A_{|n|}z_{|n|},\bar z_{|n|}\rangle+ P(\theta,I,w,\bar w, \sigma).
\end{equation}
Our goal is to prove that, for most values of parameter
 $\sigma \in \Cal O$ (in Lebesgue measure
sense), the Kirchhoff equations still admit quasi-periodic solutions provided  that $\|X_P\|_{s,\bar a,\rho,D(r,s),\Cal O}$ is sufficiently
small.
\ms

To this end,  we  need  to impose the following conditions on
$\omega(\sigma)$, $\Omega_n(\sigma)$,  $A_{|n|}$ and the perturbation $P$.

\bs \noindent $(A1)${\it Nondegeneracy:} The map $\sigma$ to
$\omega(\sigma)$ is a $C^1_W$ diffeomorphism between $\Cal O$ and its image. Besides, there exists a positive constant $E$ such that $|\omega|_{\Cal O}\leq E$.

\bs \noindent $(A2)${\it Asymptotics of normal frequencies:}
\begin{eqnarray}\label{asymp1} \Omega_n(\theta,\sigma)&=&\bar{\Omega}_n(\sigma)+\tilde
\Omega_n(\theta;\sigma)\nonumber\\
&=&|n|(1+c(\sigma))+|n|f(\theta,\sigma),
\end{eqnarray}
where $\bar{\Omega}_n(\sigma)=[\Omega_n]$, $\tilde\Omega_n(\theta,\sigma)=\Omega_n-\bar \Omega_n$; moreover, set $\mbox{spec}(A_{|n|})=\{d_n, d_{(-n)}\}$, one has
$$|c(\sigma)|_{\Cal O}+\|f(\theta,\sigma)\|_{D(r),\Cal O}+|d_n(\sigma)|_{\Cal O}+|d_{(-n)}(\sigma)|_{\Cal O}=O(\varepsilon_0).$$

\bs \noindent $(A3)${\it Non-resonance conditions:} The frequencies $\omega$ are Diophantine in the sense that there are constants $\gamma_0>0, \tau>\tilde b+2 (\tilde b=\nu+b)$ and an iteration parameter $\frac{\gamma_0}{2}\leq\gamma<\gamma_0$ such that $ |k|\leq K$
\begin{eqnarray*}
&&|\la k,\omega\ra|\geq \frac{\gamma}{|k|^\tau},\quad 0\neq k=(k_1,k_2)\in \Z^{\nu+b}:=\Z^{\tilde b},\\
&&|\la k,\omega\ra\pm (\bar \Omega_n+d_n)|\geq \frac{\gamma_0}{K^\tau},\ |n|\le EK,\\
&&|\la k,\omega\ra \pm((\bar \Omega_n+d_n)+(\bar \Omega_m+d_m))|\geq \frac{\gamma_0}{K^\tau},\ |n|,|m|\le EK,\\
&&|\la k,\omega\ra \pm((\bar \Omega_n+d_n)-(\bar \Omega_m+d_m))|\geq \frac{\gamma_0}{K^\tau},\ |n|,|m|\le EK,\ |k|+||n|-|m||\neq 0,\\
&&|\la k,\omega\ra\pm 2\bar \Omega_n|\geq \frac{\gamma_0\cdot |n|}{K^\tau},\ |n|>EK,
\end{eqnarray*}
where $|k|= \max\{|k_1|,|k_2|\},|k_1|=|k_{1_1}|+\cdots+|k_{1_{\nu}}|,|k_2|=|k_{2_1}|+\cdots+|k_{2_b}|.$

\bs \noindent $(A4)${\it Regularity of the perturbation:} The
perturbation $P$ is {\sl regular} and satisfies
$$\varepsilon:=\|X_P\|_{s,\bar a,\rho, D(r,s),\cal O} \leq \delta^{\frac{1}{1-\beta'}}$$
for some $\delta>0$, $0<\beta'\leq \frac{1}{4}$, $\bar a=a-1$.

\bs \noindent $(A5)${\it Special structure and decay properties of perturbation $P$:}
The perturbation $P=P^1+P^2+P^3$ satisfies a special structure as follows
\begin{eqnarray*}
P=\sum_{\alpha,\beta}P^1_{\alpha\beta}(\theta,I,\sigma)w^\alpha\bar w^{\beta}+\sum_{\alpha,\beta}P^2_{\alpha\beta}(\theta,\sigma)w^\alpha\bar w^{\beta}+\sum_{\alpha,\beta}P^3_{\alpha\beta}(\theta,I,\sigma)w^{\alpha}\bar w^{\beta},
\end{eqnarray*}
with the exponents\\

$\alpha,\beta\in \{\alpha,\beta\in \N^{\Z_1}, \sum\limits_{|n|>EK}\alpha_n+\beta_n>0, \alpha_n+\beta_n\in 2\N, \forall |n|>EK\}$ in $P^1$;

$\alpha,\beta\in \left\{\alpha,\beta\in \N^{\Z_1}, |\alpha+\beta|=\alpha_n+\beta_n=1, \forall  |n|>EK\right\}$ in $P^2$;\\

$\alpha,\beta \in$ $\left\{\alpha,\beta\in \N^{\Z_1},\alpha_n+\beta_n=0, \forall |n|>EK\right\}$ in $P^3$.
When $|n|>EK, \alpha+\beta=e_n$ in $P^2$, we have
\begin{eqnarray}{\label{decayP^2}}
\|P^2_{\alpha\beta}\|_{D(r),\Cal O}\leq c \varepsilon e^{-|n|{\bar\rho}}\quad (\bar\rho>\rho).
\end{eqnarray}

\bs \noindent $(A6)${\it {T\"{o}plitz-Lipschitz property:}}  The following limits exist
$$\|\lim_{n\rightarrow \infty}\frac{1}{|n|}\sum_{\upsilon=\pm}\frac{\partial^2 P}{\partial w^\upsilon_n\partial w^\upsilon_n}\|_{\!{}_{D(r,s),\Cal O}}\leq \varepsilon,$$
moreover, $P$ satisfies for any $ n\in\Z_1$,
$$\|\frac{1}{|n|}\sum_{\upsilon=\pm}\frac{\partial^2 P}{\partial w^\upsilon_n\partial w^\upsilon_n}-\lim_{n\rightarrow \infty}\frac{1}{|n|}\sum_{\upsilon=\pm}\frac{\partial^2 P}{\partial w^\upsilon_n\partial w^\upsilon_n}\|_{\!{}_{D(r,s),\Cal O}}\leq \frac{\varepsilon}{|n|},$$
where $ w^+_n:=w_n, w^-_n:=\bar w_n.$

\bs \noindent $(A7)$ The function $\tilde{\Omega}_n(\theta,\sigma)$ is analytic on some strip $D(r)=\{\theta:|{\rm Im}\theta|<r\}$ around the torus $\T^{\tilde b}$ with $[\tilde{\Omega}_n]=0$ and satisfies
$$ \|\tilde{\Omega}_n\|_{r,2\tau+2,\Cal O}=\sum_{k\in \T^{\tilde b}}|\tilde{\Omega}_{kn}|_{\Cal O}\cdot|k|^{2\tau+2}\cdot e^{|k|r}\leq \delta_0(\gamma_0-\gamma)|n|,\quad  \forall n\in \Z_1,$$
with some constant $\delta_0>0$ and the same $\tau$ as before.

\bs \noindent Now we are ready to state our KAM Theorem.

  \begin{Theorem}\label{KAM}
Assume that $H=N+P$ satisfies $(A1)-(A7)$,
Let $\gamma>0$ small enough, there is a positive constant
$\varepsilon_0=\varepsilon_0(\nu,b,\tau,\gamma,r,s,\rho,\delta_0,E,K)\leq \delta^{\frac{1}{1-\beta'}}$ such
that if $\|X_P\|_{s,\bar a,\rho,D(r,s),\Cal O}=\varepsilon\leq \varepsilon_0$, then the
following holds true: There exist a Cantor set $\Cal
O_\gamma\subset\Cal O$ with ${\rm meas}(\Cal O\setminus \Cal
O_\gamma)=O(\gamma)$ and two maps ( $ C^{\infty}$ in $\theta$ and $C^1_{W}$ in
$\sigma$)
$$\Psi: \T^{\nu+b}\times \Cal O_\gamma\to D(r,s),\ \ \ \
 \omega_*:\Cal O_\gamma\to \R^{\nu+b},$$ where $\Psi$ is close to the trivial embedding
$\Psi_0:\T^{\nu+b}\times \Cal O\to \T^{\nu+b}\times\{0,0,0\}$ and $
\omega_*$ is close to the unperturbed frequency
$\omega$, such that for any $\sigma\in \Cal O_\gamma$ and $\theta\in
\T^{\nu+b}$, the curve $t\to \Psi(\theta+\omega_*(\sigma) t,\sigma)$ is a linearly stable
quasi-periodic solution of the Hamiltonian equations governed by
$H=N+P$.
\end{Theorem}

\begin{Remark} Compared to Montalto\cite{M}, our obtained solutions are at least $C^\infty$, while the obtained solutions in \cite{M} are $C^k$ ($k$ finite); Compared to Chen--Geng\cite{CG2}, we generalizes
the result of \cite{CG2} from the simple--eigenvalue case to the double--eigenvalue case.  \end{Remark}

\section{Application to the One Dimensional Forced Kirchhoff Equations under Periodic Boundary Conditions}

We consider  one dimensional Kirchhoff equations, by scaling $u\rightarrow \epsilon^{\frac{1}{3}}u$, we have
\begin{equation}\label{3.1}
u_{tt}-(1+\varepsilon\int_{0}^{2\pi} |u_x|^2 dx)u_{xx}+M_\xi u+\varepsilon g(\bar{\omega}t,x) =0,\quad  \varepsilon=\epsilon^{\frac{2}{3}}
 \end{equation}
with periodic boundary conditions $u(t,x+2\pi)=u(t,x)$.

 Here we assume that the operator
$A=-\partial_{xx}+M_\xi$ with periodic boundary conditions has
 eigenvalues $\{\mu_n\}$ satisfying
\begin{eqnarray*}
\tilde{\omega}_j&=&\lambda_j=\sqrt{\mu_{i_j}}=\sqrt{i_j^2+\xi_j},\  1\le j\le b,\ \
\Omega_n=\lambda_n=\sqrt{\mu_n}=|n|,\  n\in \Bbb Z\setminus\{i_1, \cdots, i_b\}
\end{eqnarray*}
 and the corresponding  eigenfunctions $\phi_n(x), n\in \Bbb Z$. We assume $\sigma=(\bar \omega,\xi_1,\cdots,\xi_b)$ is a parameter
taking on a closed set $\Cal O\subset\R^{\tilde b}$ of the
positive measure.

  Introducing $v=u_t$, $(\ref{3.1})$ reads
\beq\label{nonlinearKirchhoff1}\left\{\begin{array}{l}
u_t=v,\\
v_t=-Au+\varepsilon(\int_0^{2\pi}|u_x|^2 dx)u_{xx}-\varepsilon g(\bar \omega t,x),
\end{array}
\right.\eeq
the associated  Hamiltonian function
\[H=\frac{1}{2}\int_0^{2\pi}v^2dx+\frac{1}{2}(Au,u)+\varepsilon(\frac{1}{2}\int_0^{2\pi}|u_x|^2dx)^2+\varepsilon\int_0^{2\pi}g(\bar \omega t,x)udx,\]
where $(\cdot, \cdot)$ is the inner product in $L^2(\Bbb T)$. Then we introduce sequences $q=(q_n)_{n\in\Z},p=(p_n)_{n\in\Z}$,
\[u(x)=\sum_{n\in\Z}\frac{q_n}{\sqrt{\lambda_n}}\phi_n(x),\quad v(x)=\sum_{n\in\Z}\sqrt{\lambda_n}{p_n}\phi_n(x),\quad  \lambda_n=\sqrt{\mu_n},\]
this is equivalent to the lattice Hamiltonian equations
\beq\left\{\begin{array}{l}
\dot{q}_n=\lambda_n p_n,\\
\dot{p}_n=-\lambda_n q_n-\varepsilon\frac{\partial G}{\partial q_n},
\end{array}
\right.\eeq
and the corresponding Hamiltonian function
\[H(p,q)=\frac{1}{2}\sum_{n\in\Z}\lambda_n(p^2_n+q^2_n)+\varepsilon G(q),\]
\[G(q)=\frac{1}{4}\sum_{n,m\in\Z}\frac{n^2m^2}{\lambda_n \lambda_m}q^2_nq^2_m+\sum_{n\in\Z}g_n(\bar \omega t)\frac{q_n}{\sqrt{\lambda_n}}.\]
with the Fourier coefficients $\{g_n(\bar \omega t)\}$ of the function $g(x,\bar \omega t)$.

We switch to complex variables
\[w_n=\frac{q_n+{\rm i} p_n}{\sqrt{2}},\quad \bar{w}_n=\frac{q_n-{\rm i} p_n}{\sqrt{2}},\]
hence we obtain
 \[H=\sum_{n\in\Z}\lambda_n w_n\bar{w}_n +\varepsilon G(w,\bar{w}),\]
 \begin{eqnarray}\label{G}
 G(w,\bar
w)
&=&\frac{1}{4}\sum_{n,m\in\Z}\frac{n^2m^2}{\lambda_n\lambda_m}(\frac{w_n+\bar {w}_n}{\sqrt{2}})^2(\frac{w_m+\bar {w}_m}{\sqrt{2}})^2+\sum_{n\in\Z}g_n(\bar \omega t)\frac{w_n+\bar w_n} {\sqrt{2\lambda_n}}\nonumber\\
&=&G^1+G^2=\sum_{\alpha,\beta}G^1_{\alpha,\beta}w^{\alpha}\bar{w}^{\beta}
+\sum_{\alpha,\beta}G^2_{\alpha,\beta}w^{\alpha}\bar{w}^{\beta}
\end{eqnarray}
with $\alpha,\beta \in\{\alpha,\beta\in \N^{\Z},|\alpha+\beta|=4,\alpha_n+\beta_n \in 2\N,\forall n\in \Z\}$ in $G^1_{\alpha\beta}$, $\alpha,\beta \in\{\alpha,\beta\in \N^{\Z},|\alpha+\beta|=1\}$ in $G^2_{\alpha\beta}$.

 Moreover, the perturbation $G$ in $(\ref{G})$
has  the following  regularity property.
\begin{Lemma}\label{3.2}
For any fixed $a> 0,\rho >0$,the space $l_1^{a,\rho}$ is Banach algebra with respect to convolution of sequences, and
\[\|p\ast q\|_{a,\rho}\leq c\|p\|_{a,\rho}\|q\|_{a,\rho}.\]
\end{Lemma}
\proof  See \cite{P1}.\qed

\begin{Lemma}\label{Decay}
Suppose $g(x,\bar{\omega}t)$ is analytic with $|{\rm Im}x|<\bar \rho$, then
the coefficients $\{g_n\}_{n\in \Z}$ have the estimate $$\sup\limits_{t\in \R}|g_n(\bar \omega t)|\leq ce^{-|n|\bar\rho},\  \forall n\in \Z.$$ Thus one have
 $$|G^2_{\alpha\beta}|\leq c |n|^{-\frac{1}{2}}e^{-|n|\bar\rho},\quad \alpha+\beta=e_n.$$
\end{Lemma}
\proof   We expand $g(x,\bar{\omega}t)$  into Fourier series
\begin{equation}\label{F-S}
g(x,\bar{\omega}t)=\sum_{n\in \Z}\hat{g}_n(\bar{\omega}t)e^{{\rm i}nx},
\end{equation}

on the other hand, $g(x,\bar{\omega}t)=\sum\limits_{n\in \Z}g_n(\bar \omega t)\phi_n(x),$ then it is clear that $$g(x,\bar{\omega}t)=\sum_{n\in \Z}\hat{g}_n(\bar{\omega}t)e^{{\rm i}nx}=\sum_{n\in \Z}g_n(\bar \omega t)\phi_n(x).$$

 Since $g(x,\bar{\omega}t)$ is analytic in $x$, so $g(x,\bar{\omega}t)$ is bounded and $(\ref{F-S})$ is uniformly convergent
 $$\sup_{(x,t)\in [0,2\pi]\times\R}|g(x,\bar{\omega}t)|<c, \ \ \ \Rightarrow \ \ \ \sup_{t\in \R}\sum_{n\in \Z} |\hat{g}_{n}(\bar{\omega}t)|e^{|n|\bar\rho}<c, $$
 then we can obtain
 $$  \sup_{t\in \R}|\hat{g}_{n}(\bar{\omega}t)|e^{|n|\bar\rho}<c,\ \ \ \Rightarrow \ \ \sup_{t\in \R}|\hat{g}_{n}(\bar{\omega}t)|< ce^{-|n|\bar\rho},\ \ \forall n\in\Z,$$
and the coefficients ${g_n(\bar \omega t)}$ in the basis $\{\phi_n(x),n\in \Z\}$ satisfy
$$\sup_{t\in \R}|g_n(\bar \omega t)|\leq c\sup_{t\in \R}|\hat{g}_{n}(\bar{\omega}t)|< ce^{-|n|\bar\rho},\ \ \forall n\in \Z,$$ where c is some constant and may be different in the above formulas. Finally from (\ref{G}), we obtain the decay property  $|G^2_{\alpha\beta}|\leq c |n|^{-\frac{1}{2}}e^{-|n|\bar\rho}$ .\qed

\begin{Lemma}\label{regularity}
For any fixed $a\geq 1$, $0<\rho<\bar \rho$, the gradient $G_{\bar w}$ is real analytic as a map
 in a neighborhood of the origin with
 \beq \|G^1_{\bar w}\|_{\bar a,\rho}\le c\|w\|_{a,\rho}^3,\quad \|G^2_{\bar w}\|_{\bar a,\rho}\le c,\quad \bar a=a-1.
\label{3.3} \eeq
\end{Lemma}
\proof
 In $(\ref{G})$, we have
 \[G^1(w,\bar w)
=\frac{1}{16}\sum_{n,m\in \Z}\frac{n^2m^2}{\lambda_n\lambda_m}(w_n+\bar {w}_n)^2(w_m+\bar {w}_m)^2,\]
hence
\begin{eqnarray*}
G^1_{\bar w_n}=\frac{n}{8\sqrt{\lambda_n}}\sum_{m\in \Z}\frac{nm^2}{\sqrt{\lambda_n}\lambda_m}(w_n+\bar {w}_n)(w_m+\bar {w}_m)^2=\frac{n}{8\sqrt{\lambda_n}}h_n,
\end{eqnarray*}
where $h_n:=\sum_{m\in \Z}\frac{nm^2}{\sqrt{\lambda_n}\lambda_m}(w_n+\bar {w}_n)(w_m+\bar {w}_m)^2$ and defining $v=(v_n)_{n\in \Z}=((\tilde{w}*\tilde{w}*\tilde{w})_n)_{n\in \Z}$, $ \tilde{w}=(\sqrt{|n|}\cdot w_n)_{n\in \Z}$, we know
\begin{eqnarray*}
\|G^1_{\bar w}\|_{\bar a,\rho}&=&\sum_{n\in \Z}|\frac{n}{8\sqrt{\lambda_n}}h_n|\cdot n^{\bar a}e^{n\rho}
\leq c\sum_{n\in \Z}|h_n|\cdot |n|^{(\bar a+\frac{1}{2})}e^{n\rho}\nonumber\\
&\leq & c\sum_{n\in \Z}|v_n|\cdot |n|^{(\bar a+\frac{1}{2})}e^{n\rho}
\leq  c\|v\|_{\bar a+\frac{1}{2},\rho}\leq  c\|\tilde w\|^3_{\bar a+\frac{1}{2},\rho}\leq c\|w\|_{a,\rho}^3.
\end{eqnarray*}
By the above lemma $|G^2_{\alpha\beta}|\leq c |n|^{-\frac{1}{2}}e^{-|n|\bar\rho}$ and
 \begin{eqnarray}
\|G^2_{\bar w}\|_{\bar a,\rho}&=&\sum_{n\in \Z}\sum_{|\alpha|+|\beta-e_n|=0}|G^2_{\alpha\beta}||n|^{\bar a}e^{|n|\rho}
\leq c\sum_{n\in \Z}||n|^{-\frac{1}{2}}e^{-|n|\bar\rho}||n|^{\bar a}e^{|n|\rho}\nonumber\\
&\leq& c\sum_{n\in \Z}|n|^{\bar a-\frac{1}{2}}e^{-|n|(\bar\rho-\rho)}\leq  c,\nonumber
 \end{eqnarray}
 where we  let $\rho< \bar\rho,\bar a=a-1$, the sum will be bounded, so the lemma follows.\qed

 Next we first introduce auxiliary action-angle variables
 $(\bar{\theta},\bar{I})\in \T^{\nu}\times \R^{\nu} $ satisfying
\[\frac{d\bar \theta}{dt}=\frac{\partial H}{\partial \bar I}=\bar \omega ,\quad \frac {d\bar I}{dt}=-\frac{\partial H}{\partial \bar \theta},\quad
{\rm i}\frac {dw_n}{dt}=\frac{\partial H}{\partial \bar w_n},\ \ \
{\rm i}\frac {d\bar w_n}{dt}=-\frac{\partial H}{\partial w_n}, \
n\in\Z,\]
then we introduce the internal action-angle variables $(\tilde\theta,
I)=((\tilde \theta_1,\cdots, \tilde \theta_b),( I_1, \cdots,  I_b))\in \T^b\times \R^b $
in the
$(w_{i_1}, \cdots, w_{i_b},\bar w_{i_1}, \cdots, \bar w_{i_b}
)$-space by letting,
\[w_{i_j}= \sqrt{I_j}e^{{\rm i}\tilde{\theta}_j}, \bar {w}_{i_j}= \sqrt{I_j}e^{-{\rm i}\tilde{\theta}_j},\quad j=1, \cdots,b,\]
so the system becomes
\begin{eqnarray}
\frac{d\bar{\theta}_j}{dt}&=&\bar{\omega}_j,\quad \frac {d\bar{I}_j}{dt}=-
 P_{\bar{\theta}_j},\quad j=1, \cdots, \nu,\nonumber\\
\frac{d\tilde{\theta}_j}{dt}&=&\tilde{\omega}_j +  P_{I_j},\quad \frac {dI_j}{dt}=- P_{\tilde{\theta}_j},\quad j=1, \cdots, b,\label{3.4}\\
\frac {dw_n}{dt}&=&-{\rm i}(\Omega_nw_n+  P_{\bar w_n}),\ \ \
\frac {d\bar w_n}{dt}={\rm i}(\Omega_n\bar w_n+ P_{w_n}), \
n\in\Z_1,\nonumber
\end{eqnarray}
where $P$ is just $\varepsilon G$ with the $(w_{i_1}, \cdots, w_{i_b}, \bar
w_{i_1}, \cdots, \bar w_{i_b} )$-variables
expressed in terms of the $(\tilde\theta,I)$ variables.
The Hamiltonian associated to (\ref{3.4}) with respect to the
symplectic structure $\displaystyle d\bar I\wedge d\bar \theta+dI\wedge d\tilde\theta +{\rm
i}\sum_{n\in\Z_1} dw_n \wedge d\bar w_n$ is given by
\begin{equation}\label{Hamiltonian}
H=\la\bar \omega, \bar I\ra+\la\tilde \omega(\xi), I\ra+\sum_{{n\in \Z_1}} \Omega_n(\sigma)w_n\bar
w_n+ P(\theta, I, w, \bar w, \sigma,\varepsilon),
\end{equation}
where $\omega=(\bar{\omega},\tilde{\omega}),\tilde{\omega}(\xi)=(\sqrt{i_1^2+\xi_1},\cdots,
\sqrt{i_b^2+\xi_b}),\Omega_n=|n|,n\in \Z_1.$

Next let us verify  that $H=N+P$ satisfies the assumptions $(A1)-(A7)$ in the initial step.

Verification of $(A1)$:
\begin{eqnarray*}
|\omega|_{\Cal O}&=&\sup_{\sigma\in \Cal O}\sup_{1\leq j\leq \tilde b}\{|\omega_j|+|\frac{\partial \omega_j}{\partial \sigma}|\}\nonumber\\
&=&\max\{|\bar \omega|+1,\sup_{\sigma\in \Cal O}\sup_{1\leq j\leq b}\{|\sqrt{i_j^2+\xi_j}|+|\frac{1}{2\sqrt{i_j^2+\xi_j}}|\}\}\leq E_0,\nonumber
\end{eqnarray*}
\[\frac{\partial \omega}{\partial\sigma}=\left(
\begin{array}{ccc}
I_{\nu\times \nu}&0 \\
0&Diag(\frac{1}{2\sqrt{|i_1|^2+\xi_1}},\ldots,\frac{1}{2\sqrt{|i_b|^2+\xi_b}})\\
\end{array}
\right),det(\frac{\partial \omega}{\partial\sigma})\neq0,\]
so $\omega(\sigma)$ is a $C^1_W$ diffeomorphism and there exists a positive constant $E_0$ such that $|\omega|_{\Cal O}\leq E_0$.

Verification of $(A2)$: According to the form of $N$ and $G$ in the initial step, it is obviously that
$$\bar{\Omega}_n=|n|,\quad \tilde{\Omega}_n =0, A_{|n|}=0,\quad c(\sigma)=f(\theta,\sigma)=d_n=d_{(-n)}=0.$$ Then $(A2)$ is automatically satisfied.

Verification of $(A3)$: In the initial step in Section 4, the small divisors have three kinds of form
$|\langle k,\omega\rangle|\leq \frac{\gamma}{|k|^\tau},\ 0\neq k=(k_1,k_2)\in \Z^{\nu+b}$;\
$|\langle k,\omega\rangle\pm\bar{\Omega}_n|\leq \frac{\gamma}{K_0^\tau},\ |k|\leq K_0, |n|\leq E_0K_0$;\
$|\langle k,\omega\rangle\pm 2\bar{\Omega}_n|\leq \frac{\gamma\cdot |n|}{K_0^\tau},\ |k|\leq K_0$.\\
For the first one $|k|\neq 0$, we have
\[\mid\frac{\partial\langle k,\omega\rangle}{\partial\sigma}\mid=|(k_{1_1},\cdots,k_{1_\nu},\frac{k_{2_1}}{2\sqrt{|i_1|^2+\xi_1}},\cdots,
\frac{k_{2_b}}{2\sqrt{|i_b|^2+\xi_b}})|\geq c|k|,\]
so for any fixed $k$,
\[\textmd{ meas}\{\sigma:|\langle k,\omega\rangle
|<\frac{\gamma}{|k|^\tau}\}\leq c\frac{\gamma}{|k|^{\tau+1}}.\]
then if $\tau \geq\tilde b$,
\[\sum_{|k|\neq 0} \textmd{meas}\{\sigma:|\langle k,\omega\rangle|<\frac{\gamma}{|k|^\tau}\}\leq c\sum_{|k|\neq 0}\frac{\gamma}{|k|^{\tau+1}}<c \gamma.\]
For the second, if $|\bar{\Omega}_n|\geq c|k|+1$, then
$|\langle k,\omega\rangle\pm \bar{\Omega}_n|\geq |\bar{\Omega}_n|-c|k|\geq 1$, there will be no small divisors. Otherwise,
if $1\leq |\bar{\Omega}_n|\leq c|k|+1$, $0\neq|k|\leq K_0$, then
\[\mid\frac{\partial(\langle k,\omega\rangle\pm \bar{\Omega}_n)}{\partial\sigma}\mid\geq c|k|\geq c.\]
For fixed $|k|\leq K_0,|n|\leq E_0K_0,$
\[\textmd{ meas}\{\sigma:|\langle k,\omega\rangle\pm \Omega_n|<\frac{\gamma}{K_0^\tau}\}\leq c\frac{\gamma}{K_0^{\tau}},\]
similarly for the last one, by the same argument, we have for fixed $|k|\leq K_0,|n|\leq C|k|+1,$
\[\textmd{ meas}\{\sigma:|\langle k,\omega\rangle\pm 2\Omega_n|<\frac{\gamma\cdot |n|}{K_0^\tau}\}\leq c\frac{\gamma\cdot |n|}{K_0^{\tau}},\]
then if $\tau >\tilde b+2$,
\[\sum_{0<|k|\leq K_0,n} \textmd{meas}\{\sigma:|\langle k,\omega\rangle\pm\Omega_n|<\frac{\gamma}{|k|^\tau}\}\leq\sum_{0<|k|\leq K_0,\atop |\Omega_n|\leq c|k|+1}cK_0^{\tilde b+1}\frac{\gamma}{K_0^{\tau}}\leq c\frac{\gamma}{K_0^{\tau-\tilde b-1}}<c \gamma,\]
\[\sum_{0<|k|\leq K_0,n} \textmd{meas}\{\sigma:|\langle k,\omega\rangle\pm 2\Omega_n|<\frac{\gamma\cdot |n|}{|k|^\tau}\}\leq\sum_{0<|k|\leq K_0,\atop |2\Omega_n|\leq c|k|+1}cK_0^{\tilde b+1}\frac{\gamma\cdot |n|}{K_0^{\tau}}\leq c\frac{\gamma}{K_0^{\tau-\tilde b-2}}<c \gamma.\]
so there exists a subset $\cal{O}_\gamma \subset \cal{O}$ with meas$(\cal{O}\backslash\cal{O}_\gamma)=O(\gamma)$ such that for any $\sigma\in \cal{O}_\gamma$, the non-resonance conditions in the initial step are satisfied.

Verification of $(A4)$: In fact, the regularity of $P$ holds true:
\begin{Lemma}\label{regularityP}  For any
  $\varepsilon>0$ sufficiently small and $s\ll1$, if $|I|<s^2$ and $\|w\|_{a,\rho}<s$,
 then we have
 \beq \|X_P\|_{s,\bar a, \rho, D(r,s),\Cal O}\le \varepsilon,\qquad \bar a=a-1.
\label{3.16P} \eeq
\end{Lemma}
\proof
According to Lemma 3.2,
$$\varepsilon\|G^1_{\bar w}\|_{\bar{a},\rho}\le c\varepsilon\|w\|_{a,\rho}^3,\  \varepsilon\|G^2_{\bar w}\|_{\bar{a},\rho}\le c\varepsilon.$$
Denote $P^1+P^3, P^2$ instead of $\varepsilon G^1, \varepsilon G^2$ respectively after the transformation of the action-angle variables, then we have
\begin{eqnarray*}
\sum_{ n\in\Z_1} \|P^1_{w_n}\|_{D(r,s),\cal{O}}|n|^{\bar{a}}e^{|n|\rho}+\|P^1_{\bar{w}_n}\|_{D(r,s),\cal{O}}|n|^{\bar{a}}e^{|n|\rho}&=&
\|P^1_w\|_{\bar{a},\rho}+\|P^1_{\bar{w}}\|_{\bar{a},\rho}\\
&\leq& c\varepsilon\|w\|^3_{a,\rho}\leq c \varepsilon (|I|^\frac{3}{2}+\|w\|^3_{a,\rho}).
\end{eqnarray*}
It is obvious that $\sup\limits_{\|w\|_{a,\rho}\leq 2s\atop\|\bar w\|_{a,\rho}\leq 2s}\|G^1\|_{D(r),\cal{O}}\leq cs^4$,
thus
$\|P^1\|_{\!{}_{D(2r,2s),\cal{O}}}\leq c\varepsilon s^4.$ According to Cauchy estimates,
$\|P^1_I\|_{\!{}_{D(r,s),\cal{O}}}\leq c\varepsilon s^2,\|P^1_\theta\|_{\!{}_{D(r,s),\cal{O}}}\leq c \varepsilon s^4,$
hence
\begin{eqnarray}
\|X_{P^1}\|_{\!{}_{s,\bar a,\rho, D(r,s),\cal O}}&=& \|P^1_I\|_{\!{}_{ D(r,s)
, \cal O}}+ \frac
1{s^2}\|P^1_\theta\|_{\!{}_{ D(r,s),\cal O}}\nonumber\\
&+& \frac 1s\sum_{n\in\Z_1} \|P^1_{w_n}\|_{\!{}_{ D(r,s),\cal
O}}|n|^{\bar a}e^{|n|\rho}+\frac 1s\sum_{n\in\Z_1} \|P^1_{\bar
w_n}\|_{\!{}_{ D(r,s),\cal O}}|n|^{\bar a}e^{|n|\rho}
\nonumber\\
&\leq&c \varepsilon s^2+c\varepsilon s^2+c \varepsilon\frac{1}{s}(|I|^\frac{3}{2}+\|w\|^3_{a,\rho})
\leq c\varepsilon.\nonumber
\end{eqnarray}
With the similar arguments, we have
$$\|P^2_w\|_{\bar{a},\rho}+\|P^2_{\bar{w}}\|_{\bar{a},\rho}\le c\varepsilon ,\sup\limits_{\|w\|_{a,\rho}\leq 2s\atop\|\bar w\|_{a,\rho}\leq 2s}\|G^2\|_{D(r),\cal{O}}\leq cs,$$ then
$\|P^2_\theta\|_{\!{}_{D(r,s),\cal{O}}}\leq c \varepsilon s$ and
\begin{eqnarray*}
\|X_{P^2}\|_{\!{}_{ s,\bar a,\rho, D(r,s), \cal O}}&=& \frac{1}{s^2}\|P^2_\theta\|_{\!{}_{ D(r,s)
, \cal O}}+\frac 1s\sum_{n\in\Z_1} \|P^2_{w_n}\|_{\!{}_{ D(r,s) , \cal
O}}|n|^{\bar a}e^{|n|\rho}\\
&+& \frac 1s\sum_{n\in\Z_1}\|P^2_{\bar
w_n}\|_{\!{}_{ D(r,s) , \cal O}}|n|^{\bar a}e^{|n|\rho}
\leq c\varepsilon,\nonumber\end{eqnarray*}
it is easy to prove $\|X_{P^3}\|_{\!{}_{ s,\bar a,\rho, D(r,s), \cal O}}\leq c\varepsilon$, hence
$$\|X_P\|_{\!{}_{s,\bar a,\rho, D(r,s),\Cal O}}\le \|X_{P^1}\|_{\!{}_{s,\bar a,\rho, D(r,s),\Cal O}}+\|X_{P^2}\|_{\!{}_{s,\bar a,\rho, D(r,s),\Cal O}}+\|X_{P^3}\|_{\!{}_{s,\bar a,\rho, D(r,s),\Cal O}}\leq \varepsilon.$$
  The verification of $(A4)$ is accomplished.\qed

Verification of $(A5)$: Observing the form of the perturbation $G$ in $(\ref{G})$,  $P$ can also be written as follows $P=P^1+P^2+P^3$,
\begin{eqnarray}
P^1&=&\frac{\varepsilon}{8}\sum_{1\leq j\leq b\atop n\in \Z_1}\frac{i_j^2|n|}{\sqrt{i_j^2+\xi_j}}I_j(e^{2{\rm i}\tilde\theta_j}+2+e^{-2{\rm i}\tilde\theta_j})(w_n^2+2w_n\bar w_n+\bar w_n^2)\nonumber\\
&+&\frac{\varepsilon}{16}\sum_{m, n\in\Z_1}|n||m|(w_n^2+2w_n\bar w_n+\bar w_n^2)(w_m^2+2w_m\bar w_m+\bar w_m^2),\label{p1}\\
P^2&=&\varepsilon\sum_{n\in\Z_1}g_n(\bar \omega t)\frac{w_n+\bar w_n}{\sqrt{2|n|}},\label{p2}\\
P^3&=&\frac{\varepsilon}{16}\sum_{1\leq j,k\leq b}\frac{i_j^2 i_k^2}{\sqrt{(i_j^2+\xi_j)(i_k^2+\xi_k)}}I_jI_k(e^{2{\rm i}\tilde\theta_j}+2+e^{-2{\rm i}\tilde\theta_j})(e^{2{\rm i}\tilde\theta_k}+2+e^{-2{\rm i}\tilde\theta_k})\nonumber\\
&+&\varepsilon\sum_{1\leq j\leq b}g_j(\bar \omega t)\frac{e^{{\rm i}\tilde\theta_j}+e^{-{\rm i}\tilde\theta_j}}{\sqrt{2(i_j^2+\xi_j)}}\sqrt{I_j},\label{p3}\end{eqnarray}
the exponents of $w, \bar w$ in $P^1,P^2,P^3$ respectively satisfies the assumption $(A5)$ and by Lemma $\ref{Decay}$, the decay properties of $P^2$ can be satisfied automatically.\\

Verification of $(A6)$:
 According to the perturbation $P$ in the assumption $(A5)$, we just need to consider the T\"oplitz-Lipschitz property of the first term $P^1$,
 when $n\in \Z_1$, the second order derivatives of $P^1=P^1(\theta,I,w,\bar w, \varepsilon)$
\begin{eqnarray*}
\frac{1}{|n|}\sum_{\upsilon=\pm}\frac{\partial^2P^1}{\partial w^\upsilon_n\partial w^\upsilon_n}
&=&\frac{3\varepsilon}{4}\sum_{1\leq j\leq b}\frac{i_j^2}{\sqrt{i_j^2+\xi_j}}I_j(e^{2{\rm i}\tilde\theta_j}+2+e^{-2{\rm i}\tilde\theta_j})+\frac{3\varepsilon}{4}|n|(w_n^2+2w_n\bar w_n\\
&+&\bar w_n^2)
+\frac{3\varepsilon}{4}\sum_{m\neq n\in\Z_1}|m|(w_m^2+2w_m\bar w_m+\bar w_m^2),\nonumber
\end{eqnarray*}
it is obvious that the first and third sum are independent of $n$ and  uniformly convergent in the form $\|\cdot\|_{D(r,s),\Cal O}$ due to the set of indices $1\leq j\leq b$ is finite and $\|w\|_{a,\rho}\leq s,\|\bar w\|_{a,\rho}\leq s$, we deduce $|w_n|,|\bar w_n|\leq s|n|^{-a}e^{-|n|\rho},a\geq 1$, then
\begin{eqnarray*}
&&\|\lim_{n\rightarrow\infty}\frac{1}{|n|}\sum_{\upsilon=\pm}\frac{\partial^2P^1}{\partial w^\upsilon_n\partial w^\upsilon_n}\|_{\!{}_{D(r,s),\Cal O}}=\frac{3\varepsilon}{4}\|\lim_{n\rightarrow\infty}|n|(w_n^2+2w_n\bar w_n+\bar w_n^2)\|_{\!{}_{D(r,s),\Cal O}}\leq \varepsilon,\\
&&\|\frac{1}{|n|}\sum_{\upsilon=\pm}\frac{\partial^2P^1}{\partial w^\upsilon_n\partial w^\upsilon_n}- \lim_{n\rightarrow\infty}\frac{1}{|n|}\sum_{\upsilon=\pm}\frac{\partial^2P^1}{\partial w^\upsilon_n\partial w^\upsilon_n}\|_{\!{}_{D(r,s),\Cal O}}
\leq 3\varepsilon(s|n|^{-a+1}e^{-|n|\rho})^2\leq \frac{\varepsilon}{|n|}.
\end{eqnarray*}

Verification of $(A7)$:
 In the assumption $(A2)$, we know $\tilde\Omega_n=0$, so $\|\tilde \Omega_n\|_{r,2\tau+2,\Cal O}=0$, $(A7)$ is satisfied.

To this point, we have verified  all the initial assumptions of Theorem
\ref{KAM}. By applying Theorem
\ref{KAM}, we get Theorem \ref{maintheorem}. In the next sections, we will show explicitly how to construct an iterative KAM algorithm to prove Theorem \ref{KAM}.

\section{ KAM Step}

  A KAM iteration involves an infinite sequence of transformation and each step makes the perturbation smaller than that of the previous
one at the cost of excluding a small set of parameters. We have
to prove  the convergence of the iteration and estimate the
measure of the excluded set after infinite KAM steps.\\
 In our paper, due to the special structure and the decay property of the perturbation, it is necessary to show the initial KAM step clearly to see how those coupled terms appear and the coefficients of them inherit some decay property from the perturbation $g(\bar\omega t,x)$. Thanks to these special properties, it is feasible for us to implement KAM iteration and prove the convergence of the iteration and  measure estimate.
\subsection{Normal form}\label{4.1}

 In order to perform the KAM iteration, we will first write the Hamiltonian into a normal form and  fix the positive constant $\bar \rho>0$ in the whole KAM iteration. Denote $E_{-1}=K_{-1}=1$. Choosing $\varepsilon_0=\varepsilon$, $\varepsilon_1\sim \varepsilon_0^{\frac{4}{3}}$, $K_0\sim|\ln\varepsilon_0|$, $r_0=r,E_0=E,s=s_0,\rho=\rho_0=\frac{\bar \rho}{2}$. Let $\rho_1<\rho_0<\bar\rho$ and $s_0$ be such that $0< s_1<\min\{\varepsilon_1,s_0\}.$\\
 Recalling that $H$ in (\ref{Hamiltonian}), $$H=N+P =\la\bar{\omega},\bar{I}\ra+\la\tilde{\omega}(\xi),I\ra+\sum_{n\in\Z_1}\Omega_n(\sigma)w_n \bar
w_n+ P(\theta,I,w,\bar w,\sigma,\varepsilon),$$
 where $P=P^1+P^2+P^3$ with $P^1, P^2, P^3$ in $(\ref{p1}),(\ref{p2}),(\ref{p3})$ satisfies the assumption $(A5)$.\\
Let the truncation $R$ be as follows
 \begin{eqnarray}
 R&=&\sum_{|l|\leq 1}P_{l00}(\theta,\sigma)I^l+\sum_{n\in \Z_1\atop |n|\leq E_0K_0}(P_{n}^{10}(\theta,\sigma)w_n+P_{n}^{01}(\theta,\sigma)\bar w_n)\nonumber\\
 &+&\sum_{n\in \Z_1}(P^{20}_{nn}(\theta,\sigma)w_nw_n+P^{11}_{nn}(\theta,\sigma)w_n\bar w_n+P^{02}_{nn}(\theta,\sigma)\bar w_n\bar w_n)\nonumber\\
 &=&R_0+R_1+R_2.\label{R0}
\end{eqnarray}
To handle the term $R$, we will first construct a symplectic transformation $\Phi_0=\phi_{F_0}^1$,
\begin{eqnarray}
\Gamma_{K_0}F_0=F_0&=&\sum_{|l|\leq 1}F_{l00}(\theta,\sigma)I^l +\sum_{n\in \Z_1\atop |n|\leq E_0K_0}(F_{n}^{10}(\theta,\sigma)w_n+F_{n}^{01}(\theta,\sigma)\bar w_n) \nonumber\\
 &+&\sum_{n\in \Z_1}(F^{20}_{nn}(\theta,\sigma)w_nw_n+F^{02}_{nn}(\theta,\sigma)\bar w_n\bar w_n)\nonumber\\
 &=&F^{0}+F^{1}+F^{2},\label{F0}
\end{eqnarray}
where $[F_{l00}]=0$, so the terms $[P_{l00}](|l|\leq 1),P_{nn}^{11}(\theta,\sigma)w_n\bar w_n$  will be added to the normal form part of the new Hamiltonian. More precisely, let $F_0$ satisfy the homological equation
\[\{N,F_0\}+R=\sum_{|l|\leq 1}[P_{l00}]I^l+\sum_{n\in \Z_1}P_{nn}^{11}(\theta,\sigma)w_n\bar w_n,\]
where $N=\la\bar{\omega},\bar{I}\ra+\la\tilde{\omega}(\xi),I\ra+\sum\limits_{n\in\Z_1}\Omega_n(\sigma)w_n \bar
w_n$. Moreover,it is clear that
\begin{eqnarray}
 P-R&=&(1-\Gamma_{K_0})(\sum_{|l|\leq 1}P_{l00}(\theta,\sigma)I^l+\sum_{|n|\geq E_0K_0\atop\alpha_n+\beta_n=1}P_n(\theta,\sigma)w_n^{\alpha_n}\bar w_n^{\beta_n}\nonumber\\
 &+&\sum_{|\alpha+\beta|=2\atop \alpha_n+\beta_n\in 2\N}P_{0\alpha\beta}(\theta,\sigma)w^\alpha \bar w^\beta)
 +O(|I|^2+|I||w|^2+|w|^4),\nonumber
\end{eqnarray}
where by Lemma $\ref{Decay}$ and $(\ref{decayP^2})$ in the assumption $(A5)$
 \begin{eqnarray}\label{decay}
\|P_n(\theta,\sigma)\|_{D(r),\Cal O}\leq c \varepsilon e^{-|n|\bar{\rho}}, \ |n|> E_0K_0.
 \end{eqnarray}
It thus follows from (\ref{decay}) and the Cauchy inequality, one can make $\rho_1<\rho_0,s_1\ll s_0$ small enough such that $\|X_{P-R}\|_{s_0,\bar a, \rho_1,D(r,s_1),\Cal O}\leq \varepsilon_1$.\\
In section 3, we have proved that this homological equation is solvable with $|k|\leq K_0$ on the parameter set with meas$(\Cal O_0\backslash\Cal O_1)\leq c\gamma$:
$$\Cal O_1=
\left\{\sigma\in \Cal O_0:
\begin{array}{rcl}
&&|\la k,\omega\ra|\geq \frac{\gamma}{|k|^\tau},0\neq k=(k_1,k_2)\in \Z^{\nu+b}\\
&&|\la k,\omega\ra\pm \Omega_n|\geq \frac{\gamma }{K_0^\tau},k=(k_1,k_2)\in \Z^{\nu+b},n\in \Z_1,|n|\leq E_0K_0\\
&&|\la k,\omega\ra\pm 2\Omega_n|\geq \frac{\gamma\cdot |n|}{K_0^\tau},k=(k_1,k_2)\in \Z^{\nu+b},n\in \Z_1
\end{array}
\right\}.$$
In this way, we obtain the transformation $\Phi_0$ which transforms the Hamiltonian to
\[H_1=H\circ\Phi_0=N_1+P_1,\]
where
\begin{eqnarray}
N_1&=&\la\bar{\omega},\bar{I}\ra+\la\tilde{\omega}_1(\sigma),I\ra+\sum_{n\in\Z_1 }\Omega_n^1(\theta,\sigma)w_n \bar w_n,\nonumber\\
&&\tilde{\omega}_1(\sigma)=\tilde{\omega}(\sigma)+[P_{l00}],(|l|=1),\quad  \Omega_n^1(\theta,\sigma)=\Omega_n(\sigma)+P^{11}_{nn}(\theta,\sigma),\nonumber\\
P_1&=&\sum_{\alpha,\beta}P^1_{\alpha\beta}(\theta,I,\sigma)w^\alpha\bar w^\beta+\sum_{\alpha,\beta}P^2_{\alpha,\beta}(\theta,\sigma)w^{\alpha}\bar w^{\beta}+\sum_{\alpha,\beta}P^3_{\alpha,\beta}(\theta,I,\sigma)w^{\alpha}\bar w^{\beta}\nonumber\\
&=&P_1^1+P_1^2+P_1^3,\label{P0}
\end{eqnarray}
with $l\in \N^b$, $\alpha,\beta \in \{\alpha,\beta\in \N^{\Z_1},\sum\limits_{|n|>E_0K_0}\alpha_n+\beta_n>0, \alpha_n+\beta_n\in 2\N, \forall |n|>E_0K_0\}$ in $P^1_1$, $\alpha,\beta \in \{\alpha,\beta\in \N^{\Z_1}, |\alpha+\beta|=\alpha_n+\beta_n=1,\forall |n|>E_0K_0\}$ in $P^2_1$, $\alpha,\beta \in \{\alpha,\beta\in \N^{\Z_1}, \alpha_n+\beta_n=0, \forall |n|>E_0K_0\}$ in $P^3_1$. Due to the special structure of $P_1=P_1^1+P_1^2+P^3_1$,  by Lemma $\ref{Decay}$ and Lemma $
\ref{Lem7.3}$,
if $|n|,|m|\leq E_0K_0, \alpha+\beta=e_n+e_m$ in $P_1^3$, we have
$$\|P^3_{\alpha\beta}\|_{D(r),\Cal O}\leq c\varepsilon |n||m| e^{-(|n|+|m|)\bar\rho}\leq c\varepsilon e^{-(|n|+|m|)\rho_1};$$
if $|n|>E_0K_0,\alpha+\beta=e_n$ in $P_1^2$, we have
$$\|P^2_{\alpha\beta}\|_{D(r),\Cal O}\leq c \varepsilon e^{-|n|\bar{\rho}}.$$

 Indeed, the terms $w_nw_m$, $w_n\bar w_m$,$\bar w_n\bar w_m$ consequently appear
 due to the Poisson bracket of $ P-R$ and $F^1$, with $R$ in (4.1) and $F^1$ defined in (4.2). Specifically,
 among the $\{ P-R, F^1\}$, there are some terms like
\begin{eqnarray}
\{ (w_n+\bar {w}_n)^2(w_m+\bar {w}_m)^2 , w_n+\bar w_n\}=4(w_n+\bar
w_n)(w_m+\bar w_m)^2, |n| \leq E_0K_0,m\in \Z_1 \nonumber
\end{eqnarray}
included in $P^1_1$ defined in (4.4), with  $\alpha,\beta \in
\{\alpha,\beta\in
\N^{\Z_1},\sum\limits_{|n|>E_0K_0}\alpha_n+\beta_n>0,
\alpha_n+\beta_n\in 2\N, \forall |n|>E_0K_0\}$. This means that the``
pair-property" will be preserved in $P^1_1$ only if the
spatial-index $n$ is large.

 Furthermore, among the $\{\{ P-R, F^1\}, F^1\}$, there are some terms like
\begin{eqnarray}
\{ \{(w_n+\bar {w}_n)^2(w_m+\bar {w}_m)^2 , w_n+\bar w_n\},w_m+\bar w_m\}=
8(w_n+\bar w_n)(w_m+\bar w_m), |n|, |m| \leq E_0K_0\nonumber
\end{eqnarray}
included in $P^3_1$ defined in (4.4), with $\alpha,\beta \in
\{\alpha,\beta\in \N^{\Z_1}, \alpha_n+\beta_n=0, \forall
|n|>E_0K_0\}$. This means the variables $w_n, \bar w_n $ contained in
$P^3_1$ exist only if all the spatial-indices $n$ are less than
$E_0K_0$. The terms $w_nw_m$, $w_n\bar w_m$,$\bar w_n\bar w_m$ need
to be added to the truncation $R$ in the next step. Their coefficients
 are all bounded due to the exponential decay property of the
coefficients of order 1 in $w,\bar w$. This means the ``pair
property " is totally destroyed in $P^3_1$.

Besides, the term $P^2_1$ only contains the first-order terms, like
$w_n,\bar w_n, |n|>E_0K_0 $ coming from $ P-R$. So in conclusion, the
perturbation $P_1$ also has the special form which is the assumption
$(A5)$.

So at the $\nu$-step of the KAM iteration,  we consider a
Hamiltonian vector field
$$ H_\nu=N_\nu+ P_\nu,\ \ \nu\geq 1,$$
where $N_\nu$ is a ``generalized normal form" and $P_{\nu}$ is defined in $D(r_\nu, s_\nu)\times \Cal O_{\nu}$.

 We then construct a map
$$\Phi_\nu: D (r_{\nu+1}, s_{\nu+1})\times\Cal O_{\nu+1}
 \to D(r_\nu, s_\nu)
\times\Cal O_{\nu},
$$
so that the vector field $X_{H_\nu \circ\Phi_\nu}$ defined on
$D(r_{\nu+1}, s_{\nu+1})$ satisfies
\[\|X_{P_{\nu+1}}\|_{s_{\nu+1},\bar a,\rho_{\nu+1},D(r_{\nu+1}, s_{\nu+1}),\Cal
O_{\nu+1}}=\|X_{H_\nu\circ\Phi_\nu}-X_{N_{\nu+1}}\|_{s_{\nu+1},\bar a,\rho_{\nu+1},D(r_{\nu+1}, s_{\nu+1})\times\Cal
O_{\nu+1}}\le \varepsilon_{\nu}^\kappa,\] $\kappa>1$, with some new normal form $N_{\nu +1}$.

\sss
To simplify notations, in what follows, the quantities
without subscripts refer to quantities at the  $\nu^{\rm th}$
step, while the quantities with subscripts $-$, $+$ respectively denote the
corresponding quantities at the $(\nu-1)^{\rm th}$, $(\nu+1)^{\rm th}$  step. Let us
then consider Hamiltonian function
\begin{eqnarray*} \label{h}
H&=&N+ P\\
&\equiv&\la\bar{\omega},\bar{I}\ra+\la\tilde{\omega}(\sigma),I\ra
 + \sum_{n\in \Z_1}\Omega_{n}(\theta,\sigma)w_n \bar w_n
+\sum_{|n|\leq E_{-}K_{-}}\langle A_{|n|}z_{|n|},\bar z_{|n|}\rangle +P(\theta, I, w,\bar w,
\sigma,\varepsilon)\\
&\equiv&\la\bar{\omega},\bar{I}\ra+\la\tilde{\omega}(\sigma),I\ra
 +\sum_{|n|\leq E_{-}K_{-}}\langle [\Omega_{n}(\theta,\sigma)I_2+A_{|n|}]z_{|n|},\bar z_{|n|}\rangle\\
 &+& \sum_{|n|>E_{-}K_{-}}\Omega_{n}(\theta,\sigma)w_n \bar w_n
 +P(\theta, I, w,\bar w,
\sigma,\varepsilon)\end{eqnarray*}
 defined in $D(r, s)\times\Cal O$ with
$\|X_P\|_{s,\bar a,\rho,D( r,s),\Cal O}\leq \varepsilon. \label{4.4}$
Because $A_{|n|}$ is real symmetric matrix, there exists an orthogonal matrix $Q_{|n|}$ such that
$$Q_{|n|}^TA_{|n|}Q_{|n|}=\Lambda_{|n|}=\left(\begin{array}{ccc}d_n(\sigma)&0\\
0&d_{(-n)}(\sigma)
\end{array}\right), Q_{|n|}^TI_2Q_{|n|}=I_2.$$
We still denote the variables with $|n|\leq E_{-}K_{-}$ by $w_n, w_{(-n)}$ without confusion. Hence our Hamiltonian function
\begin{eqnarray*}H&=&N+ P\\
&\equiv&
\la\bar{\omega},\bar{I}\ra+\la\tilde{\omega}(\sigma),I\ra
 +\sum_{|n|\leq E_{-}K_{-}}(\Omega_{n}(\theta,\sigma)+d_n)w_n\bar w_n\\
 &+& \sum_{|n|>E_{-}K_{-}}\Omega_{n}(\theta,\sigma)w_n \bar w_n
 +P(\theta, I, w,\bar w,
\sigma,\varepsilon).\end{eqnarray*}
\begin{remark}\label{RK}
Note that we introduce the orthogonal matrices in order to only simplify the notations for solving the homological equations, alternatively, we should solve vector or matrix homological equations. However, the essential small divisor difficulties are the same, hence, we intend to solve the scalar homological equations. In fact, we return to the original coordinates for $N_+$ and $P_+$.
 \end{remark}
 Next we will describe how  to construct a set $\Cal O_+\subset \Cal O$
and a change of variables
 $\Phi: D_+\times\Cal O_+=D(r_+, s_+)\times \Cal O_+\to  D(r,s)  \times \Cal O$
 such that the transformed
Hamiltonian $H_+=N_++P_+\equiv H\circ \Phi$ satisfies all the
above iterative assumptions with new parameters $s_+,r_+,\rho_+,
\varepsilon_+,$ and with $\sigma\in \Cal O_+$.

\subsection{Solving the Homological Equations}\label{4.2}

According to $(A5)$, expanding $P=P^1+P^2+P^3$ into the  Taylor series
\begin{eqnarray*}
P^1&=&\sum_{\alpha,\beta} P^1_{\alpha\beta}(\theta,I,\sigma)w^\alpha\bar w^\beta=\sum_{l,\alpha,\beta} P^1_{l\alpha\beta}(\theta,\sigma)I^lw^\alpha\bar w^\beta,\\
P^2&=&\sum_{\alpha,\beta}P^2_{\alpha\beta}(\theta,\sigma)w^{\alpha}\bar w^{\beta},\nonumber\\
P^3&=&\sum_{\alpha,\beta}P^3_{\alpha\beta}(\theta,I,\sigma)w^{\alpha} \bar w^{\beta}=\sum_{l,\alpha,\beta} P^3_{l\alpha\beta}(\theta,\sigma)I^lw^\alpha\bar w^\beta,
\end{eqnarray*}
with $l\in \N^b$, $\alpha,\beta \in\{\alpha,\beta\in\N^{\Z_1}, \sum\limits_{|n|>EK}\alpha_n+\beta_n>0, \alpha_n+\beta_n\in 2\N, \forall |n|>E_-K_- \}$ in $P^1$, $\alpha,\beta \in\{\alpha,\beta\in\N^{\Z_1}, |\alpha+\beta|= \alpha_n+\beta_n=1, \forall |n|>E_-K_- \}$ in $P^2$, $\alpha,\beta \in\{\alpha,\beta\in\N^{\Z_1}, \alpha_n+\beta_n=0, \forall |n|>E_-K_- \}$ in $P^3$. In addition, by the assumption (A5) and Lemma 7.3,
when $|n|\leq E_-K_-,\alpha+\beta=e_n$ or $|n|, |m|\leq E_-K_-, \alpha+\beta=e_n+e_m$ in $P^3$, we have
\begin{eqnarray}{\label{decay11}}
&&\|P^3_{l\alpha\beta}\|_{D(r),\Cal O}\leq c\varepsilon  e^{-|n|\rho},\ \ \|P^3_{l\alpha\beta}\|_{D(r),\Cal O}\leq c\varepsilon e^{-(|n|+|m|)\rho},
\end{eqnarray}
when $|n|>E_-K_-,\alpha+\beta=e_n$ in $P^2$, we have
\begin{eqnarray}{\label{decay1}}
\|P^3_{\alpha\beta}\|_{D(r),\Cal O}\leq c \varepsilon e^{-|n|\bar{\rho}}.
\end{eqnarray}
\begin{remark}Compared to Chen-Geng \cite{CG2}, the homological equations for solving $P^1+P^2$ are the same, the difference is to the homological equations for solving $P^3$. In \cite{CG2}, the term $(P^3)_{nn}^{11}(\theta,\sigma)w_n\bar w_n$ is put into the generalized normal form $N$, while in this paper, we only put the term $[(P^3)_{nn}^{11}(\theta,\sigma)]w_n\bar w_n+[(P^3)_{n(-n)}^{11}(\theta,\sigma)]w_n\bar w_{(-n)}+[(P^3)_{(-n)n}^{11}(\theta,\sigma)]w_{(-n)}\bar w_n$ into the generalized normal form $N$. Hence our normal frequencies $\Omega_n$, $A_{|n|}$ satisfy assumption $(A2)$. \end{remark}

Let $R$ be the truncation of $P$ given by
\begin{eqnarray*}
R&=& R_0+R_1+R_2,\ \ R_2=R_{2,<}+R_{2,>},\\
R_0&=&\sum_{|l|\leq 1}P_{l00}(\theta,\sigma)I^l,\ \
R_1=\sum_{n\in \Z_1\atop |n|\leq EK}(P^{10}_{n}(\theta,\sigma)w_n+P^{01}_{n}(\theta,\sigma)\bar w_n),\\
R_{2,<}&=&\sum_{n,m\in \Z_1 \atop |n|, |m|\leq E_-K_-}(P^{20}_{nm}(\theta,\sigma)w_n w_m+P^{11}_{nm}(\theta,\sigma)w_n\bar w_m+P^{02}_{nm}(\theta,\sigma)\bar w_n\bar w_m)\\
&+&\sum_{n\in \Z_1\atop E_-K_-<|n|\leq EK}(P_{nn}^{20}(\theta,\sigma)w_nw_n+P_{nn}^{11}(\theta,\sigma)w_n\bar w_n+P_{nn}^{02}(\theta,\sigma)\bar w_n\bar w_n),\\
R_{2,>}&=&\sum_{|n|>EK}(P_{nn}^{20}(\theta,\sigma)w_nw_n+P_{nn}^{11}(\theta,\sigma)w_n\bar w_n+P_{nn}^{02}(\theta,\sigma)\bar w_n\bar w_n),
 \end{eqnarray*}
 where
 $P_{l00}=P^1_{l\alpha\beta}+P^3_{l\alpha\beta}$ with $\alpha=\beta=0$;
 $P^{10}_{n}=P^3_{0\alpha\beta}$ with $\alpha=e_n,\beta=0, |n|\leq E_-K_-$;
 $P^{10}_{n}=P^2_{\alpha\beta}$ with $\alpha=e_n,\beta=0,E_-K_-<|n|\leq EK$;
 $P^{01}_{n}=P^3_{0\alpha\beta}$ with $\alpha=0,\beta=e_n, |n|\leq E_-K_-$;
 $P^{01}_{n}=P^2_{\alpha\beta}$ with $\alpha=0,\beta=e_n,E_-K_-<|n|\leq EK$;
 $P^{11}_{nm}=P^3_{0\alpha\beta}$ with $\alpha=e_n, \beta=e_m, |n|,|m| \leq E_-K_-$;
 $P^{20}_{nm}=P^3_{0\alpha\beta}$ with $\alpha=e_n+e_m, \beta=0, |n|,|m| \leq E_-K_-$;
$P^{02}_{nm}=P^3_{0\alpha\beta}$ with $\alpha=0, \beta=e_n+e_m, |n|,|m| \leq E_-K_-$;
$P^{11}_{nn}=P^1_{0\alpha\beta}$ with $\alpha=e_n, \beta=e_n, |n|>E_-K_-$;
$P^{20}_{nn}=P^1_{0\alpha\beta}$ with $\alpha=2e_n, \beta=0, |n|>E_-K_-$;
$P^{02}_{nn}=P^1_{0\alpha\beta}$ with $\alpha=0, \beta=2e_n, |n|>E_-K_-$.

 Next,we will look for an $F$ defined in a domain $D_+$ such that the time one map $\phi_F^1$ of the Hamiltonian vector field $X_F$ defines a map from $D_+\to D$ and transforms $H$ into $H_+$.
 More precisely, by second order Taylor formula, we have
\begin{eqnarray*}
H\circ\phi^1_F
&=&N+\{N,F\}+R+\int_0^1(1-t)\{\{N,F\},F\}\circ\phi^t_F dt\\
&+&\int^1_0\{R,F\}\circ\phi^t_F dt+(P-R)\circ\phi^1_F\nonumber\\
&=&N_++P_++\{N,F\}+R-\sum_{|l|\leq 1}[P_{l00}]I^l-\sum_{ |n|\leq EK}[P_{nn}^{11}(\theta,\sigma)]w_n\bar
w_n\nonumber\\
&-&\sum_{ |n|> EK}P_{nn}^{11}(\theta,\sigma)w_n\bar
w_n-\sum_{ |n|\leq E_-K_-}([P_{n(-n)}^{11}(\theta,\sigma)]w_n\bar
w_{(-n)}+[P_{(-n)n}^{11}(\theta,\sigma)]w_{(-n)}\bar
w_n)\nonumber\\
&+&\sum_{n\in \Z_1}\la\partial_{\tilde{\theta}}\Omega_n,\partial_IF_0\ra w_n\bar w_n.
\end{eqnarray*}
 We shall find a function $F(\theta, I, w,\bar w,\sigma)$ of the form
\begin{eqnarray*}
&&F=F_0+F_1+F_2,\quad F_2=F_{2,<}+F_{2,>},\\
&&F_0=\sum_{|l|\le 1} F_{l00}(\theta,\sigma)I^l,\ \ F_1=\sum_{n\in \Z_1\atop |n|\leq EK}(F^{10}_{n}(\theta,\sigma)w_n+F^{01}_{n}(\theta,\sigma)\bar w_n),\\
&&F_{2,<}=\sum_{n,m\in \Z_1\atop |n|,|m|\leq E_-K_-}(F^{20}_{nm}(\theta,\sigma)w_n w_m+F^{11}_{nm}(\theta,\sigma)w_n \bar w_m+F^{02}_{nm}(\theta,\sigma)\bar w_n\bar w_m)\\
&\;&+\sum_{n\in \Z_1 \atop E_-K_-<|n|\leq EK}(F^{20}_{nn}(\theta,\sigma)w_nw_n+F^{11}_{nn}(\theta,\sigma)w_n\bar w_n+F^{02}_{nn}(\theta,\sigma)\bar w_n\bar w_n),\\
&&F_{2,>}=\sum_{n\in \Z_1 \atop |n|>EK}(F^{20}_{nn}(\theta,\sigma)w_nw_n+F^{02}_{nn}(\theta,\sigma)\bar w_n\bar w_n),
\label{4.12}
 \end{eqnarray*}
 with $[F_0]=0$, $[F^{11}_{nm}(\theta,\sigma)]=0(|n|=|m|\leq E_-K_-)$, $[F^{11}_{nn}(\theta,\sigma)]=0(E_-K_-<|n|\leq EK)$ satisfying the equation
\begin{eqnarray}\label{4.13}
\{N,F\}+R&=&\sum_{|l|\leq 1}[P_{l00}]I^l- \sum_{n\in \Z_1}\la\partial_{\tilde{\theta}}\Omega_n,\partial_IF_0\ra w_n\bar w_n\nonumber\\
&+&\sum_{ |n|\leq E_-K_-}([P_{n(-n)}^{11}(\theta,\sigma)]w_n\bar
w_{(-n)}+[P_{nn}^{11}(\theta,\sigma)]w_n\bar
w_n+[P_{(-n)n}^{11}(\theta,\sigma)]w_{(-n)}\bar
w_n)\nonumber\\
&+&\sum_{E_-K_-<|n|\leq EK}[P_{nn}^{11}(\theta,\sigma)]w_n\bar
w_n +\sum_{|n|>EK}P_{nn}^{11}(\theta,\sigma)w_n\bar
w_n\nonumber
\end{eqnarray}
We denote that $\partial_\omega=\sum\limits_{1\leq j\leq \nu}\bar{\omega}_j\frac{\partial}{\partial_{\bar{\theta}_j}}+\sum\limits_{\nu+1\leq j\leq \nu+b}\tilde{\omega}_j\frac{\partial}{\partial_{\tilde{\theta}_j}}$,
and get the nine equations
\begin{eqnarray*}
&&\partial_\omega F_{l00} +P_{l00}=[P_{l00}],\ |l|\leq 1,\\
&&\partial_\omega F_{n}^{10}-{\rm i}(\Omega_n+d_n)F^{10}_{n}+P^{10}_{n}=0,\ n\in \Z_1,|n|\leq EK,\\
&&\partial_\omega F_{n}^{01}+{\rm i}(\Omega_n+d_n)F^{01}_{n}+P^{01}_{n}=0,\ n\in \Z_1,|n|\leq EK,\\
&&\partial_\omega F_{nm}^{20}-{\rm i}(\Omega_n+d_n)F^{20}_{nm}-{\rm i}(\Omega_m+d_m)F^{20}_{nm}+P^{20}_{nm}=0,\ n,m\in \Z_1,|n|,|m|\leq E_-K_-,\\
&&\partial_\omega F_{nm}^{11}-{\rm i}(\Omega_n+d_n)F^{11}_{nm}+{\rm i}(\Omega_m+d_m)F^{11}_{nm}+P^{11}_{nm}=0,|n|,|m|\leq E_-K_-,\\
&&\partial_\omega F_{nm}^{02}+{\rm i}(\Omega_n+d_n)F^{02}_{nm}+{\rm i}(\Omega_m+d_m)F^{02}_{nm}+P^{02}_{nm}=0,\ n,m\in \Z_1,|n|,|m|\leq E_-K_-,\\
&&\partial_\omega F_{nn}^{20}-2{\rm i}\Omega_nF^{20}_{nn}+P^{20}_{nn}=0,\ n\in \Z_1,|n|> E_-K_-,\\&&\partial_\omega F_{nn}^{11}+P^{11}_{nn}=0,\ n\in \Z_1,E_-K_-<|n|\leq EK,\\
&&\partial_\omega F_{nn}^{02}+2{\rm i}\Omega_nF^{02}_{nn}+P^{02}_{nn}=0,\ n\in \Z_1,|n|> E_-K_-,
 \end{eqnarray*}
in order to make the range of $n,m$ consistent in the above all equations, so it is feasible to combine the last three equations with the fourth, fifth and sixth equations respectively when $E_-K_-<|n|\leq EK$. Hence we rewrite them in the following
\begin{eqnarray}
&&\partial_\omega F_{l00} +P_{l00}=[P_{l00}],\ |l|\leq 1,\nonumber\\
&&\partial_\omega F_{n}^{10}-{\rm i}(\Omega_n+d_n)F^{10}_{n}+P^{10}_{n}=0, |n|\leq EK,\nonumber\\
&&\partial_\omega F_{n}^{01}+{\rm i}(\Omega_n+d_n)F^{01}_{n}+P^{01}_{n}=0, |n|\leq EK,\nonumber\\
&&\partial_\omega F_{nm}^{20}-{\rm i}(\Omega_n+d_n)F^{20}_{nm}-{\rm i}(\Omega_m+d_m)F^{20}_{nm}+P^{20}_{nm}=0, |n|,|m|\leq EK,\nonumber\\
&&\partial_\omega F_{nm}^{11}-{\rm i}(\Omega_n+d_n)F^{11}_{nm}+{\rm i}(\Omega_m+d_m)F^{11}_{nm}+P^{11}_{nm}=0,|n|,|m|\leq EK,
\label{fll}\\
&&\partial_\omega F_{nm}^{02}+{\rm i}(\Omega_n+d_n)F^{02}_{nm}+{\rm i}(\Omega_m+d_m)F^{02}_{nm}+P^{02}_{nm}=0,|n|,|m|\leq EK, \nonumber\\
&&\partial_\omega F_{nn}^{20}-2{\rm i}\Omega_nF^{20}_{nn}+P^{20}_{nn}=0,\ n\in \Z_1,|n|> EK,\nonumber\\
&&\partial_\omega F_{nn}^{02}+2{\rm i}\Omega_nF^{02}_{nn}+P^{02}_{nn}=0,\ n\in \Z_1,|n|> EK.\nonumber
\end{eqnarray}

\subsection{Estimation on the coordinate transformation}

\begin{Lemma}\label{4.3}
 Suppose that uniformly on $\Cal O_+$, $\Z^{\nu+b}=\Z^{\tilde{b}}$, $|k|\leq K, n,m\in \Z_1$,
\begin{eqnarray}
&&|\la k,\omega(\sigma)\ra|\geq \frac{\gamma}{|k|^\tau},\ k\in \Z^{\tilde b},\ |k|\neq0,\label{4.21}\\
&&|\la k,\omega(\sigma)\ra\pm (\bar{\Omega}_n+d_n)|\geq
\frac{\gamma_0}{K^\tau}, \ |n|\leq EK,\\
&&|\la k,\omega(\sigma)\ra\pm((\bar{\Omega}_n+d_n)+(\bar{\Omega}_m+d_m))|\geq
\frac{\gamma_0}{K^\tau}, \ |n|,|m|\leq EK,\\
&&|\la k,\omega(\sigma)\ra\pm((\bar{\Omega}_n+d_n)-(\bar{\Omega}_m+d_m))|\geq
\frac{\gamma_0}{K^\tau},  \begin{array}{ll}&|k|+||n|-|m||\neq 0,\\ &|n|,|m|\leq EK,\end{array}\label{c11}\\
&&|\la k,\omega(\sigma)\ra\pm 2\bar{\Omega}_n|\geq
\frac{\gamma_0\cdot |n|}{K^\tau},\ |n|>EK,\\
&&\|\tilde{\Omega}_n\|_{r,2\tau+2,\Cal O}\leq \delta_0(\gamma_0-\gamma)|n|, \label{4.23}
\end{eqnarray}
with constants $\tau\geq \tilde b=\nu+b$. If $\delta_0 $ is sufficiently small, then the linearized equation $\{N,F\}+R=\hat N$ has a solution $F$,
which is regular on $D(r,s)\times\Cal O_+$ and satisfies for $0<5\varrho<r$ the estimates
\begin{eqnarray*}
\|X_{F}\|_{s,a,\rho,D(r-3\varrho,s),\Cal O_+}&\leq&
\frac{cE^2K^{2\tau+2}}{\gamma_0^2\varrho^{\tilde{b}+1}}\cdot e^{\frac{8 E^2\delta_0(\gamma_0-\gamma)Kr}{\gamma^2}}
\|X_{R}\|_{s,\bar a,\rho,D(r,s), \Cal O_+},
\end{eqnarray*}
where the constants $c$ may be different and dependent only on $\tilde{b}$.
Besides, the error term $\hat R$ has the norm estimate
\begin{eqnarray*}
\|X_{\hat R}\|_{s,a,\rho,D(r-5\varrho,s),\Cal O_+}&\leq&
\frac{cE^2K^{2\tau+2}\delta_0}{\gamma_0\varrho^{2\tilde b+1}}e^{-K\varrho}\cdot e^{\frac{16 E^2\delta_0(\gamma_0-\gamma)Kr}{\gamma^2}}
\|X_{R}\|_{s,\bar a,\rho,D(r,s), \Cal O_+}.
\end{eqnarray*}
\end{Lemma}
\proof
Firstly we consider the most complicated equations in $(\ref{fll})$ with $|n|, |m|\leq EK$
\begin{eqnarray}\label{C1}
\partial_\omega F^{11}_{nm}(\theta,\sigma)-{\rm i}((\Omega_n(\theta,\sigma)+d_n)-(\Omega_m(\theta,\sigma)+d_m)) F^{11}_{nm}(\theta,\sigma)+ P^{11}_{nm}(\theta,\sigma)=0,
\end{eqnarray}
Let $ \partial_\omega T^{11}_{nm}(\theta,\sigma)=\Gamma_K(\tilde{\Omega}_n(\theta,\sigma)-\tilde{\Omega}_m(\theta,\sigma)$), $F^{11}_{nm}=e^{{\rm i}T^{11}_{nm}}\tilde F_{nm}^{11}$, $P^{11}_{nm}=e^{{\rm i}T^{11}_{nm}}\tilde P_{nm}^{11}$, then $(\ref{C1})$ is transformed into
\begin{eqnarray*}
\partial_\omega \tilde{F}^{11}_{nm}-{\rm i}((\bar{\Omega}_n(\sigma)+d_n)-(\bar{\Omega}_m(\sigma)+d_m))\tilde{ F}^{11}_{nm}-(1-\Gamma_{K})(\tilde{\Omega}_n(\theta)-\tilde{\Omega}_m(\theta))\tilde{ F}^{11}_{nm}+ \tilde{P}^{11}_{nm}=0.
\end{eqnarray*}
We only solve the truncation equation
\begin{eqnarray}\label{11}
\Gamma_{K}(\partial_\omega \tilde{F}^{11}_{nm}-{\rm i}((\bar{\Omega}_n(\sigma)+d_n)-(\bar{\Omega}_m(\sigma)+d_m))\tilde{ F}^{11}_{nm}+ \tilde{P}^{11}_{nm})=0,
\end{eqnarray}
and the error term is
\begin{eqnarray}{\label{er1}}
\hat R^{11}_{nm}=e^{{\rm i}T^{11}_{nm}}[(1-\Gamma_K)(e^{-{\rm i}T^{11}_{nm}}P_{nm}^{11})+{\rm i}(1-\Gamma_K)(\tilde{\Omega}_n-\tilde{\Omega}_m)e^{-{\rm i}T^{11}_{nm}}F_{nm}^{11}].
\end{eqnarray}
To solve the equation $(\ref{11})$, we expand $\tilde{F}^{11}_{nm},\tilde{P}^{11}_{nm}$ into Fourier series
\begin{eqnarray*}
\Gamma_{K}\tilde{F}^{11}_{nm}(\theta,\sigma)=\sum_{|k|\leq K}\tilde{F}^{11}_{knm}\kth,\quad
\Gamma_{K}\tilde{P}^{11}_{nm}(\theta,\sigma)=\sum_{|k|\leq K}\tilde{P}^{11}_{knm}\kth,
\end{eqnarray*}
and substitute them into the equation $(\ref{11})$
\begin{eqnarray*}
{\rm i}\la k,\omega\ra \tilde{F}^{11}_{knm}(\sigma)-{\rm i}((\bar{\Omega}_n(\sigma)+d_n)-(\bar{\Omega}_m(\sigma)+d_m))\tilde{ F}^{11}_{knm}(\sigma)+ \tilde{P}^{11}_{knm}(\sigma)=0,
\end{eqnarray*}
we can easily get
\begin{eqnarray*}
\tilde{F}^{11}_{knm}(\sigma)={\rm i}\frac{\tilde{P}^{11}_{knm}(\sigma)}{\la k,\omega\ra-(\bar{\Omega}_n(\sigma)+d_n)+(\bar{\Omega}_m(\sigma)+d_m)},\ \ |k|+||n|-|m||\neq 0, \ |n|,|m|\leq EK,
\end{eqnarray*}
by the condition $(\ref{c11})$ and the assumption $(A1),(A2)$, then
\begin{eqnarray*}
&&|\tilde{F}^{11}_{knm}|_{\Cal O_+}=\sup_{\sigma\in\Cal O}\mid {\rm i}\frac{\tilde{P}^{11}_{knm}}{\la k,\omega\ra-(\bar{\Omega}_n(\sigma)+d_n)+(\bar{\Omega}_m(\sigma)+d_m)}\mid\\
&&+\mid \frac{\partial}{\partial\sigma}({\rm i}\frac{\tilde{P}^{11}_{knm}}{\la k,\omega\ra-(\bar{\Omega}_n(\sigma)+d_n)+(\bar{\Omega}_m(\sigma)+d_m)})\mid\nonumber\\
&\leq&\sup_{\sigma\in\Cal O}\left(\frac{K^\tau}{\gamma_0}(|\tilde {P}^{11}_{knm}|+|\frac{\partial}{\partial\sigma}\tilde {P}^{11}_{knm}|)+\frac{K^{2\tau}}{\gamma_0^2}|\tilde {P}^{11}_{knm}|[|\frac{\partial}{\partial\sigma}\la k,\omega\ra|
+|\frac{\partial}{\partial\sigma}((\bar{\Omega}_n(\sigma)+d_n)-(\bar{\Omega}_m(\sigma)+d_m))|]\right)\\
&\leq &\frac{K^\tau}{\gamma_0}|\tilde {P}^{11}_{knm}|_{\Cal O}+\sup_{\sigma\in\Cal O}\left(\frac{K^{2\tau}}{\gamma_0^2}|\tilde {P}^{11}_{knm}|(EK+c(n+m)\varepsilon_0)\right)\\
&\leq&\frac{cEK^{2\tau+1}}{\gamma_0^2}|\tilde {P}^{11}_{knm}|_{\Cal O},
\end{eqnarray*}
and the estimate of the function $\tilde F^{11}_{nm}$ is
\begin{eqnarray}{\label{4.24}}
\|\tilde{F}^{11}_{nm}\|_{D(r-\varrho),\Cal O_+}&\leq& \sum_{|k|\leq K}|\tilde F^{11}_{knm}|_{\Cal O}e^{|k|(r-\varrho)}
\leq \frac{cEK^{2\tau+1}}{\gamma_0^2}\sum_{|k|\leq K}|\tilde{P}^{11}_{knm}|_{\Cal O}e^{|k|(r-\varrho)}\nonumber\\
&\leq& \frac{cEK^{2\tau+1}}{\gamma_0^2}\cdot \frac{(2+2e)^{\tilde b}}{\varrho^{\tilde b}}\|\tilde {P}^{11}_{nm}\|_{D(r),\Cal O}\nonumber\\
&\leq& \frac{c(\tilde b)EK^{2\tau+1}}{\gamma_0^2\varrho^{\tilde b}}\|\tilde {P}^{11}_{nm}\|_{D(r),\Cal O},
\end{eqnarray}
where $ c(\tilde{b})=c\cdot(2+2e)^{\tilde{b}}$ is a constant.\\
In the following we will estimate $F^{11}_{nm}$. Since
$\partial_\omega T^{11}_{nm}(\theta,\sigma)=\Gamma_{K}(\tilde{\Omega}_n(\theta,\sigma)-\tilde{\Omega}_m(\theta,\sigma))$, we expand $T^{11}_{nm}(\theta,\sigma),\tilde{\Omega}_n(\theta,\sigma),\tilde{\Omega}_m(\theta,\sigma)$ into Fourier series
\begin{eqnarray*}
T^{11}_{nm}(\theta,\sigma)=\sum_{|k|\neq 0}T^{11}_{knm}(\sigma)\kth, \tilde{\Omega}_n(\theta,\sigma)=\sum_{|k|\neq 0}\tilde{\Omega}_{kn}\kth,\tilde{\Omega}_m(\theta,\sigma)=\sum_{|k|\neq 0}\tilde{\Omega}_{km}\kth,
\end{eqnarray*}
and obtain
\begin{eqnarray*}
&&{\rm i}\la k,\omega\ra T^{11}_{knm}(\sigma)=\tilde{\Omega}_{kn}-\tilde{\Omega}_{km},\quad
T^{11}_{knm}(\sigma)=\frac{\tilde{\Omega}_{kn}-\tilde{\Omega}_{km}}{{\rm i}\la k,\omega\ra},\ 0<|k|\leq K,\\
&&T^{11}_{nm}(\theta,\sigma)=\sum_{0<|k|\leq K}\frac{\tilde{\Omega}_{kn}-\tilde{\Omega}_{km}}{{\rm i}\la k,\omega\ra}\kth.
\end{eqnarray*}
Let $\theta=\theta_1+{\rm i}\theta_2,\theta_1,\theta_2 \in \T^{\tilde{b}}$ and we denote
\begin{eqnarray*}
&&T^{11}_{nm,1}(\theta_1,\sigma)=\sum_{0<|k|\leq K}\frac{\tilde{\Omega}_{kn}-\tilde{\Omega}_{km}}{{\rm i}\la k,\omega\ra}e^{{\rm i}\la k,\theta_1\ra},\\
&&T^{11}_{nm,2}(\theta,\sigma)=T^{11}_{nm}(\theta,\sigma)-T^{11}_{nm,1}(\theta_1,\sigma)=\sum_{0<|k|\leq K}\frac{\tilde{\Omega}_{kn}-\tilde{\Omega}_{km}}{{\rm i}\la k,\omega\ra}e^{{\rm i}\la k,\theta_1\ra}(e^{-\la k,\theta_2\ra}-1),
\end{eqnarray*}
since $\tilde \Omega_n,\tilde \Omega_m$ is real analytic, so is $T^{11}_{nm,1}(\theta_1,\sigma)$. \\
Meanwhile, by the condition $(\ref{4.21}),(\ref{4.23})$ and the assumption $(A1)$, we have
\begin{eqnarray*}
|\frac{\tilde{\Omega}_{kn}-\tilde{\Omega}_{km}}{{\rm i}\la k,\omega\ra}|_{\Cal O_+}&\leq& \sup_{\sigma\in \Cal O}(\frac{|k|^\tau}{\gamma}(|\tilde{\Omega}_{kn}-\tilde{\Omega}_{km}|+|\frac{\partial}{\partial\sigma}(\tilde{\Omega}_{kn}-\tilde{\Omega}_{km})|)\\
&+&\frac{|k|^{2\tau}}{\gamma^2}|\frac{\partial}{\partial\sigma}\la k,\omega\ra|\cdot|\tilde{\Omega}_{kn}-\tilde{\Omega}_{km}|)\\
&\leq& \frac{E |k|^{2\tau+1}}{\gamma^2}|\tilde{\Omega}_{kn}-\tilde{\Omega}_{km}|_{\Cal O},
\end{eqnarray*}
 and the estimate of the transformation $T^{11}_{nm}(\theta,\sigma)$
\begin{eqnarray}
&&\|{\mathrm{Im}} T^{11}_{nm}(\theta,\sigma)\|_{D(r),\Cal O_+}=\|{\mathrm{Im }}T^{11}_{nm,2}(\theta,\sigma)\|_{D(r),\Cal O_+}\nonumber\\
&\leq& \sum_{0<|k|\leq K}|\frac{\tilde{\Omega}_{kn}-\tilde{\Omega}_{km}}{{\rm i}\la k,\omega\ra}|_{\Cal O_+}\cdot|e^{-\la k,\theta_2\ra}-1|\nonumber\\
&\leq& \frac{E}{\gamma^2}\sum_{0<|k|\leq K}|k|^{2\tau+1}|\tilde{\Omega}_{kn}-\tilde{\Omega}_{km}|_{\Cal O}\cdot e^{|k|r}\cdot |k|r\nonumber\\
&\leq& \frac{Er}{\gamma^2}\cdot (\|\tilde{\Omega}_{n}\|_{r,2\tau+2,\Cal O}+\|\tilde{\Omega}_{m}\|_{r,2\tau+2,\Cal O})\nonumber\\
&\leq& \frac{2E^2\delta_0(\gamma_0-\gamma)Kr}{\gamma^2}.\label{T11}
\end{eqnarray}
Then we can easily obtain the estimate of $F_{nm}^{11}$
\begin{eqnarray*}
\|F_{nm}^{11}\|_{D(r-2\varrho),\Cal O_+}&=&\|e^{{\rm i}T^{11}_{nm}}\tilde F_{nm}^{11}\|_{D(r-2\varrho),\Cal O_+}\\
&\leq& e^{2\|{\mathrm{Im}} T^{11}_{nm}(\theta,\sigma)\|_{D(r),\Cal O_+}}\cdot \|\tilde F_{nm}^{11}\|_{D(r-2\varrho),\Cal O_+}\\
&\leq& e^{\frac{4 E^2\delta_0(\gamma_0-\gamma)Kr}{\gamma^2}}\|\tilde F_{nm}^{11}\|_{D(r-2\varrho),\Cal O_+},
\end{eqnarray*}
and similarly the estimate of $\tilde P_{nm}^{11}$
\begin{eqnarray*}
\|\tilde P_{nm}^{11}\|_{D(r),\Cal O_+}&=&\|e^{-{\rm i}T^{11}_{nm}}P_{nm}^{11}\|_{D(r),\Cal O_+}\leq e^{2\|{\mathrm{Im}} T(\theta,\sigma)\|_{D(r),\Cal O}}\cdot \|P_{nm}^{11}\|_{D(r),\Cal O}\\
&\leq &e^{\frac{4 E^2\delta_0(\gamma_0-\gamma)Kr}{\gamma^2}}\| P_{nm}^{11}\|_{D(r),\Cal O},
\end{eqnarray*}
so finally associated with $(\ref{4.24})$ we obtain
\begin{eqnarray}{\label{F11}}
\|F_{nm}^{11}\|_{D(r-2\varrho),\Cal O_+}&\leq& e^{\frac{4 E^2\delta_0(\gamma_0-\gamma)Kr}{\gamma^2}}\|\tilde F_{nm}^{11}\|_{D(r-2\varrho),\Cal O_+}\nonumber\\
&\leq & e^{\frac{4 E^2\delta_0(\gamma_0-\gamma)Kr}{\gamma^2}}\cdot\frac{cEK^{2\tau+1}}{\gamma_0^2\varrho^{\tilde{b}}}\|\tilde{P}^{11}_{nm}\|_{D(r),\Cal O_+}\nonumber\\
&\leq&\frac{cEK^{2\tau+1}}{\gamma_0^2\varrho^{\tilde{b}}}\cdot e^{\frac{8 E^2\delta_0(\gamma_0-\gamma)Kr}{\gamma^2}}\|{P}^{11}_{nm}\|_{D (r),\Cal O},
\end{eqnarray}
where $ c(\tilde{b})=c\cdot(2+2e)^{\tilde{b}}$.\\
Besides, we need to estimate the error term $\hat R^{11}_{nm}$ in $(\ref{er1})$. First, for any analytic function $h(\theta,\sigma)$ defined in $D(r)\times \Cal O$, we give an inequality
$$ \|(1-\Gamma_K)h(\theta,\sigma)\|_{D(r-2\varrho),\Cal O}\leq ce^{-K\varrho}\|h\|_{D(r),\Cal O}, \ c=\frac{(2+2e)^{\tilde b}}{\varrho^{\tilde b}}.$$
Indeed, this inequality can be easily proved
\begin{eqnarray*}
&&\|(1-\Gamma_K)h(\theta,\sigma)\|_{D(r-2\varrho),\Cal O}
=\|\sum_{|k|>K}h_{k}(\sigma)\kth\|_{D(r-2\varrho),\Cal O}\\
&\leq&\sum_{|k|>K}|h_k|_{\Cal O}e^{|k|(r-2\varrho)}
\leq e^{-K\varrho}\sum_{|k|>K}|h_k|_{\Cal O}e^{|k|(r-\varrho)}\\
&\leq&\frac{(2+2e)^{\tilde b}}{\varrho^{\tilde b}}e^{-K\varrho}\|h\|_{D(r),\Cal O}=ce^{-K\varrho}\|h\|_{D(r),\Cal O}.
\end{eqnarray*}
In this way, by $(\ref{T11}),(\ref{F11})$, the estimate of the error term $\hat R^{11}_{nm}$ is
\begin{eqnarray}{\label{R11}}
&&\|e^{{\rm i}T^{11}_{nm}}[(1-\Gamma_K)(e^{-{\rm i}T^{11}_{nm}}P_{nm}^{11})+{\rm i}(1-\Gamma_K)(\tilde{\Omega}_n-\tilde{\Omega}_m)e^{-{\rm i}T^{11}_{nm}}F_{nm}^{11}]\|_{D(r-4\varrho),\Cal O_+}\nonumber\\
&\leq&  \frac{(2+2e)^{\tilde b}}{\varrho^{\tilde b}}e^{-K\varrho}\cdot e^{4\|\mathrm{Im}T^{11}_{nm}\|_{D(r-4\varrho),\Cal O_+}}(\|P^{11}_{nm}\|_{D(r-2\varrho),\Cal O}\nonumber\\
&+&\|\tilde{\Omega}_n-\tilde{\Omega}_m\|_{D(r-2\varrho),\Cal O}\|F_{nm}^{11}\|_{D(r-2\varrho),\Cal O_+})\nonumber\\
&\leq& \frac{(2+2e)^{\tilde b}}{\varrho^{\tilde b}}e^{-K\varrho}\cdot e^{\frac{8 E^2\delta_0(\gamma_0-\gamma)Kr}{\gamma^2}}
(\|P^{11}_{nm}\|_{D(r-2\varrho),\Cal O}\nonumber\\
&+&(n+m)\delta_0(\gamma_0-\gamma)\cdot \frac{cEK^{2\tau+1}}{\gamma_0^2\varrho^{\tilde{b}}}\cdot e^{\frac{8 E^2\delta_0(\gamma_0-\gamma)Kr}{\gamma^2}}\|{P}^{11}_{nm}\|_{D (r),\Cal O})\nonumber\\
&\leq&  \frac{cE^2K^{2\tau+2}\delta_0}{\gamma_0\varrho^{2\tilde b}}e^{-K\varrho}\cdot e^{\frac{16 E^2\delta_0(\gamma_0-\gamma)Kr}{\gamma^2}}\|P^{11}_{nm}\|_{D(r),\Cal O},
\end{eqnarray}
the estimate of the $F^{20}_{nm},F^{02}_{nm}$ and their error term $\hat R^{20}_{nm}, \hat R^{02}_{nm}$ can be  similarly obtained. According to all the above estimates of terms in $F_{2,<}$, we now compute the vector field norm of $X_{F_{2,<}}$ namely
\begin{eqnarray*}
&&\|X_{F_{2,<}}\|_{s, a,\rho, D(r-3\varrho,s),\Cal O_+}\\
&=&\frac{1}{s^2}\|(F_{2,<})_{\theta}\|_{D(r-3\varrho,s),\Cal O_+}+\frac{1}{s}(\sum_{n\leq EK}\|(F_{2,<})_{w_n}\|_{D(r-3\varrho),\Cal O_+}n^{a}e^{n\rho}\\
&+&\|(F_{2,<})_{\bar {w}_n}\|_{D(r-3\varrho),\Cal O_+}n^{a}e^{n\rho}).
\end{eqnarray*}
For the first term $\|(F_{2,<})_{\theta}\|_{D(r-3\varrho,s),\Cal O_+} $, we have
\begin{eqnarray*}
&&\|(F_{2,<})_{\theta}\|_{D(r-3\varrho,s),\Cal O_+}=\sum_{1\leq j\leq \tilde b}\|(F_{2,<})_{\theta_j}\|_{D(r-3\varrho,s),\Cal O_+},\\
&&\|(F_{2,<})_{\theta_j}\|_{D(r-3\varrho,s),\Cal O_+}\\
&=&\sup_{\|w\|_{a,\rho}<s\atop\|\bar w\|_{a,\rho}<s}\sum_{|n|,|m|\leq EK}
[\|(F^{20}_{nm})_{\theta_j}\|_{D(r-3\varrho),\Cal O_+}|w_n||w_m|+\|(F^{02}_{nm})_{\theta_j}\|_{D(r-3\varrho),\Cal O_+}|\bar w_n||\bar w_m|]\\
&+&\sum_{|n|,|m|\leq EK}\|(F^{11}_{nm})_{\theta_j}\|_{D(r-3\varrho),\Cal O_+}|w_n||\bar w_m|,
\end{eqnarray*}
by Lemma $\ref{Lem7.2}$ and the  estimate in $(\ref{F11})$, one have
$$\|(F^{11}_{nm})_{\theta_j}\|_{D(r-3\varrho),\Cal O_+}\leq \varrho^{-1}\|F^{11}_{nm}\|_{D(r-2\varrho),\Cal O_+}\leq
\frac{cEK^{2\tau+1}}{\gamma_0^2\varrho^{\tilde{b}+1}}\cdot e^{\frac{8 E^2\delta_0(\gamma_0-\gamma)Kr}{\gamma^2}}\|{P}^{11}_{nm}\|_{D (r),\Cal O}$$
and the $\|(F^{20}_{nm})_{\theta_j}\|_{D(r-3\varrho),\Cal O_+}$, $\|(F^{02}_{nm})_{\theta_j}\|_{D(r-3\varrho),\Cal O_+}$
have the same estimate by the similar argument. Then the estimate of $\|(F_{2,<})_{\theta}\|_{D(r-3\varrho,s),\Cal O_+}$ is obtained
\begin{eqnarray*}
&&\|(F_{2,<})_{\theta}\|_{D(r-3\varrho,s),\Cal O_+}\\
&\leq& \frac{cEK^{2\tau+1}}{\gamma_0^2\varrho^{\tilde{b}+1}}\cdot e^{\frac{8 E^2\delta_0(\gamma_0-\gamma)Kr}{\gamma^2}}\sup_{\|w\|_{a,\rho}<s\atop\|\bar w\|_{a,\rho}<s}\sum_{n,m\leq EK}
[\|P^{20}_{nm}\|_{D(r),\Cal O}|w_n||w_m|\\
&+&\|P^{02}_{nm}\|_{D(r),\Cal O}|\bar w_n||\bar w_m|]+\sum_{n\neq m\atop n,m\leq EK}\|P^{11}_{nm}\|_{D(r),\Cal O}|w_n||\bar w_m|.
\end{eqnarray*}
Similarly, the norms of the term $(F_{2,<})_{w_n},(F_{2,<})_{\bar w_n}$ respectively satisfy
\begin{eqnarray*}
&&\|(F_{2,<})_{w_n}\|_{D(r-3\varrho,s),\Cal O_+}\\
&\leq &\sup_{\|w\|_{a,\rho}<s\atop\|\bar w\|_{a,\rho}<s}\sum_{|m|\leq EK}(\|F^{20}_{nm}\|_{D(r-3\varrho),\Cal O_+}|w_m|+\|F^{11}_{nm}\|_{D(r-3\varrho),\Cal O_+}|\bar w_m|)\\
&\leq& \frac{cEK^{2\tau+1}}{\gamma_0^2\varrho^{\tilde{b}}}\cdot e^{\frac{8 E^2\delta_0(\gamma_0-\gamma)Kr}{\gamma^2}}\sup_{\|w\|_{a,\rho}<s\atop\|\bar w\|_{a,\rho}<s}\sum_{|m|\leq EK}(\|P^{20}_{nm}\|_{D(r),\Cal O}|w_m|
+\|P^{11}_{nm}\|_{D(r),\Cal O}|\bar w_m|),\\
&&\|(F_{2,<})_{\bar w_n}\|_{D(r-3\varrho,s),\Cal O_+}\\
&\leq &\sup_{\|w\|_{a,\rho}<s\atop\|\bar w\|_{a,\rho}<s}\sum_{|m|\leq EK}(\|F^{02}_{nm}\|_{D(r-3\varrho),\Cal O_+}|\bar w_m|+\|F^{11}_{mn}\|_{D(r-3\varrho),\Cal O_+}|w_m|)\\
&\leq& \frac{cEK^{2\tau+1}}{\gamma_0^2\varrho^{\tilde{b}}}\cdot e^{\frac{8 E^2\delta_0(\gamma_0-\gamma)Kr}{\gamma^2}}\sup_{\|w\|_{a,\rho}<s\atop\|\bar w\|_{a,\rho}<s}\sum_{|m|\leq EK}(\|P^{02}_{nm}\|_{D(r),\Cal O}|\bar w_m|
+\|P^{11}_{mn}\|_{D(r),\Cal O}|w_m|).
\end{eqnarray*}
Associated with the above estimates of the terms $(F_{2,<})_{\theta},(F_{2,<})_{w_n},(F_{2,<})_{\bar w_n}$, we finally get the norm of the vector field $X_{F_{2,<}}$
\begin{eqnarray}{\label{XF1}}
&&\|X_{F_{2,<}}\|_{s, a,\rho, D(r-3\varrho,s),\Cal O_+}\nonumber\\
&\leq& \frac{cEK^{2\tau+1}}{\gamma_0^2\varrho^{\tilde{b}+1}}\cdot e^{\frac{8 E^2\delta_0(\gamma_0-\gamma)Kr}{\gamma^2}}\cdot\frac{1}{s^2}\sup_{\|w\|_{a,\rho}<s\atop\|\bar w\|_{a,\rho}<s}(\sum_{|n|,|m|\leq EK}
[\|P^{20}_{nm}\|_{D(r),\Cal O}|w_n||w_m|\nonumber\\
&+&\|P^{02}_{nm}\|_{D(r),\Cal O}|\bar w_n||\bar w_m|]+\sum_{|n|,|m|\leq EK}\|P^{11}_{nm}\|_{D(r),\Cal O}|w_n||\bar w_m|)\nonumber\\
&+&\frac{cEK^{2\tau+1}}{\gamma_0^2\varrho^{\tilde{b}}}\cdot e^{\frac{8 E^2\delta_0(\gamma_0-\gamma)Kr}{\gamma^2}}\cdot \frac{1}{s}\sup_{\|w\|_{a,\rho}<s\atop\|\bar w\|_{a,\rho}<s}(\sum_{|n|,|m|\leq EK}[\|P^{20}_{nm}\|_{D(r),\Cal O}|w_m|n^ae^{n\rho}\nonumber\\
&+&\|P^{02}_{nm}\|_{D(r),\Cal O}|\bar w_m|n^ae^{n\rho}]
+\sum_{|n|,|m| \leq EK}\|P^{11}_{mn}\|_{D(r),\Cal O}|w_m|n^ae^{n\rho})\nonumber\\
&\leq& \frac{cE^2K^{2\tau+2}}{\gamma_0^2\varrho^{\tilde{b}+1}} e^{\frac{8 E^2\delta_0(\gamma_0-\gamma)Kr}{\gamma^2}}\|X_{R_{2,<}}\|_{s,\bar a,\rho,D(r,s),\Cal O}.
\end{eqnarray}
 With the similar arguments of $F_{2,<}$, the error term $\hat {R}_{2,<}$ is represented as
 $$\hat R_{2,<}=\sum_{|n|,|m|\leq Ek}(\hat {R}^{20}_{nm}w_nw_m+\hat {R}^{02}_{nm}\bar w_n\bar w_m)+\sum_{|n|,|m|\leq EK}\hat R^{11}_{nm}w_n\bar w_m,$$
 where $\hat {R}^{11}_{nm}$  defined in $ (\ref{er1})$ and $\hat R^{20}_{nm},\hat R^{02}_{nm}$ have the similar formulas
$$\hat R^{20}_{nm}=e^{{\rm i}T^{20}_{nm}}[(1-\Gamma_K)(e^{-{\rm i}T^{20}_{nm}}P_{nm}^{20})+{\rm i}(1-\Gamma_K)(\tilde{\Omega}_n+\tilde{\Omega}_m)e^{-{\rm i}T^{20}_{nm}}F_{nm}^{20}],$$
$$\hat R^{02}_{nm}=e^{{\rm i}T^{02}_{nm}}[(1-\Gamma_K)(e^{-{\rm i}T^{02}_{nm}}P_{nm}^{02})-{\rm i}(1-\Gamma_K)(\tilde{\Omega}_n+\tilde{\Omega}_m)e^{-{\rm i}T^{02}_{nm}}F_{nm}^{02}].$$
 We repeat the same calculation process of $ X_{F_{2,<}}$ and finally get the estimate of the vector field $X_{\hat {R}_{2,<}}$
\begin{eqnarray}{\label{xer1}}
\|X_{\hat {R}_{2,<}}\|_{\!{}_{s,\bar a,\rho,D(r-5\varrho,s),\Cal O_+}}\leq \frac{cE^2K^{2\tau+2}\delta_0}{\gamma_0\varrho^{2\tilde b+1}}e^{-K\varrho} e^{\frac{16 E^2\delta_0(\gamma_0-\gamma)Kr}{\gamma^2}}\|X_{R_{2,<}}\|_{\!{}_{s,\bar a,\rho,D(r,s),\Cal O}}.
\end{eqnarray}
For $|n|> EK$, we have to solve two equations
\begin{eqnarray*}
&&\partial_\omega F_{n}^{20}-2{\rm i}\Omega_nF^{20}_{n}+P^{20}_{n}=0,\quad n\in \Z_1,|n|> EK,\\
&&\partial_\omega F_{n}^{02}+2{\rm i}\Omega_nF^{02}_{n}+P^{02}_{n}=0,\quad n\in \Z_1,|n|> EK.
\end{eqnarray*}
It is sufficient to solve the first one and the second can be similarly solved. For the first one, we solve the truncation equation
\begin{eqnarray}\label{C2}
-{\rm i}\partial_\omega F^{20}_{nn}-2\bar \Omega_nF^{20}_{nn}-2\Gamma_{K}(\tilde{\Omega}_nF^{20}_{nn})=\Gamma_{K}({\rm i}P^{20}_{nn}), \ \Gamma_{K}F^{20}_{nn}=F^{20}_{nn},
\end{eqnarray}
and the error term $\hat R^{20}=\sum\limits_{|n|>EK}\hat{R}^{20}_{nn}w_nw_n$ with the elements defined by
 \begin{eqnarray}\label{er3}
 \hat{R}^{20}_{nn}=(1-\Gamma_{K})({\rm i}P^{20}_{nn}+2\tilde\Omega_nF^{20}_{nn}),\ |n| > EK.
\end{eqnarray}
We expand $F^{20}_{nn}(\theta,\sigma)$, $\tilde \Omega_n(\theta,\sigma)$, $P^{20}_{nn}(\theta,\sigma)$ into Fourier series
$$F^{20}_{nn}=\sum_{|k|\leq K}F^{20}_{knn}\kth,\ \tilde \Omega_n=\sum_{|k|\neq 0}\tilde \Omega_{kn}\kth, \ P^{20}_{nn}=\sum_{k\in \Z^{\tilde b}}P^{20}_{knn}\kth$$
and the equation $(\ref{C2})$ is represented as
$$\sum_{|k|\leq K}(\la k,\omega\ra -2\bar \Omega_n)F^{20}_{knn}\kth-2\sum_{|k|\leq K}(\sum_{|l|\leq |k|}\tilde \Omega_{k-l,n}F^{20}_{lnn})\kth={\rm i}\sum_{|k|\leq K}P^{20}_{knn}\kth. $$
We introduce the following denotations for simplicity,
$$\Lambda_n=\mathrm{diag}(\la k,\omega\ra-2\bar \Omega_n)_{|k|\leq K},\ \ D_n=(-2\tilde \Omega_{k-l,n})_{|k|,|l|\leq K},$$
$$ \hat F^{20}_{n}=(F^{20}_{knn})_{|k|\leq K},\ \ \hat P^{20}_{n}=({\rm i}P^{20}_{knn})_{|k|\leq K},$$
so the above equation is equivalence to
\begin{eqnarray}{\label{C3}}
(\Lambda_n+D_n) \hat F^{20}_{n}= \hat P^{20}_{n},\ |n|>EK.
\end{eqnarray}
According to the assumption $(A1),(A2)$,  $|\la k, \omega\ra|\leq EK<|n|$, $|\bar \Omega_n|\geq |n|-c\varepsilon_0 |n|\geq \frac{3|n|}{4}$ if $\varepsilon_0\ll \frac{1}{4}$ small enough, it is clear that
$$|\la k, \omega\ra-2\bar\Omega_n|\geq 2|\bar \Omega_n|-|\la k,\omega\ra|\geq  2|\bar \Omega_n|-|n|\geq \frac{|n|}{2}.$$
Moreover, we denote a matrix $A_{\tilde{r}}=\mathrm{diag}(e^{|k|\tilde{r}})_{|k|\leq K}$ with $0<\tilde r< r$ and multiply $(\ref{C3})$ in the left by $A_{\tilde{r}}$
$$(\Lambda_n+A_{\tilde{r}}D_nA^{-1}_{\tilde{r}})A_{\tilde{r}}\hat F^{20}_{n}= A_{\tilde{r}}\hat P^{20}_{n}.$$
It is obvious that the matrix norm of $\Lambda_n^{-1}$ is
\begin{eqnarray}{\label{4.25}}
\|\Lambda_n^{-1}\|_{\Cal O}&=&\max_{|k|\leq K}\sup_{\sigma\in \Cal O}(|\frac{1}{\la k,\omega\ra-2\bar \Omega_n}|+|\frac{\partial}{\partial\sigma}\frac{1}{\la k,\omega\ra-2\bar \Omega_n}|)\nonumber\\
&\leq& \frac{2}{|n|}+\frac{4}{n^2}(KE+2c\varepsilon_0|n|)\leq \frac{2}{|n|}+\frac{4}{n^2}\cdot \frac{3|n|}{2}\leq \frac{8}{|n|}.
\end{eqnarray}
By the condition $(\ref{4.23})$, the norm of $A_{\tilde r}D_nA^{-1}_{\tilde r}$ is
\begin{eqnarray}{\label{4.26}}
\|A_{\tilde r}D_nA^{-1}_{\tilde r}\|_{\Cal O}&=&\max_{|k|\leq K}\sup_{\sigma\in \Cal O}\sum_{|l|\leq K}2(e^{(|l|-|k|)\tilde r}|\tilde \Omega_{l-k,n}|+|\frac{\partial}{\partial\sigma}(e^{(|l|-|k|)\tilde r}\tilde \Omega_{l-k,n})|)\nonumber\\
&\leq & 2\max_{|k|\leq K}\sum_{|l|\leq K}e^{(|l|-|k|)\tilde r}|\tilde \Omega_{l-k,n}|_{\Cal O}
\leq 2\sum_{|k|\leq K}e^{|k|\tilde r}|\tilde \Omega_{k,n}|_{\Cal O}\nonumber\\
&\leq & 2\|\tilde \Omega_n\|_{r,2\tau+2,\Cal O}\leq 2|n|\delta_0(\gamma_0-\gamma),
\end{eqnarray}
then associated with $(\ref{4.25}),(\ref{4.26})$, and if  $\delta_0\gamma_0\ll \frac{1}{32}$ is small enough, we have
\begin{eqnarray*}
\|\Lambda_n^{-1}(A_{\tilde r}D_nA^{-1}_{\tilde r})\|_{\Cal O}\leq \|\Lambda_n^{-1}\|_{\Cal O}\cdot\|A_{\tilde r}D_nA^{-1}_{\tilde r}\|_{\Cal O}\leq \frac{8}{|n|}\cdot2|n|\delta_0(\gamma_0-\gamma)<\frac{1}{2},
\end{eqnarray*}
with this condition, the matrix $\Lambda_n+ A_{\tilde r}D_nA^{-1}_{\tilde r}$ is invertible and its inverse matrix has the norm estimate
\begin{eqnarray*}
\|(\Lambda_n+ A_{\tilde r}D_nA^{-1}_{\tilde r})^{-1}\|_{\Cal O}\leq \|\Lambda_n^{-1}\|_{\Cal O}\cdot\frac{1}{1-\|\Lambda_n^{-1}(A_{\tilde r}D_nA^{-1}_{\tilde r})\|_{\Cal O}}\leq \frac{16}{|n|},
\end{eqnarray*}
hence the estimate of $F^{20}_{nn}$ is
\begin{eqnarray*}
&&\|F^{20}_{nn}\|_{D(\tilde r),\Cal O}\leq \sum_{|k|\leq K}|F^{20}_{knn}\|_{\Cal O}e^{|k|\tilde r}=\|A_{\tilde r}\hat F^{20}_n\|_{\Cal O}\\
&=&\|(\Lambda_n+ A_{\tilde r}D_nA^{-1}_{\tilde r})^{-1}A_{\tilde r}\hat P^{20}_{n}\|_{\Cal O}\\
&\leq&\|(\Lambda_n+ A_{\tilde r}D_nA^{-1}_{\tilde r})^{-1}\|_{\Cal O}\cdot\|A_{\tilde r}\hat P^{20}_{n}\|_{\Cal O}\\
&\leq& \frac{16}{|n|}\sum_{|k|\leq K}|P^{20}_{knn}|_{\Cal O}e^{|k|\tilde r}
\leq \frac{16(2+2e)^{\tilde b}}{|n|(r-\tilde r)^{\tilde b}}\|P^{20}_{nn}\|_{D(r),\Cal O},
\end{eqnarray*}
we take $\tilde r=r-2\varrho$ and finally get the estimate of $F^{20}_{nn}$
\begin{eqnarray}{\label{F20}}
\|F^{20}_{nn}\|_{D(r-2\varrho),\Cal O}\leq \frac{16(1+e)^{\tilde b}}{|n|\varrho^{\tilde b}}\|P^{20}_{nn}\|_{D(r),\Cal O}.
\end{eqnarray}
In addition, the error term $\hat R^{20}_{nn}$ in $(\ref{er3})$ has the following estimate
\begin{eqnarray}{\label{er2}}
&&\|(1-\Gamma_K)(2\tilde \Omega_mF^{20}_{nn}+{\rm i}P^{20}_{nn})\|_{D(r-4\varrho),\Cal O}\nonumber\\
&\leq& e^{-K\varrho}\cdot \frac{(1+e)^{\tilde b}}{\varrho^{\tilde b}}\|2\tilde \Omega_mF^{20}_{nn}+{\rm i}P^{20}_{nn}\|_{D(r-2\varrho),\Cal O}\nonumber\\
&\leq& e^{-K\varrho}\cdot \frac{(1+e)^{\tilde b}}{\varrho^{\tilde b}}(2\|\tilde \Omega_n\|_{D(r-2\varrho),\Cal O}\cdot \|F^{20}_{nn}\|_{D(r-2\varrho),\Cal O}+\|P^{20}_{nn}\|_{D(r-2\varrho),\Cal O})\nonumber\\
&\leq& e^{-K\varrho}\cdot \frac{(1+e)^{\tilde b}}{\varrho^{\tilde b}}(2|n|\delta_0\cdot \frac{16(1+e)^{\tilde b}}{|n|\varrho^{\tilde b}}\|P^{20}_{nn}\|_{D(r),\Cal O}+\|P^{20}_{nn}\|_{D(r),\Cal O})\nonumber\\
&\leq& e^{-K\varrho}\cdot \frac{2(1+e)^{2\tilde b}}{\varrho^{2\tilde b}}\|P^{20}_{nn}\|_{D(r),\Cal O},
\end{eqnarray}
the estimates of $F^{02}_{nn}$ and the error term $\hat R^{02}_{nn}$ can be similarly estimated with the same results.\\
To get the norm estimates of the vector field $X_{F_{2,>}}, X_{\hat R_{2,>}}$, we can repeat the above proof process of $X_{F_{2,<}}, X_{\hat R_{2,<}}$ and obtain
\begin{eqnarray}
\|X_{F_{2,>}}\|_{s,a,\rho,D(r-3\rho,s),\Cal O}\leq \frac{16(1+e)^{\tilde b}}{\varrho^{\tilde b+1}}\|X_{R_{2,>}}\|_{s,\bar a, \rho,D(r,s),\Cal O},\label{XF2}\\
\|X_{\hat R_{2,>}}\|_{s,\bar a,\rho,D(r-5\rho,s),\Cal O}\leq e^{-K\varrho}\cdot \frac{2(1+e)^{2\tilde b}}{\varrho^{2\tilde b+1}}\|X_{R_{2,>}}\|_{s,\bar a,\rho,D(r,s),\Cal O}.\label{xer2}
\end{eqnarray}
Hence, the norm of the vector field $X_{F_2}$ is obtained with the estimates in $(\ref{XF1})$, $(\ref{XF2})$
\begin{eqnarray}{\label{xf2}}
&&\|X_{F_{2}}\|_{s,a,\rho,D(r-3\rho,s),\Cal O_+}\leq \|X_{F_{2,<}}\|_{s,a,\rho,D(r-3\rho,s),\Cal O_+}+\|X_{F_{2,>}}\|_{s,a,\rho,D(r-3\rho,s),\Cal O_+}\nonumber\\
&&\leq\frac{cE^2K^{2\tau+2}}{\gamma_0^2\varrho^{\tilde{b}+1}}\cdot e^{\frac{8 E^2\delta_0(\gamma_0-\gamma)Kr}{\gamma^2}}\|X_{R_{2,<}}\|_{s,\bar a,\rho,D(r,s),\Cal O}+
\frac{16(1+e)^{\tilde b}}{\varrho^{\tilde b+1}}\|X_{R_{2,>}}\|_{s,\bar a, \rho,D(r,s),\Cal O}\nonumber\\
&&\leq \frac{cE^2K^{2\tau+2}}{\gamma_0^2\varrho^{\tilde{b}+1}}\cdot e^{\frac{8 E^2\delta_0(\gamma_0-\gamma)Kr}{\gamma^2}}\|X_{R_{2}}\|_{s,\bar a,\rho,D(r,s),\Cal O}
\end{eqnarray}
and the vector norm field of the error term $\hat{R}_{2}$ is obtained with the estimates in $(\ref{xer1})$, $(\ref{xer2})$
\begin{eqnarray}{\label{XER2}}
&&\|X_{\hat R_{2}}\|_{s,\bar a,\rho,D(r-5\rho,s),\Cal O_+}\nonumber\\
&\leq& \|X_{\hat R_{2,<}}\|_{s,\bar a,\rho,D(r-5\rho,s),\Cal O_+}+\|X_{\hat R_{2,>}}\|_{s,\bar a,\rho,D(r-5\rho,s),\Cal O_+}\nonumber\\
&\leq& \frac{cE^2K^{2\tau+2}\delta_0}{\gamma_0\varrho^{2\tilde b+1}}e^{-K\varrho}\cdot e^{\frac{16 E^2\delta_0(\gamma_0-\gamma)Kr}{\gamma^2}}\|X_{R_{2,<}}\|_{\!{}_{s,\bar a,\rho,D(r,s),\Cal O}}\nonumber\\
&+&e^{-K\varrho}\cdot \frac{2(1+e)^{2\tilde b}}{\varrho^{2\tilde b+1}}\|X_{R_{2,>}}\|_{\!{}_{s,\bar a,\rho,D(r,s),\Cal O}}\nonumber\\
&\leq& \frac{cE^2K^{2\tau+2}\delta_0}{\gamma_0\varrho^{2\tilde b+1}}e^{-K\varrho}\cdot e^{\frac{16 E^2\delta_0(\gamma_0-\gamma)Kr}{\gamma^2}}\|X_{R_{2}}\|_{\!{}_{s,\bar a,\rho,D(r,s),\Cal O}}.
\end{eqnarray}
Similarly we can get the estimates of $X_{F_0},X_{F_1}$ and the error term $X_{\hat R_{1}}$
\begin{eqnarray}
&&\|X_{F_0}\|_{\!{}_{s,a,\rho,D(r-3\varrho,s),\Cal O_+}}\leq
\frac{cEK^{2\tau+1}}{\gamma^2\varrho^{\tilde{b}+1}}
\|X_{R_{0}}\|_{\!{}_{s,\bar a,\rho,D(r,s),\Cal O}},\label{XF0}\\
&&\|X_{F_1}\|_{\!{}_{s,a,\rho,D(r-3\varrho,s),\Cal O_+}}\leq
\frac{cE^2K^{2\tau+2}}{\gamma_0^2\varrho^{\tilde{b}+1}}e^{\frac{4E^2\delta_0(\gamma_0-\gamma)Kr}{\gamma^2}}
\|X_{R_{1}}\|_{\!{}_{s,\bar a,\rho,D(r,s),\Cal O}},\\
&&\|X_{\hat R_{1}}\|_{\!{}_{s,\bar a, \rho,D(r-5\varrho,s),\Cal O_+}}\leq
\frac{cE^2K^{2\tau+2}\delta_0}{\gamma_0\varrho^{2\tilde b+1}}e^{-K\varrho} e^{\frac{8 E^2\delta_0(\gamma_0-\gamma)Kr}{\gamma^2}}\|X_{R_{1}}\|_{\!{}_{s,\bar a,\rho,D(r,s),\Cal O}},\label{XhR1}
\end{eqnarray}
then finally we obtain
\begin{eqnarray}
\|X_{F}\|_{\!{}_{s,a,\rho,D(r-3\varrho,s),\Cal O_+}}
&\leq&\frac{cE^2K^{2\tau+2}}{\gamma_0^2\varrho^{\tilde{b}+1}}\cdot e^{\frac{8 E^2\delta_0(\gamma_0-\gamma)Kr}{\gamma^2}}\|X_{R}\|_{\!{}_{s,\bar a,\rho,D(r,s),\Cal O}}, \label{XF}\\
\|X_{\hat R}\|_{\!{}_{s,\bar a,\rho,D(r-5\varrho,s),\Cal O_+}}
&\leq&\frac{cE^2K^{2\tau+2}\delta_0}{\gamma_0\varrho^{2\tilde b+1}}e^{-K\varrho} e^{\frac{16 E^2\delta_0(\gamma_0-\gamma)Kr}{\gamma^2}}\|X_{R}\|_{\!{}_{s,\bar a,\rho,D(r,s),\Cal O}}.\label{XER}
\end{eqnarray}

\subsection{Estimation for the new normal form}\label{no4.3}

  The map $\phi_F^1$ defined above transforms $H$
into $H_+=N_+ +  P_+$. As mentioned in Remark\ref{RK}, we return to the original coordinates, here the generalized normal form $N_+$ is
\begin{eqnarray*}
N_+= N+\hat{N}\nonumber,\quad \hat{N}&=&\la\hat{\omega},I\ra+\sum_{n\in \Z_1\atop |n|\leq EK}\hat{\Omega}_nw_n\bar w_n\\
&+&\sum_{n\in \Z_1\atop |n|\leq EK}(\hat{a_{n(-n)}}w_n\bar w_{(-n)}+\hat{a_{(-n)n}}w_{(-n)}\bar w_n)
+\sum_{n\in \Z_1\atop |n|> EK}\hat{\Omega}_nw_n\bar w_n,\nonumber\\
\hat{\omega}=[R_{l00}],(|l|=1),\quad \hat{\Omega}_n&=&P^{11}_{nn}-\la\partial_{\tilde{\theta}}\Omega_{n},\partial_I F_0\ra=P^{11}_{nn}-\la\partial_{\tilde{\theta}}\tilde{\Omega}_{n},\partial_I F_0\ra,
\label{N_+form}
 \end{eqnarray*}
 We rewrite $N_+$ as follows:
 $$
N_+=\la\bar \omega,\bar I\ra+\la\tilde\omega_+,I\ra+\sum_{n\in \Z_1}{\Omega}_n^+w_n\bar w_n+\sum_{|n|\leq EK}\langle A^+_{|n|}z_{|n|}, \bar z_{|n|}\rangle,$$
where $$A_{|n|}^+= (\overline{A_{|n|}^+})^T=\left(\begin{array}{ccc}a_{nn}^+(\sigma)&a_{n(-n)}^+(\sigma)\\
a_{(-n)n}^+(\sigma)&a_{(-n)(-n)}^+(\sigma)
\end{array}\right), z_{|n|}=\left(\begin{array}{c}w_n\\
w_{(-n)}
\end{array}\right),\bar z_{|n|}=\left(\begin{array}{c}\bar w_n\\
\bar w_{(-n)}
\end{array}\right).$$
It is obvious that
 \[|\hat{\omega}|_{\Cal
O_+}\leq c\|X_R\|_{s,\bar a,\rho,D(r,s),\Cal O}.\]
Then we estimate $\hat{\Omega}=(\hat{\Omega}_n:n\in\Z_1)$
\begin{eqnarray}
\|\la \partial_{\tilde \theta}\tilde{\Omega}_n,\partial_I F_0\ra\|_{D(r-\varrho),\Cal O}&\leq& \|\tilde{\Omega}_n\|_{r,2\tau+2,\Cal O}\cdot\|X_{F_0}\|_{s,a,\rho,D(r-\varrho,s),\Cal O}\nonumber\\
&\leq& \frac{c\delta_0(\gamma_0-\gamma)nEK^{2\tau+1}}{\gamma^2\varrho^{\tilde{b}+1}}
\|X_{R_{0}}\|_{s,\bar a,\rho,D(r,s),\Cal O},\label{4.40}
\end{eqnarray}
associated with  $\|P^{11}_{nn}\|_{D(r-\varrho),\Cal O}\leq n\cdot\|X_R\|_{s,\bar a,\rho,D(r,s),\Cal O}$, we have
$$|\hat{\Omega}|_{-1,D(r-\varrho),\Cal O}\leq \frac{c\delta_0(\gamma_0-\gamma)EK^{2\tau+1}}{\gamma^2\varrho^{\tilde{b}+1}}
\|X_{R_{0}}\|_{s,\bar a,\rho,D(r,s),\Cal O}.$$
It follows that
 \begin{eqnarray}\label{N+}
 \|X_{\hat N}\|_{s,\bar a,\rho,D(r-2\varrho,s),\Cal O_+}\leq \frac{c\delta_0(\gamma_0-\gamma)EK^{2\tau+1}}{\gamma^2\varrho^{\tilde{b}+1}}
\|X_{R}\|_{s,\bar a,\rho,D(r,s),\Cal O}.
 \end{eqnarray}
\subsection{Estimation for the new perturbation}\label{no4.4}

Since $ P_+=\hat{R}+\int_0^1 \{(1-t)(\hat N+\hat R)+tR,F\}\circ
\phi_F^{t}dt+(P-R)\circ
\phi^1_F$, we set $R(t)=(1-t)(\hat N+\hat R)+tR$, hence
$$
X_{P_+}=X_{\hat R}+\int_0^1 (\phi_F^{t})^*X_{\{R(t),F\}} dt
+(\phi^1_F)^*X_{(P-R)}.
$$
It is obvious that the vector norm of the error term $\hat R$ has been given in $(\ref{XER})$
\begin{eqnarray*}
\|X_{\hat R}\|_{s,\bar a,\rho,D(r-5\varrho,s),\Cal O_+}
&\leq&\frac{cE^2K^{2\tau+2}\delta_0}{\gamma_0\varrho^{2\tilde b+1}}e^{-K\varrho}\cdot e^{\frac{16 E^2\delta_0(\gamma_0-\gamma)Kr}{\gamma^2}}\|X_{R}\|_{s,\bar a,\rho,D(r,s),\Cal O}.
\end{eqnarray*}
We rewrite $P-R=P_{(1)}+P_{(2)}+P_{(3)}$ as
\begin{eqnarray}
P_{(1)}&=&\sum_{\alpha,\beta}P_{\alpha\beta}(\theta)w^\alpha\bar{w}^\beta,\label{P(1)}\\
P_{(2)}&=&\sum_{|k|>K,l,\alpha,\beta\atop 2|l|+|\alpha+\beta|\leq2}
P_{kl\alpha\beta}I^lw^{\alpha}\bar{w}^{\beta}\kth,\label{P(2)}\\
P_{(3)}&=&\sum_{l,\alpha,\beta\atop 2|l|+|\alpha+\beta|>2}
P_{l\alpha\beta}(\theta,\sigma)I^lw^{\alpha}\bar{w}^{\beta},\label{P(3)}
\end{eqnarray}
with $\alpha,\beta\in \{\alpha,\beta\in \N^{\Z_1},|\alpha+\beta|=\alpha_n+\beta_n=1,\forall |n|>EK\}$ in $P_{(1)}$, $\alpha,\beta\in \{\alpha,\beta\in \N^{\Z_1},\alpha_n+\beta_n\in 2\N, \forall |n|>EK\}$ in $P_{(2)}$ and $P_{(3)}$.
Recalling the decay estimates in $(\ref{decay1})$, it is clear that
$\|P_{\alpha\beta}\|_{D(r),\Cal O}\leq c\varepsilon e^{-n\bar \rho},\ \alpha+\beta=e_n, |n|>EK$ in $P_{(1)}$, and by these conditions, we have
\begin{eqnarray}
\|X_{P_{(1)}}\|_{\eta s,\bar a,\rho_+,D(r,\eta s),\Cal O_+}&\leq& c\eta^{-1}e^{-\frac{EK\bar \rho}{2}}\|X_{R}\|_{s,\bar a, \rho,D(r,s),\Cal O},\label{XP(1)}\\
\|X_{P_{(2)}}\|_{\eta s,\bar a,\rho_+,D(r-5\varrho,4\eta s),\Cal O_+}&\leq& c\eta^{-1}e^{-K\varrho}\|X_{R}\|_{s,\bar a, \rho,D(r,s),\Cal O},\label{XP(2)}\\
\|X_{P_{(3)}}\|_{\eta s,\bar a,\rho_+,D(r-5\varrho,4\eta s),\Cal O_+}&\leq& c\eta\|X_{R}\|_{s,\bar a, \rho,D(r,s),\Cal O},\label{XP(3)}
\end{eqnarray}
then
\begin{eqnarray*}
\|X_P-X_R\|_{\eta s,\bar s,\rho_+,D(r-5\varrho,4\eta s),\Cal O_+}
&\leq& c(\eta^{-1}e^{-\frac{EK\bar \rho}{2}}+ \eta^{-1}e^{-K\varrho}+\eta)\|X_R\|_{s,\bar a,\rho,D(r,s),\Cal O}.
\end{eqnarray*}
According to  $(\ref{XF})$,
\begin{eqnarray}\label{DXF}
&&\|DX_F\|_{s,a,\rho,D(r-4\varrho,s),\Cal O_+},\ \ \|DX_F\|_{s,\bar a,\rho,D(r-4\varrho,s),\Cal O_+}\nonumber\\
&&\leq \frac{cE^2K^{2\tau+2}}{\gamma_0^2\varrho^{\tilde{b}+2}}\cdot e^{\frac{8 E^2\delta_0(\gamma_0-\gamma)Kr}{\gamma^2}}\|X_{R}\|_{s,\bar a,\rho,D(r,s),\Cal O}.
\end{eqnarray}
We assume that
\begin{eqnarray}\label{XP}
\|X_P\|_{s,\bar a,\rho,D(r,s),\Cal O}\leq \varepsilon \overset{(\ref{B})}{\leq} \frac{\eta^2}{\Cal{B}_\varrho (K\varrho)^{2\tau+2}}e^{-\frac{8 E^2\delta_0(\gamma_0-\gamma)Kr}{\gamma^2}},
\end{eqnarray}
this inequality will be verified in the section 5, where $\Cal{B}_\varrho=cE^4\varrho^{-10(\tilde{b}+\tau+1)}$($c=\gamma_0^{-4}c(\tilde{b},\tau)$) is a sufficiently large constant with a fixed $\gamma_0 >0 $, then
$$\|X_F\|_{s, a,\rho,D(r-3\varrho,s),\Cal O_+}, \|DX_F\|_{s,a,\rho,D(r-4\varrho,s),\Cal O_+},\|DX_F\|_{s,\bar a,\rho,D(r-4\varrho,s),\Cal O_+}\overset{(\ref{PHI})}{\leq} \Cal{B}_\varrho^{\frac{1}{2}}\varepsilon^{1-\beta'},$$
with some constant $0<\beta'< 1$.Then the follow $\phi^t_F$ of the vector field $X_F$ exists on $D(r-5\varrho,\frac{s}{2})$ for $-1\leq t\leq 1$, and takes this domain into $D(r-4\varrho,s)$, we obtain
\begin{eqnarray}
\|\phi^t_F-id\|_{s,a,\rho,D(r-5\varrho,\frac{s}{2}),\Cal O_+} &\leq& c\|X_F\|_{s,a,\rho,D(r-4\varrho,s),\Cal O_+}\label{4.5},\\
\|D\phi^t_F-I\|_{s,a,\rho,D(r-6\varrho,\frac{s}{4}),\Cal O_+} &\leq& c\|DX_F\|_{s,a,\rho,D(r-4\varrho,s),\Cal O_+}\label{4.6},\\
\|D\phi^t_F-I\|_{s,\bar a,\rho,D(r-6\varrho,\frac{s}{4}),\Cal O_+} &\leq&  c\|DX_F\|_{s,\bar a,\rho,D(r-4\varrho,s),\Cal O_+}\label{4.7}.
\end{eqnarray}
Also we have that for any vector field $Y$,
$$\|(D\phi^t_F)^*Y\|_{\eta s,\bar a,\rho,D(r-7\varrho,\eta s), \Cal O_+}\leq c\|Y\|_{\eta s,\bar a,\rho,D(r-5\varrho,4\eta s), \Cal O_+},$$
and with the estimates  $(\ref{N+})$, $(\ref{XER})$, we get
\begin{eqnarray*}
&&\|X_{R(t)}\|_{s,\bar a, \rho,D(r-5\varrho,s),\Cal O_+}\\
&\leq& \|X_{\hat N}\|_{s,\bar a, \rho, D(r-5\varrho,s),\Cal O_+}+\|X_{\hat{R}}\|_{s,\bar a, \rho,D(r-5\varrho,s),\Cal O_+}\nonumber\\
&\leq& \frac{cE^2K^{2\tau+2}\delta_0}{\gamma_0^2\varrho^{2\tilde b+1}}e^{-K\varrho}\cdot e^{\frac{16 E^2\delta_0(\gamma_0-\gamma)Kr}{\gamma^2}}\|X_{R}\|_{s,\bar a,\rho,D(r,s),\Cal O}.
\end{eqnarray*}
Moreover, we have
\begin{eqnarray*}
&&\|[X_{R(t)},X_F]\|_{\eta s,\bar a,\rho,D(r-6\varrho,\frac{s}{2}),\Cal O_+}\\
&\leq&
\|DX_{R(t)}\|_{s,\bar a,\rho,D(r-6\varrho,\frac{s}{2}),\Cal O_+}\cdot\|X_F\|_{s,a,\rho,D(r-6\varrho,\frac{s}{2}),\Cal O_+}\nonumber\\
&+&\|DX_F\|_{s,a,\rho,D(r-6\varrho,\frac{s}{2}),\Cal O_+}\cdot\|X_{R(t)}\|_{s,\bar a,\rho,D(r-6\varrho,\frac{s}{2}),\Cal O_+}\nonumber\\
&\leq& \frac{cE^4K^{4\tau+4}\delta_0}{\eta^2\gamma_0^4\varrho^{3\tilde b+3}}e^{-K\varrho}\cdot e^{\frac{24 E^2\delta_0(\gamma_0-\gamma)Kr}{\gamma^2}}(\|X_P\|_{s,\bar a,\rho,D(r,s),\Cal O})^2,
\end{eqnarray*}
together with the estimates of $\hat R$ and $X_P-X_R$, we finally arrive at the estimate
\begin{eqnarray}\label{XP+}
\|X_{P_+}\|_{\eta s,\bar a,\rho_+,D(r-7\varrho,\eta s),\Cal O_+}&\leq&
\frac{1}{5}(\frac{cE^4K^{4\tau+4}\delta_0}{\eta^2\gamma_0^4\varrho^{3\tilde b+3}}e^{-K\varrho}\cdot e^{\frac{24 E^2\delta_0(\gamma_0-\gamma)Kr}{\gamma^2}}\|X_P\|_{s,\bar a,\rho,D(r,s),\Cal O}\nonumber\\
&+&\frac{cE^2K^{2\tau+2}\delta_0}{\gamma_0^2\varrho^{2\tilde b+1}}e^{-K\varrho}\cdot e^{\frac{16 E^2\delta_0(\gamma_0-\gamma)Kr}{\gamma^2}}+c\eta^{-1}e^{-\frac{EK\bar \rho}{2}}\nonumber\\
&+&c\eta^{-1}e^{-K\varrho}+c\eta)\|X_P\|_{ s,\bar a,\rho,D(r,s),\Cal O}.
\end{eqnarray}
This is the bound for the new perturbation.
\subsection{Verification of $(A5)$ after one KAM iteration}
We will verify the new perturbation $P_+$ with the special structure and decay properties in $(A5)$ with $E, K, \varepsilon_+$ in place of $E_-, K_-, \varepsilon$. For simplicity we denote $D(r_+,s_+)=D_+$ with $ s_+=\eta s$ defined in Section 5 in the following calculations. Since
\begin{eqnarray*}
P_+&=&\hat R +P-R+\{P,F\}+\frac{1}{2!}\{\{N,F\},F\}+\frac{1}{2!}\{\{P,F\},F\}\\
&+& \cdots +\frac{1}{n!}\{\cdots \{N,\underbrace{F\}\cdots, F}_n\}+\frac{1}{n!}\{\cdots \{P,\underbrace{F\}\cdots, F}_n\}+\cdots,
\end{eqnarray*}
where $\hat R=\hat R_1+\hat R_2$ is the error term with the formula
\begin{eqnarray*}
\hat R_1&=&\sum_{|n|\leq EK} \hat R^{10}_{n}w_n+\hat R^{01}_{n}\bar w_n,\ \ \ \hat R_2=\sum_{|n|,|m|\leq EK}(\hat R^{20}_{nm}w_nw_m+\hat{R}^{02}_{nm}\bar w_n\bar w_m)+\\
&&\sum_{|n|\neq |m|\atop |n|,|m|\leq EK}\hat R^{11}_{nm}w_n\bar w_m
+\sum_{|n|>EK}(\hat R^{20}_{nn}w_nw_n+\hat{R}^{02}_{nn}\bar w_n\bar w_n),
\end{eqnarray*}
and $P-R=P_{(1)}+P_{(2)}+P_{(3)}$ defined in the $(\ref{P(1)}),(\ref{P(2)}),(\ref{P(3)})$, so it is obvious that the $\hat R, P-R$  both have the special structure in $(A5)$. Besides, by $(\ref{decay1})$,$(\ref{XhR1})$,$(\ref{XP})$,$(\ref{e1})-(\ref{B})$, we have when $|n|\leq EK$,
$$\|\hat R^{10}_{n}\|_{D_+,\Cal O_+},\ \|\hat R^{10}_{n}\|_{D_+,\Cal O_+} \leq ce^{-|n|\rho}B^{\frac{1}{2}}_\varrho\varepsilon^{2-\frac{8}{5}\beta'}\leq
 c\varepsilon_+e^{-|n|\rho_+}.$$
Using  $(\ref{XP(1)})$ and $(\ref{XP-R1})$, the decay property of $P-R$ can be similarly obtained.
In the following, we will consider the term $\{P,F\}$ with $F=F_0+F_1+F_2$ rewritten as
\begin{eqnarray*}
F=F_0(\theta,I)+\sum_{\alpha,\beta}F^1_{\alpha\beta}(\theta)w^{ \alpha}\bar w^{\beta}+\sum_{\alpha,\beta}F^2_{\alpha\beta}(\theta)w^{ \alpha}\bar w^{\beta},
\end{eqnarray*}
where  $\alpha,\beta\in \{\alpha,\beta\in \N^{\Z_1},|\alpha+\beta|=\alpha_n+\beta_n=1, |n|\leq EK\}$ in $F^1$,  $\alpha,\beta\in \{\alpha,\beta\in \N^{\Z_1}, |\alpha+\beta|=\alpha_n+\beta_n=2,  |n|>E_-K_-\}$ in $F^2$,
 $F^2_{e_ne_m}=F^2_{e_me_n}=0$  with $|n=m| \leq E_-K_-$. Then we calculate $\{P,F\}=\{P^1,F\}+\{P^2,F\}+\{P^3,F\}$,
\begin{eqnarray*}
\{P^1,F\}
&=&\sum_{\alpha,\beta}\frac{\partial P^1_{\alpha\beta}}{\partial I}\cdot\frac{\partial F_{0}}{\partial \tilde{\theta}}w^{\alpha}\bar w^{\beta}-\sum_{\alpha,\beta}\frac{\partial P^1_{\alpha\beta}}{\partial \tilde\theta}\cdot\frac{\partial F_{0}}{\partial I}w^{\alpha}\bar w^{\beta}\\
&+&\sum_{\alpha,\beta\atop \alpha',\beta'}\frac{\partial P^1_{\alpha\beta}}{\partial I}\cdot\frac{\partial F^1_{\alpha'\beta'}}{\partial \tilde{\theta}}w^{\alpha}\bar w^{\beta}w^{\alpha'}\bar w^{\beta'}
+\sum_{\alpha,\beta\atop \tilde\alpha,\tilde\beta}\frac{\partial P^1_{\alpha\beta}}{\partial I}\cdot\frac{\partial F^2_{\tilde\alpha\tilde\beta}}{\partial \tilde{\theta}}w^{\alpha}\bar w^{\beta}w^{\tilde\alpha}\bar w^{\tilde\beta}\\
&+&{\rm i}\sum_{\alpha,\beta \atop n\leq EK}(\alpha_nP^1_{\alpha\beta}F^1_{0e_n}w^{\alpha-e_n}\bar w^{\beta}-\beta_nP^1_{\alpha\beta}F^1_{e_n0}w^{\alpha}\bar w^{\beta-e_n})\\
&+&{\rm i}\sum_{n,\alpha,\beta \atop \tilde\alpha,\tilde\beta}P^1_{\alpha\beta}F^2_{\tilde\alpha\tilde\beta}(\alpha_n\tilde\beta_nw^{\alpha-e_n}\bar w^{\beta}w^{\tilde\alpha}\bar w^{\tilde\beta-e_n}-\beta_n\tilde\alpha_nw^{\alpha}\bar w^{\beta-e_n}w^{\tilde\alpha-e_n}\bar w^{\tilde\beta}),
\end{eqnarray*}
where $\alpha,\beta\in \{\alpha_n+\beta_n\in 2\N,\forall |n|>E_-K_-\}$,
$\alpha',\beta'\in \{|\alpha'+\beta'|=\alpha'_n+\beta'_n=1, |n|\leq EK\}$, $\tilde\alpha,\tilde\beta\in \{|\tilde\alpha+\tilde\beta|=\tilde\alpha_n+\tilde\beta_n=2,\forall |n|>E_-K_-\}$. So the exponent of $ w^{\alpha}\bar w^{\beta}w^{\alpha'}\bar w^{\beta'}$ satisfies $ \alpha+\alpha',\beta+\beta'\in\{\alpha_n+\alpha'_n+\beta_n+\beta'_n=\alpha_n+\beta_n\in 2\N,\forall |n|>EK\}$, the exponent of $w^{\alpha}\bar w^{\beta}w^{\tilde\alpha}\bar w^{\tilde\beta}$ satisfies $ \alpha+\tilde\alpha,\beta+\tilde\beta\in\{\alpha_n+\tilde\alpha_n+\beta_n+\tilde\beta_n\in 2\N,\forall |n|>EK\}$, the exponents of $w^{\alpha-e_n}\bar w^{\beta}w^{\tilde\alpha}\bar w^{\tilde\beta-e_n}$, $w^{\alpha}\bar w^{\beta-e_n}w^{\tilde\alpha-e_n}\bar w^{\tilde\beta}$ satisfy $ \alpha-e_n+\tilde\alpha,\beta+\tilde\beta-e_n\in\{\alpha_m+\alpha_m+\tilde\beta_m+\tilde\beta_m-2\delta_{nm}\in 2\N,\forall |m|>EK\}$ for any $|n|>EK$.
\begin{eqnarray*}
\{P^2,F\}&=&-\sum_{\alpha,\beta}\frac{\partial P^2_{\alpha\beta}}{\partial \tilde\theta}\cdot\frac{\partial F_{0}}{\partial I}w^{\alpha}\bar w^{\beta}+{\rm i}\sum_{|n|\leq EK}(P^2_{e_n0}F^1_{0e_n}-P^2_{0e_n}F^1_{e_n0})\\
&+&{\rm i}\sum_{n,\tilde\alpha,\tilde\beta}(\tilde\beta_nP^2_{e_n0}F^2_{\tilde\alpha\tilde\beta}w^{\tilde\alpha}\bar w^{\tilde\beta-e_n}
-\tilde\alpha_nP^2_{0e_n}F^2_{\tilde\alpha\tilde\beta}w^{\tilde\alpha-e_n}\bar w^{\tilde\beta}),
\end{eqnarray*}
where  $\alpha,\beta\in \{|\alpha+\beta|=\alpha_n+\beta_n=1,\forall |n|>E_-K_-\}$, $\tilde\alpha,\tilde\beta\in \{|\tilde\alpha+\tilde\beta|=\tilde\alpha_n+\tilde\beta_n=2,|n|>E_-K_-\}$. So the exponents of $w^{\tilde\alpha}\bar w^{\tilde\beta-e_n}, w^{\tilde\alpha-e_n}\bar w^{\tilde\beta}$ are contained in   $\{\tilde\alpha_m+\tilde\beta_m-\delta_{nm}=0 \ or \ 1,\forall |m|>EK\}$ for any $|n|>EK$.
\begin{eqnarray*}
\{P^3,F\}&=&\sum_{\alpha,\beta}\frac{\partial P^3_{\alpha\beta}}{\partial I}\cdot\frac{\partial F_{0}}{\partial \tilde{\theta}}w^{\alpha}\bar w^{\beta}-\sum_{\alpha,\beta}\frac{\partial P^3_{\alpha\beta}}{\partial \tilde\theta}\cdot\frac{\partial F_{0}}{\partial I}w^{\alpha}\bar w^{\beta}\\
&+&\sum_{\alpha,\beta\atop \alpha',\beta'}\frac{\partial P^3_{\alpha\beta}}{\partial I}\cdot\frac{\partial F^1_{\alpha'\beta'}}{\partial \tilde{\theta}}w^{\alpha}\bar w^{\beta}w^{\alpha'}\bar w^{\beta'}
+\sum_{\alpha,\beta\atop \tilde\alpha,\tilde\beta}\frac{\partial P^3_{\alpha\beta}}{\partial I}\cdot\frac{\partial F^2_{\tilde\alpha\tilde\beta}}{\partial \tilde{\theta}}w^{\alpha}\bar w^{\beta}w^{\tilde\alpha}\bar w^{\tilde\beta}\\
&+&{\rm i}\sum_{\alpha,\beta \atop |n|\leq EK}(\alpha_nP^3_{\alpha\beta}F^1_{0e_n}w^{\alpha-e_n}\bar w^{\beta}-\beta_nP^3_{\alpha\beta}F^1_{e_n0}w^{\alpha}\bar w^{\beta-e_n})\\
&+&{\rm i}\sum_{n, \alpha,\beta\atop \tilde\alpha,\tilde\beta}P^3_{\alpha\beta}F^2_{\tilde\alpha\tilde\beta}(\alpha_n\tilde\beta_nw^{\alpha-e_n}\bar w^{\beta}w^{\tilde\alpha}\bar w^{\tilde\beta-e_n}-\beta_n\tilde\alpha_nw^{\alpha}\bar w^{\beta-e_n}w^{\tilde\alpha-e_n}\bar w^{\tilde\beta}),
\end{eqnarray*}
where $\alpha,\beta\in \{\alpha_n+\beta_n=0, \forall |n|>E_-K_-\}$,
$\alpha',\beta'\in \{|\alpha'+\beta'|=\alpha'_n+\beta'_n=1, |n|\leq EK\}$, $\tilde\alpha,\tilde\beta\in \{|\tilde\alpha+\tilde\beta|=\tilde\alpha_n+\tilde\beta_n=2,\forall |n|>EK\}$. So $w^{\alpha}\bar w^{\beta}w^{\alpha'}\bar w^{\beta'}$, $w^{\alpha-e_n}\bar w^{\beta}$,$ w^{\alpha}\bar w^{\beta-e_n}$, $w^{\alpha-e_n}\bar w^{\beta}w^{\tilde\alpha}\bar w^{\tilde\beta-e_n}$, $w^{\alpha}\bar w^{\beta-e_n}w^{\tilde\alpha-e_n}\bar w^{\tilde\beta}$ disappear with $|n|>EK$. The exponent of $w^{\alpha}\bar w^{\beta}w^{\tilde\alpha}\bar w^{\tilde\beta}$ satisfies $\{\alpha_n+\beta_n+\tilde\alpha_n+\tilde\beta_n=2\in 2\N, \forall |n|>EK\}$.
 When $n>EK, |\alpha+\beta|=\alpha_n+\beta_n=1, \tilde\alpha+\tilde\beta=2e_n$, by $(\ref{decay1})$,$(\ref{XF0}),(\ref{F20})$,$(\ref{XP})$,$(\ref{e1})-(\ref{B})$, we have
\begin{eqnarray*}
&&\|\frac{\partial P^2_{\alpha\beta}}{\partial \tilde\theta}\cdot \frac{\partial F_{0}}{\partial I}\|_{D_+, \Cal O_+}\leq c e^{-|n|\bar\rho}\Cal B^{\frac{1}{2}}\varepsilon^{2-\frac{\beta'}{5}}\leq c\varepsilon_+e^{-|n|\bar \rho},\\
&&\|P^2_{e_n0}F^2_{\tilde\alpha'\tilde \beta'}\|_{D_+,\Cal O_+}, \|P^2_{0e_n}F^2_{\tilde\alpha'\tilde \beta'}\|_{D_+,\Cal O_+}\leq c e^{-|n|\bar \rho}\Cal B^{\frac{1}{2}}\varepsilon^{2-\frac{\beta'}{5}}\leq c\varepsilon_+e^{-|n|\bar \rho},
\end{eqnarray*}
together with the decay estimates of $ \hat R, P-R, \{P,F\}$,  the decay property of $P_+$  in assumption $(A5)$ has been finally verified.

\subsection{Verification of $(A6)$ after one KAM iteration}
 In the following, we have to check that the new perturbation $P_+$ satisfies $(A6)$
with $\varepsilon_+$ in place of $\varepsilon$, namely,
for $n\in \Z_1$, we need to verify
\begin{eqnarray}
&&\|\lim_{n\rightarrow\infty}\frac{1}{|n|}\sum_{\upsilon=\pm}\frac{\partial^2P_+}{\partial w^\upsilon_n\partial w^\upsilon_n}\|_{\!{}_{D(r_+,s_+),\Cal O_+}}\leq \varepsilon_+,\label{P1}\\
&&\|\frac{1}{|n|}\sum_{\upsilon=\pm}\frac{\partial^2P_+}{\partial w^\upsilon_n\partial w^\upsilon_n}-\lim_{n\rightarrow\infty}\frac{1}{|n|}\sum_{\upsilon=\pm}\frac{\partial^2P_+}{\partial w^\upsilon_n\partial w^\upsilon_n}\|_{\!{}_{D(r_+,s_+),\Cal O_+}}\leq \frac{\varepsilon_+}{|n|}.\label{P2}
\end{eqnarray}
According to the form of $P_+$ in the above subsection 4.6, it is sufficient for us to consider the three main terms $\hat R, P-R, \{P,F\}$. Due to $\hat R=\hat R_{1}+\hat R_{2}$, it is sufficient to prove that $\hat R_{2}$ of order 2 in $w,w$ satisfies $(A6)$. Similarly for the term $P-R=P_{(1)}+P_{(2)}+P_{(3)}$, it is sufficient to show that $P_{(2)},P_{(3)}$ with $\alpha,\beta\in \{ \alpha_n+\beta_n\in 2\N,\forall |n|> EK\}$ satisfy $(A6)$. Besides, we need to prove the property $(A6)$ of the term $\{P,F\}$.
Firstly, we consider the term $\hat R_{2}=\hat R_{2,<}+\hat R_{2,>}$ of order 2 in $w,w$ with the form
\begin{eqnarray*}
\hat R_{2,<}&=&\sum_{|n|,|m|\leq EK}(\hat R^{20}_{nm}w_nw_m+\hat{R}^{02}_{nm}\bar w_n\bar w_m)+\sum_{|n|,|m|\leq EK}\hat R^{11}_{nm}w_n\bar w_m,\\
\hat R_{2,>}&=&\sum_{|n|>EK}(\hat R^{20}_{nn}w_nw_n+\hat{R}^{02}_{nn}\bar w_n\bar w_n),
\end{eqnarray*}
it is clear that
\beq \sum_{\upsilon=\pm}\frac{\partial^2 \hat R_{2}}{\partial w^\upsilon_n\partial w^\upsilon_n}=
\left \{\begin{array}{lr}2\hat R^{20}_{nm}+2\hat R^{02}_{nm},  &1\leq |n=m|\leq EK,\\
2\hat R^{20}_{nn}+2\hat R^{02}_{nn}, &|n|>EK,\nonumber\\
 \end{array}\right.\eeq
by the estimates $(\ref{R11}),(\ref{er2})$, $(\ref{decay11})$, we have
\begin{eqnarray*}
&&\|\lim_{n\rightarrow \infty}\frac{1}{|n|}\sum_{\upsilon=\pm}\frac{\partial^2 \hat R_{2}}{\partial w^\upsilon_n\partial w^\upsilon_n}\|_{\!{}_{D_+,\Cal O_+}}\\
&\leq&\lim_{n\rightarrow \infty} \frac{1}{|n|}\frac{cE^2K^{2\tau+2}\delta_0}{\gamma_0^2\varrho^{2\tilde b}}e^{-K\varrho}\cdot e^{\frac{16 E^2\delta_0(\gamma_0-\gamma)Kr}{\gamma^2}}(\|P^{20}_{nm}\|_{\!{}_{D(r),\Cal O}}+\|P^{02}_{nm}\|_{\!{}_{D(r),\Cal O}})\\
&\;&+\lim_{n\rightarrow \infty}\frac{1}{|n|}e^{-K\varrho}\cdot \frac{2(1+e)^{2\tilde b}}{\varrho^{2\tilde b}}(\|P^{20}_{nn}\|_{\!{}_{D(r),\Cal O}}+\|P^{02}_{nn}\|_{\!{}_{D(r),\Cal O}})\\
&\leq&  e^{-K\varrho}\cdot \frac{2(1+e)^{2\tilde b}}{\varrho^{2\tilde b}}\|X_{P}\|_{\!{}_{s,\bar a, \rho,D(r,s),\Cal O}}
\overset{(\ref{XRt})}{\leq}\eta\varepsilon\leq\varepsilon_+,
\end{eqnarray*}
\begin{eqnarray*}
&&\|\frac{1}{|n|}\sum_{\upsilon=\pm}\frac{\partial^2 \hat R_{2}}{\partial w^\upsilon_n\partial w^\upsilon_n}-
\lim_{n\rightarrow\infty}\frac{1}{|n|}\sum_{\upsilon=\pm}\frac{\partial^2 \hat R_{2}}{\partial w^\upsilon_n\partial w^\upsilon_n}\|_{\!{}_{D_+,\Cal O_+}}\\
&\leq &\|\frac{1}{|n|}\sum_{\upsilon=\pm}\frac{\partial^2 \hat R_{2,<}}{\partial w^\upsilon_n\partial w^\upsilon_n}-\lim_{n\rightarrow\infty}\frac{1}{|n|}\sum_{\upsilon=\pm}\frac{\partial^2 \hat R_{2,<}}{\partial w^\upsilon_n\partial w^\upsilon_n}\|_{\!{}_{D_+,\Cal O_+}}\\
&+&\|\frac{1}{|n|}\sum_{\upsilon=\pm}\frac{\partial^2 \hat R_{2,>}}{\partial w^\upsilon_n\partial w^\upsilon_n}- \lim_{n\rightarrow\infty}\frac{1}{|n|}\sum_{\upsilon=\pm}\frac{\partial^2 \hat R_{2,>}}{\partial w^\upsilon_n\partial w^\upsilon_n}\|_{\!{}_{D_+,\Cal O_+}}\\
&\leq& \|\frac{1}{|n|}\sum_{\upsilon=\pm}\frac{\partial^2 \hat R_{2,<}}{\partial w^\upsilon_n\partial w^\upsilon_n}\|_{\!{}_{D_+,\Cal O_+}}+0\\
&\leq& \frac{1}{|n|}\frac{cE^2K^{2\tau+2}\delta_0}{\gamma_0^2\varrho^{2\tilde b}}e^{-K\varrho}\cdot e^{\frac{16 E^2\delta_0(\gamma_0-\gamma)Kr}{\gamma^2}}\|X_{P}\|_{\!{}_{s,\bar a, \rho,D(r,s),\Cal O}}
\overset{(\ref{XRt})}{\leq}\frac{\eta\varepsilon}{|n|}\leq \frac{\varepsilon_+}{|n|}.
\end{eqnarray*}
For the term $P_{(\geq2)}=P_{(2)}+P_{(3)}$, we observe that $P_{(2)}$ in $(\ref{P(2)})$ with the indices $l,\alpha,\beta$ satisfying $2|l|+|\alpha+\beta|\leq 2$ and $\alpha,\beta\in \{ \alpha_n+\beta_n\in 2\N,\forall |n|> EK\}$, the second derivatives in $w,\bar w$ of the terms with $|l|=1$ in $P_{(2)}$ disappear so it can be specifically written as
\beq
\sum_{\upsilon=\pm}\frac{\partial^2  P_{(2)}}{\partial w^\upsilon_n\partial w^\upsilon_n}=\sum_{|k|>K}(2P_{k02e_n0}
+P_{k0e_ne_n}
+2P_{k002e_n})\kth=\left \{\begin{array}{lr}
P_{(2,n<)},&|n|\leq EK,\\
P_{(2,n>)},&|n|>EK,\nonumber\\
\end{array}\right.\eeq
when $|n|\leq EK$, by $(\ref{decay11})$,  the coefficients $P_{k02e_n0},P_{k0e_ne_n}, P_{k002e_n}$ in the norm $\|\cdot\|_{D(r),\Cal O}$ are all bounded, when $|n|>EK$,
$$ \|P_{k02e_n0}\|_{D(r),\Cal O},\|P_{k02e_n0}\|_{D(r),\Cal O},\|P_{k02e_n0}\|_{D(r),\Cal O}\leq c |n|\|X_P\|_{s,\bar a, \rho,D(r,s),\Cal O}.$$
Similarly considering $P_{(3)}$ in $(\ref{P(3)})$ with $2|l|+|\alpha+\beta|> 2$ and $\alpha,\beta\in \{ \alpha_n+\beta_n\in 2\N,\forall |n|> EK\}$, we can rewrite $P_{(3)}=P_{(3),1}+P_{(3),2}$,
\begin{eqnarray*}
P_{(3),1}=\sum_{l,\alpha,\beta}
P_{l\alpha\beta}I^lw^{\alpha}\bar w^{\beta},\ \  P_{(3),2}=\sum_{l,\alpha'\beta'}
P_{l\alpha\beta}I^lw^{\alpha'}
\bar w^{\beta'},
\end{eqnarray*}
where $\alpha,\beta\in \{\alpha_n+\beta_n=0,\forall |n|> EK\}$ in $P_{(3),1}$, $\alpha',\beta'\in \{\sum\limits_{|n|>EK}\alpha'_n+\beta'_n >0, \alpha'_n+\beta'_n\in 2\N,\forall |n|> EK\}$ in $P_{(3),2}$. Besides, due to the decay property $(\ref{decay1})$ and Lemma 7.2, $P_{l\alpha\beta}$ in $P_{(3),1}$ in the norm $\|\cdot\|_{D(r),\Cal O}$ are all bounded for any $l,\alpha,\beta\in \{\alpha_n+\beta_n=0,\forall |n|> EK\}$. Then we calculate
\beq \frac{1}{|n|}\sum_{\upsilon=\pm}\frac{\partial^2 P_{(3)}}{\partial w^\upsilon_n\partial w^\upsilon_n}=\frac{1}{|n|}\left \{\begin{array}{lr}
P_{(3,n<)},&|n|\leq EK,\\
P_{(3,n>)},&|n|>EK,\nonumber\\
\end{array}\right.\eeq
where $P_{(3,n<)}=P_{(3,n<),1}+P_{(3,n<),2}$,
\begin{eqnarray*}
P_{(3,n<),1}&=&\sum_{\upsilon=\pm}\frac{\partial^2 P_{(3),1}}{\partial w^\upsilon_n\partial w^\upsilon_n}, \
P_{(3,n<),2}=\sum_{\upsilon=\pm}\frac{\partial^2 P_{(3),2}}{\partial w^\upsilon_n\partial w^\upsilon_n}, \ \ |n|\leq EK,\\
P_{(3,n>)}&=&\sum_{\upsilon=\pm}\frac{\partial^2 P_{(3),2}}{\partial w^\upsilon_n\partial w^\upsilon_n}, \ |n|>EK.
\end{eqnarray*}
In this way, associated with the estimates $(\ref{XP(2)}),(\ref{XP(3)})$,  we have
\begin{eqnarray*}
&&\|\lim_{n\rightarrow \infty}\frac{1}{|n|}\sum_{\upsilon=\pm}\frac{\partial^2P_{(\geq2)}}{\partial w^\upsilon_n\partial w^\upsilon_n}\|_{\!{}_{D_+,\cal{O_+}}}\\
&\leq& \|\lim_{n\rightarrow \infty}\frac{1}{|n|}\sum_{\upsilon=\pm}\frac{\partial^2P_{(3)}}{\partial w^\upsilon_n\partial w^\upsilon_n}\|_{\!{}_{D_+,\cal{O_+}}}
+\|\lim_{n\rightarrow \infty}\frac{1}{|n|}\sum_{\upsilon=\pm}\frac{\partial^2P_{(2)}}{\partial w^\upsilon_n\partial w^\upsilon_n}\|_{\!{}_{D_+,\cal{O_+}}}\\
&\leq& \|\lim_{n\rightarrow \infty}\frac{1}{|n|}P_{(2,n>)}\|_{\!{}_{D_+,\cal{O_+}}}+\|\lim_{n\rightarrow \infty}\frac{1}{|n|}P_{(3,n<),2}\|_{\!{}_{D_+,\cal{O_+}}}
+\|\lim_{n\rightarrow \infty}\frac{1}{|n|}P_{(3,n>)}\|_{\!{}_{D_+,\cal{O_+}}}\\
&\leq&\|X_{P_{(2)}}\|_{\!{}_{s_+,\bar a, \rho_+, D_+,\Cal O_+}}+\|X_{P_{(3)}}\|_{\!{}_{s_+,\bar a, \rho_+, D_+,\Cal O_+}}\\
&\leq&(c\eta^{-1}e^{-K\varrho}+c\eta) \|X_{P}\|_{\!{}_{s,\bar a, \rho,D(r,s),\Cal O}}
\overset{(\ref{XP-R1}),(\ref{XP-R2})}{\leq}\eta\varepsilon\leq\varepsilon_+,\\
&&\|\frac{1}{|n|}\sum_{\upsilon=\pm}\frac{\partial^2P_{(\geq2)}}{\partial w^\upsilon_n\partial w^\upsilon_n}-  \lim_{n\rightarrow\infty}\frac{1}{|n|}\sum_{\upsilon=\pm}\frac{\partial^2P_{(\geq2)}}{\partial w^\upsilon_n\partial w^\upsilon_n}\|_{\!{}_{D_+,\Cal O_+}}\\
&\leq & \|\frac{1}{|n|}\sum_{\upsilon=\pm}\frac{\partial^2P_{(2)}}{\partial w^\upsilon_n\partial w^\upsilon_n}- \lim_{n\rightarrow\infty}\frac{1}{|n|}\sum_{\upsilon=\pm}\frac{\partial^2P_{(2)}}{\partial w^\upsilon_n\partial w^\upsilon_n}\|_{\!{}_{D_+,\Cal O_+}}\\
&+& \|\frac{1}{|n|}\sum_{\upsilon=\pm}\frac{\partial^2P_{(3)}}{\partial w^\upsilon_n\partial w^\upsilon_n}- \lim_{n\rightarrow\infty}\frac{1}{|n|}\sum_{\upsilon=\pm}\frac{\partial^2P_{(3)}}{\partial w^\upsilon_n\partial w^\upsilon_n}\|_{\!{}_{D_+,\Cal O_+}}\\
&\leq & \|\frac{1}{|n|}\sum_{\upsilon=\pm}\frac{\partial^2P_{(2,n<)}}{\partial w^\upsilon_n\partial w^\upsilon_n}\|_{\!{}_{D_+,\Cal O_+}}
+ \|\frac{1}{|n|}\sum_{\upsilon=\pm}\frac{\partial^2P_{(3,n<),1}}{\partial w^\upsilon_n\partial w^\upsilon_n}\|_{\!{}_{D_+,\Cal O_+}}\\
&\leq& \frac{1}{|n|}(c\eta^{-1}e^{-K\varrho}+c\eta) \|X_{P}\|_{\!{}_{s,\bar a, \rho,D(r,s),\Cal O}}\overset{(\ref{XP-R1}),(\ref{XP-R2})}{\leq}\frac{\eta\varepsilon}{|n|}\leq\frac{\varepsilon_+}{|n|}.
\end{eqnarray*}
In the following we consider $\{P,F\}=\{P^1,F\}+\{P^2,F\}+\{P^3,F\}$,
\begin{eqnarray}
&&\sum_{\upsilon=\pm}\frac{\partial^2 \{P^1,F\}}{\partial w^\upsilon_n\partial w^\upsilon_n}\nonumber\\
&=&\frac{\partial^3 P^1}{\partial w_n\partial w_n\partial I}\cdot\frac{\partial F}{\partial\tilde{\theta}}+\frac{\partial^2 P^1}{\partial w_n\partial I}\cdot\frac{\partial^2 F}{\partial w_n\partial\tilde{\theta}}+\frac{\partial P^1}{\partial I}\cdot\frac{\partial^3 F}{\partial w_n\partial w_n\partial\tilde{\theta}}-\frac{\partial F_0}{\partial I}\cdot\frac{\partial^3 P^1}{\partial w_n\partial w_n\partial\tilde{\theta}}\nonumber\\
&+&{\rm i}\sum_{m\in \Z_1}(\frac{\partial^3 P^1}{\partial w_n\partial w_n\partial w_m}\frac{\partial F}{\partial\bar{w}_m}+\frac{\partial^2 P^1}{\partial w_n\partial w_m}\frac{\partial^2 F}{\partial w_n\partial\bar{w}_m}
-\frac{\partial^3 P^1}{\partial w_n\partial w_n\partial \bar w_m}\frac{\partial F}{\partial w_m}\nonumber\\
&-&\frac{\partial^2 P^1}{\partial w_n\partial \bar w_m}\frac{\partial^2 F}{\partial w_n\partial w_m})
+\frac{\partial^3 P^1}{\partial w_n\partial \bar w_n\partial I}\cdot\frac{\partial F}{\partial\tilde{\theta}}
+\frac{\partial^2 P^1}{\partial w_n\partial I}\cdot\frac{\partial^2 F}{\partial\bar w_n\partial\tilde{\theta}}+\frac{\partial^2 P^1}{\partial \bar w_n\partial I}\cdot\frac{\partial^2 F}{\partial w_n\partial\tilde{\theta}}\nonumber\\
&+&\frac{\partial P^1}{\partial I}\cdot\frac{\partial^3 F}{\partial w_n\partial\bar w_n\partial\tilde{\theta}}
-\frac{\partial F_0}{\partial I}\cdot\frac{\partial^3 P^1}{\partial w_n\partial \bar w_n\partial\tilde{\theta}}
+{\rm i}\sum_{m\in \Z_1}(\frac{\partial^3 P^1}{\partial w_n\partial \bar w_n\partial w_m}\frac{\partial F}{\partial\bar{w}_m}\nonumber\\
&+&\frac{\partial^2 P^1}{\partial w_n\partial w_m}\frac{\partial^2 F}{\partial \bar w_n\partial\bar{w}_m}+\frac{\partial^2 P^1}{\partial \bar w_n\partial w_m}\frac{\partial^2 F}{\partial w_n\partial\bar{w}_m}
-\frac{\partial^3 P^1}{\partial w_n\partial \bar w_n\partial \bar w_m}\frac{\partial F}{\partial w_m}\nonumber\\
&-&\frac{\partial^2 P^1}{\partial w_n\partial \bar w_m}\frac{\partial^2 F}{\partial \bar w_n\partial w_m}-\frac{\partial^2 P^1}{\partial \bar w_n\partial \bar w_m}\frac{\partial^2 F}{\partial w_n\partial w_m})
+\frac{\partial^3 P^1}{\partial \bar w_n\partial \bar w_n\partial I}\cdot\frac{\partial F}{\partial\tilde{\theta}}
+\frac{\partial^2 P^1}{\partial \bar w_n\partial I}\cdot\frac{\partial^2 F}{\partial \bar w_n\partial\tilde{\theta}}\nonumber\\
&+&\frac{\partial P^1}{\partial I}\cdot\frac{\partial^3 F}{\partial \bar w_n\partial \bar w_n\partial\tilde{\theta}}-\frac{\partial F_0}{\partial I}\cdot\frac{\partial^3 P^1}{\partial \bar w_n\partial \bar w_n\partial\tilde{\theta}}
+{\rm i}\sum_{m\in \Z_1}(\frac{\partial^3 P^1}{\partial \bar w_n\partial \bar w_n\partial w_m}\frac{\partial F}{\partial\bar{w}_m}\nonumber\\
&+&\frac{\partial^2 P^1}{\partial \bar w_n\partial w_m}\frac{\partial^2 F}{\partial \bar w_n\partial\bar{w}_m}-\frac{\partial^3 P^1}{\partial \bar w_n\partial \bar w_n\partial \bar w_m}\frac{\partial F}{\partial w_m}
-\frac{\partial^2 P^1}{\partial \bar w_n\partial \bar w_m}\frac{\partial^2 F}{\partial \bar w_n\partial w_m}).\label{P1F}
\end{eqnarray}
For the term $P^2$ with $\alpha,\beta\in\{|\alpha+\beta|=\alpha_n+\beta_n=1,\forall |n|>E_-K_-\}$, it is obvious that $\sum\limits_{\upsilon=\pm}\frac{\partial^2 \{P^2,F\}}{\partial w^\upsilon_n\partial w^\upsilon_n}$ vanishes for any $n\in Z_1$. Similarly, with $\alpha,\beta\in \{\alpha_n+\beta_n=0,\forall |n|>E_-K_-\}$ in $P^3$, we can get the same formula as $(\ref{P1F})$ with $ P^3$ in place of $P^1$ and the sum index $m$ is limited to less than $E_-K_-.$
\begin{Lemma}\label{Lem4.5} Let $D_+=D(r_+,s_+)$ with $r_+=\frac{r}{2}, s_+=\eta s$ defined in section 5, for any $n\in \Z_1$ and a constant $0<\beta'\leq \frac{1}{4}$, we get some estimates in the following
\begin{eqnarray*}
&&\|\frac{\partial P^1}{\partial I}\cdot\frac{\partial^3 F}{\partial w_n\partial w_n\partial\tilde{\theta}}\|_{D_+,\Cal O_+},\|\frac{\partial P^1}{\partial I}\cdot\frac{\partial^3 F}{\partial w_n\partial \bar w_n\partial\tilde{\theta}}\|_{D_+,\Cal O_+},\|\frac{\partial P^1}{\partial I}\cdot\frac{\partial^3 F}{\partial \bar w_n\partial \bar w_n\partial\tilde{\theta}}\|_{D_+,\Cal O_+},\\
&&\|\frac{\partial^2 P^1}{\partial w_n\partial I}\cdot\frac{\partial^2 F}{\partial w_n\partial\tilde{\theta}}\|_{D_+,\Cal O_+},\ \|\frac{\partial^2 P^1}{\partial w_n\partial I}\cdot\frac{\partial^2 F}{\partial \bar w_n\partial\tilde{\theta}}\|_{D_+,\Cal O_+},\ \|\frac{\partial^2 P^1}{\partial \bar w_n\partial I}\cdot\frac{\partial^2 F}{\partial w_n\partial\tilde{\theta}}\|_{D_+,\Cal O_+},\\
 &&\|\frac{\partial^2 P^1}{\partial \bar w_n\partial I}\cdot\frac{\partial^2 F}{\partial \bar w_n\partial\tilde{\theta}}\|_{D_+,\Cal O_+}
\leq\Cal B_\varrho^{\frac{1}{2}}\varepsilon^{2-\beta'},\\
&&\|\frac{\partial^3 P^1}{\partial w_n\partial w_n\partial I}\cdot\frac{\partial F}{\partial\tilde{\theta}}\|_{D_+,\Cal O_+},\ \|\frac{\partial^3 P^1}{\partial w_n\partial \bar w_n\partial I}\cdot\frac{\partial F}{\partial\tilde{\theta}}\|_{D_+,\Cal O_+},\ \|\frac{\partial^3 P^1}{\partial \bar w_n\partial \bar w_n\partial I}\cdot\frac{\partial F}{\partial\tilde{\theta}}\|_{D_+,\Cal O_+},\\
&&\|\frac{\partial^3 P^1}{\partial w_n\partial w_n\partial\tilde{\theta} }\cdot\frac{\partial F}{\partial I}\|_{D_+,\Cal O_+},\ \|\frac{\partial^3 P^1}{\partial w_n\partial \bar w_n\partial\tilde{\theta} }\cdot\frac{\partial F}{\partial I}\|_{D_+,\Cal O_+}, \ \|\frac{\partial^3 P^1}{\partial \bar w_n\partial \bar w_n\partial\tilde{\theta} }\cdot\frac{\partial F}{\partial I}\|_{D_+,\Cal O_+},\\
&&\leq |n|\Cal B_\varrho^{\frac{1}{2}}\varepsilon^{2-\beta'},
\end{eqnarray*}
\begin{eqnarray*}
&&\|\frac{\partial^2 P^1}{\partial w_n\partial w_m}\frac{\partial^2 F}{\partial w_n\partial\bar{w}_m}\|_{D_+,\Cal O_+},\|\frac{\partial^2 P^1}{\partial w_n\partial \bar w_m}\frac{\partial^2 F}{\partial \bar w_n\partial w_m}\|_{D_+,\Cal O_+},\|\frac{\partial^2 P^1}{\partial \bar w_n\partial w_m}\frac{\partial^2 F}{\partial w_n\partial\bar{w}_m}\|_{D_+,\Cal O_+},\\
&&\|\frac{\partial^2 P^1}{\partial \bar w_n\partial \bar w_m}\frac{\partial^2 F}{\partial \bar w_n\partial w_m}\|_{D_+,\Cal O_+},\|\frac{\partial^2 P^1}{\partial w_n\partial w_m}\frac{\partial^2 F}{\partial \bar w_n\partial\bar{w}_m}\|_{D_+,\Cal O_+},\|\frac{\partial^2 P^1}{\partial \bar w_n\partial w_m}\frac{\partial^2 F}{\partial \bar w_n\partial\bar{w}_m}\|_{D_+,\Cal O_+},\\
&&\|\frac{\partial^2 P^1}{\partial w_n\partial \bar w_m}\frac{\partial^2 F}{\partial w_n\partial w_m}\|_{D_+,\Cal O_+},\|\frac{\partial^2 P^1}{\partial \bar w_n\partial \bar w_m}\frac{\partial^2 F}{\partial w_n\partial w_m}\|_{D_+,\Cal O_+},\\
&&\leq |m|^{a}|n|^{-\bar a}e^{-2|n|\rho}\Cal B_\varrho^{\frac{1}{2}}\varepsilon^{2-\beta'},|n\neq m|\leq E_-K_-;\ or
\ \leq |n|\Cal B_\varrho^{\frac{1}{2}}\varepsilon^{2-\beta'},\ n=m\in \Z_1,\\
&&\|\frac{\partial^3 P^1}{\partial w_n\partial w_n\partial w_m}\frac{\partial F}{\partial\bar{w}_m}\|_{D_+,\Cal O_+},\|\frac{\partial^3 P^1}{\partial w_n\partial \bar w_n\partial w_m}\frac{\partial F}{\partial\bar{w}_m}\|_{D_+,\Cal O_+},\|\frac{\partial^3 P^1}{\partial \bar w_n\partial \bar w_n\partial w_m}\frac{\partial F}{\partial\bar{w}_m}\|_{D_+,\Cal O_+},\\
&&\|\frac{\partial^3 P^1}{\partial w_n\partial w_n\partial \bar w_m}\frac{\partial F}{\partial w_m}\|_{D_+,\Cal O_+},\|\frac{\partial^3 P^1}{\partial w_n\partial \bar w_n\partial \bar w_m}\frac{\partial F}{\partial w_m}\|_{D_+,\Cal O_+},\|\frac{\partial^3 P^1}{\partial \bar w_n\partial \bar w_n\partial \bar w_m}\frac{\partial F}{\partial w_m}\|_{D_+,\Cal O_+},\\
&&\|\frac{\partial^3 P^1}{\partial w_n\partial w_n\partial w_m}\frac{\partial F}{\partial\bar{w}_m}\|_{D_+,\Cal O_+},\|\frac{\partial^3 P^1}{\partial w_n\partial \bar w_n\partial w_m}\frac{\partial F}{\partial\bar{w}_m}\|_{D_+,\Cal O_+},\|\frac{\partial^3 P^1}{\partial \bar w_n\partial \bar w_n\partial w_m}\frac{\partial F}{\partial\bar{w}_m}\|_{D_+,\Cal O_+},\\
&&\|\frac{\partial^3 P^1}{\partial w_n\partial w_n\partial \bar w_m}\frac{\partial F}{\partial w_m}\|_{D_+,\Cal O_+},\|\frac{\partial^3 P^1}{\partial w_n\partial \bar w_n\partial \bar w_m}\frac{\partial F}{\partial w_m}\|_{D_+,\Cal O_+}, \|\frac{\partial^3 P^1}{\partial \bar w_n\partial \bar w_n\partial \bar w_m}\frac{\partial F}{\partial w_m}\|_{D_+,\Cal O_+},\\
&&\leq |n||m|^{-a-\bar a}e^{-2|m|\rho}\Cal B_\varrho^{\frac{1}{2}}\varepsilon^{2-\beta'},\ n,m\in\Z_1.
\end{eqnarray*}
\end{Lemma}

\proof
In the above inequality estimates, we mainly consider the following six kinds of terms respectively and the others can be obtained by the similar arguments.\\
$(1)$ For the term $\frac{\partial P^1}{\partial I}\cdot\frac{\partial^3 F}{\partial w_n\partial w_n\partial\tilde{\theta}}=\frac{\partial P^1}{\partial I}\cdot\frac{\partial F^{20}_{nn}}{\partial\tilde{\theta}},$
it is obvious that $\frac{\partial P^2}{\partial I}$ is at least of order 2 in $w,\bar w$, associated with Lemma $\ref{Lem7.2}$, we have
$$\|\frac{\partial P^1}{\partial I}\|_{D_+,\Cal O_+}\leq \|X_P\|_{s,\bar a,\rho,D(r,s),\cal O},\ \|\frac{\partial F^{20}_{nn}}{\partial\tilde{\theta}}\|_{D_+,\Cal O_+}\leq \frac{c}{\varrho}\|F^{20}_{nn}\|_{D_+,\Cal O_+},$$
then
\begin{eqnarray*}
\|\frac{\partial P^1}{\partial I}\cdot\frac{\partial^3 F}{\partial w_n\partial w_n\partial\tilde{\theta}}\|_{D_+,\Cal O_+}&\leq& \frac{c}{\varrho}\|X_P\|_{s,\bar a,\rho,D(r,s),\cal O}\cdot\|F^{20}_{nn}\|_{D_+,\Cal O_+}\\
&\leq& \|X_P\|_{s,\bar a,\rho,D(r,s),\cal O}\cdot\|X_F\|_{s,a,\rho,D(r-3\varrho,s),\Cal O_+}\\
&\overset{(\ref{XP}),(\ref{PHI})}{\leq}& \Cal B_\varrho^{\frac{1}{2}}\varepsilon^{2-\beta'}.
\end{eqnarray*}
$(2)$ For the term $\frac{\partial^2 P^1}{\partial w_n\partial I}\cdot\frac{\partial^2 F}{\partial w_n\partial\tilde{\theta}}$, associated with Lemma $\ref{Lem7.2}$ and the definition of the vector field norm, we have
$$ \|\frac{\partial^2 P^1}{\partial w_n\partial I}\|_{D_+,\Cal O_+}\leq \frac{c}{s}|n|^{a}e^{|n|\rho}\|\frac{\partial P^1}{\partial I}\|_{D_+,\Cal O_+}\leq   \frac{c}{s}|n|^{a}e^{|n|\rho}\|X_P\|_{s,\bar a,\rho,D(r,s),\cal O},$$
$$\|\frac{\partial^2 F}{\partial w_n\partial\tilde{\theta}}\|_{D_+,\Cal O_+}\leq \frac{c}{\varrho}\|\frac{\partial F}{\partial w_n}\|_{D_+,\Cal O_+}\leq \frac{c|n|^{-a}e^{-|n|\rho}}{\varrho}\|X_F\|_{s,a,\rho,D(r-3\varrho,s),\Cal O_+},$$
then
\begin{eqnarray*}
\|\frac{\partial^2 P^1}{\partial w_n\partial I}\cdot\frac{\partial^2 F}{\partial w_n\partial\tilde{\theta}}\|_{D_+,\Cal O_+}
&\leq& \frac{c}{s\varrho}\|X_F\|_{s,a,\rho,D(r-3\varrho,s),\Cal O_+}\cdot\|X_P\|_{s,\bar a,\rho,D(r,s),\cal O}\\
&\leq& \Cal B_\varrho^{\frac{1}{2}}\varepsilon^{2-\beta'}.
\end{eqnarray*}
$(3)$ For the term $\frac{\partial^3 P^1}{\partial w_n\partial w_n\partial I}\cdot\frac{\partial F}{\partial\tilde{\theta}}$, by the assumption $(A6)$ of $P$,  we have for any $n\in \Z_1$
$$\|\frac{\partial^3 P^1}{\partial w_n\partial w_n\partial I}\|_{D_+,\Cal O_+}\leq c|n|\cdot\|X_P\|_{s,\bar a,\rho,D(r,s),\cal O},$$
and by the definition of the vector field norm,
\begin{eqnarray*}\|\frac{\partial F}{\partial\tilde{\theta}}\|_{D_+,\Cal O_+}&\leq& s^2\|X_F\|_{s,a,\rho,D(r-3\varrho,s),\Cal O_+},\\
\|\frac{\partial^3 P^2}{\partial w_n\partial w_n\partial I}\cdot\frac{\partial F}{\partial\tilde{\theta}}\|_{D_+,\Cal O_+}&\leq& cs^2|n|\|X_P\|_{s,\bar a,\rho,D(r,s),\cal O}\cdot\|X_F\|_{s,a,\rho,D(r-3\varrho,s),\Cal O_+}\\
&\leq& |n|\Cal B_\varrho^{\frac{1}{2}}\varepsilon^{2-\beta'}.
\end{eqnarray*}
$(4)$ For the term $\frac{\partial^2 P^1}{\partial w_n\partial w_m}\frac{\partial^2 F}{\partial w_n\partial\bar{w}_m}$, associated with Lemma $\ref{Lem7.2}$ and the definition of the vector field norm, we have if $|n\neq m|\leq E_-K_-$,
$$\|\frac{\partial^2 P^1}{\partial w_n\partial w_m}\|_{D_+,\Cal O_+}\leq \frac{c}{s}|m|^{a}e^{|m|\rho}\|\frac{\partial P^1}{\partial w_n}\|_{D_+,\Cal O_+}\leq c|m|^{a}e^{|m|\rho}|n|^{-\bar a}e^{-|n|\rho}\|X_P\|_{s,\bar a,\rho,D(r,s),\cal O},$$
$$\|\frac{\partial^2 F}{\partial w_n\partial\bar{w}_m}\|_{D_+,\Cal O_+}\leq \|F^{11}_{nm}\|_{D_+,\Cal O_+}\overset{(\ref{decay11}),(\ref{F11}),(\ref{PHI})}{\leq}e^{-|n|\rho}e^{-|m|\rho}\Cal B^{\frac{1}{2}}_{\varrho}\varepsilon^{1-\beta'},$$
if $|m=n|\leq E_-K_-$ or $|m|\leq E_-K_-,|n|>E_-K_-$ or $|m|>E_-K_-,|n|\leq E_-K_-$ or $|m|,|n|>E_-K_-$, $\frac{\partial^2 F}{\partial w_n\partial\bar{w}_m}$ vanishes, namely $\|\frac{\partial^2 F}{\partial w_n\partial\bar{w}_m}\|_{D_+,\Cal O_+}=0$,
hence
\begin{eqnarray*}
\|\frac{\partial^2 P^1}{\partial w_n\partial w_m}\frac{\partial^2 F}{\partial w_n\partial\bar{w}_m}\|_{D_+,\Cal O_+}
&\leq & \Cal B^{\frac{1}{2}}_{\varrho}\varepsilon^{1-\beta'}|m|^{a}|n|^{-\bar a}e^{-2|n|\rho}\|X_P\|_{s,\bar a,\rho,D(r,s),\cal O}\\
&\overset{(\ref{XP})}{\leq}& |m|^{a}|n|^{-\bar a}e^{-2|n|\rho}\Cal B_\varrho^{\frac{1}{2}}\varepsilon^{2-\beta'},\ \  |n\neq m|\leq E_-K_-.
\end{eqnarray*}
$(5)$ For the term $\frac{\partial^2 P^1}{\partial w_n\partial w_m}\frac{\partial^2 F}{\partial \bar w_n\partial\bar{w}_m}$, if $|n\neq m|\leq E_-K_-$,  associated with the estimates in $(\ref{decay11})$, then
$$\|\frac{\partial^2 P^1}{\partial w_n\partial w_m}\|_{D_+,\Cal O_+}\leq \frac{c}{s}|m|^{a}e^{|m|\rho}\|\frac{\partial P^1}{\partial w_n}\|_{D_+,\Cal O_+}\leq c|m|^{a}e^{|m|\rho}|n|^{-\bar a}e^{-|n|\rho}\|X_P\|_{s,\bar a,\rho,D(r,s),\cal O},$$
 $$\|\frac{\partial^2 F}{\partial \bar w_n\partial\bar{w}_m}\|_{D_+,\Cal O_+}\leq \|F^{02}_{nm}\|_{D_+,\Cal O_+}\overset{(\ref{decay11}),(\ref{F11}),(\ref{PHI})}{\leq}e^{-|n|\rho}e^{-|m|\rho}\Cal B^{\frac{1}{2}}_{\varrho}\varepsilon^{1-\beta'},$$
if $|m=n|\leq E_-K_-,$
$$\|\frac{\partial^2 P^1}{\partial w_n\partial w_n}\|_{D_+,\Cal O_+}\leq |n|\|X_P\|_{s,\bar a,\rho,D(r,s),\cal O},\ \|\frac{\partial^2 F}{\partial \bar w_n\partial\bar{w}_n}\|_{D_+,\Cal O_+}\leq \|F^{02}_{nn}\|_{D_+,\Cal O_+},$$
if $|m|\leq E_-K_-,|n|>E_-K_-$ or $|m|>E_-K_-,|n|\leq E_-K_-$, $\frac{\partial^2 F}{\partial \bar w_n\partial\bar{w}_m}$ vanishes, namely
$$\|\frac{\partial^2 F}{\partial \bar w_n\partial\bar{w}_m}\|_{D_+,\Cal O_+}=0$$
if $|n|,|m|> E_-K_-$, $\frac{\partial^2 F}{\partial \bar w_n\partial\bar{w}_m}$ exists if and only if $n=m$, then we have
$$\|\frac{\partial^2 P^1}{\partial w_n\partial w_n}\|_{D_+,\Cal O_+}\leq |n|\|X_P\|_{s,\bar a,\rho,D(r,s),\cal O},\ \|\frac{\partial^2 F}{\partial \bar w_n\partial\bar{w}_n}\|_{D_+,\Cal O_+}\leq \|F^{02}_{nn}\|_{D_+,\Cal O_+},$$
and get the estimates
\begin{eqnarray*}
\|\frac{\partial^2 P^1}{\partial w_n\partial w_m}\frac{\partial^2 F}{\partial \bar w_n\partial\bar{w}_m}\|_{D_+,\Cal O_+}
&\leq& |m|^{a}|n|^{-\bar a}e^{-2|n|\rho}\Cal B_\varrho^{\frac{1}{2}}\varepsilon^{2-\beta'},\ |n\neq m|\leq E_-K_-,\\
\|\frac{\partial^2 P^1}{\partial w_n\partial w_m}\frac{\partial F}{\partial \bar w_n\partial\bar{w}_m}\|_{D_+,\Cal O_+}
&\leq& |n|\|X_P\|_{s,\bar a,\rho,D(r,s),\cal O}\cdot\|F^{02}_{nn}\|_{D_+,\Cal O_+}\\
&\leq& |n|\Cal B_\varrho^{\frac{1}{2}}\varepsilon^{2-\beta'},\ n=m.
\end{eqnarray*}
$(6)$ For the term $ \frac{\partial^3 P^1}{\partial w_n\partial w_n\partial w_m}\frac{\partial F}{\partial\bar{w}_m}$, using the assumption $(A6)$ of $P$ and the definition of the vector field norm, we have
$$\|\frac{\partial^3 P^1}{\partial w_n\partial w_n\partial w_m}\|_{D_+,\Cal O_+}\leq c|n||m|^{-\bar a}e^{-|m|\rho}\|X_P\|_{s,\bar a,\rho,D(r,s),\cal O},$$
$$\|\frac{\partial F}{\partial\bar{w}_m}\|_{D_+,\Cal O_+}\leq s|m|^{-a}e^{-|m|\rho}\|X_F\|_{s,a,\rho,D(r-3\varrho,s),\Cal O_+},$$
hence for any $n,m \in \Z_1$, we have
\begin{eqnarray*}
&&\|\frac{\partial^3 P^1}{\partial w_n\partial w_n\partial w_m}\frac{\partial F}{\partial\bar{w}_m}\|_{D_+,\Cal O_+}\leq |n||m|^{-a-\bar a}e^{-2|m|\rho}\Cal B_\varrho^{\frac{1}{2}}\varepsilon^{2-\beta'}.
\end{eqnarray*}\qed

In the above lemma, if the term $P^1$ is replaced with $P^3$, we can get the same results or even better. So it is sufficient for us to calculus $\frac{1}{|n|}\sum\limits_{\upsilon=\pm}\frac{\partial^2\{P^1,F\}}{\partial w^\upsilon_n\partial w^\upsilon_n}$ and $$\frac{1}{|n|}\sum\limits_{\upsilon=\pm}\frac{\partial^2\{P^1,F\}}{\partial w^\upsilon_n\partial w^\upsilon_n}-\lim\limits_{n\rightarrow \infty}\frac{1}{|n|}\sum\limits_{\upsilon=\pm}\frac{\partial^2\{P^1,F\}}{\partial w^\upsilon_n\partial w^\upsilon_n}.$$
The estimates of the term $\frac{1}{|n|}\sum\limits_{\upsilon=\pm}\frac{\partial^2\{P^3,F\}}{\partial w^\upsilon_n\partial w^\upsilon_n}$ and
$$\frac{1}{|n|}\sum\limits_{\upsilon=\pm}\frac{\partial^2\{P^3,F\}}{\partial w^\upsilon_n\partial w^\upsilon_n}-\lim\limits_{n\rightarrow \infty}\frac{1}{|n|}\sum\limits_{\upsilon=\pm}\frac{\partial^2\{P^3,F\}}{\partial w^\upsilon_n\partial w^\upsilon_n}$$
 can be obtained with the same arguments.\\
By Lemma $\ref{Lem4.5}$, we obtain the estimate of $ (\ref{P1F})$ with the careful calculations
\begin{eqnarray*}
&&\frac{1}{|n|}\|\sum_{\upsilon=\pm}\frac{\partial^2\{P^1,F\}}{\partial w^\upsilon_n\partial w^\upsilon_n}\|_{\!{}_{D_+,\Cal O_+}}\\
&\leq&\frac{1}{|n|}\Big(\|\frac{\partial^3 P^1}{\partial w_n\partial w_n\partial I}\cdot\frac{\partial F}{\partial\tilde{\theta}}\|_{\!{}_{D_+,\Cal O_+}}+\|\frac{\partial^2 P^1}{\partial w_n\partial I}\cdot\frac{\partial^2 F}{\partial w_n\partial\tilde{\theta}}\|_{\!{}_{D_+,\Cal O_+}}\\
&+&\|\frac{\partial P^1}{\partial I}\cdot\frac{\partial^3 F}{\partial w_n\partial w_n\partial\tilde{\theta}}\|_{\!{}_{D_+,\Cal O_+}}
+\|\frac{\partial F_0}{\partial I}\cdot\frac{\partial^3 P^1}{\partial w_n\partial w_n\partial\tilde{\theta}}\|_{\!{}_{D_+,\Cal O_+}}\\
&+&\sum_{m\in\Z_1}\big(\|\frac{\partial^3 P^1}{\partial w_n\partial w_n\partial w_m}\frac{\partial F}{\partial\bar{w}_m}\|_{\!{}_{D_+,\Cal O_+}}
+\|\frac{\partial^2 P^1}{\partial w_n\partial w_m}\frac{\partial^2 F}{\partial w_n\partial\bar{w}_m}\|_{\!{}_{D_+,\Cal O_+}}\\
&+&\|\frac{\partial^3 P^1}{\partial w_n\partial w_n\partial \bar w_m}\frac{\partial F}{\partial w_m}\|_{\!{}_{D_+,\Cal O_+}}
+\|\frac{\partial^2 P^1}{\partial w_n\partial \bar w_m}\frac{\partial^2 F}{\partial w_n\partial w_m}\|_{\!{}_{D_+,\Cal O_+}}\big)\\
&+&\|\frac{\partial^3 P^1}{\partial w_n\partial \bar w_n\partial I}\cdot\frac{\partial F}{\partial\tilde{\theta}}\|_{\!{}_{D_+,\Cal O_+}}
+\|\frac{\partial^2 P^1}{\partial w_n\partial I}\cdot\frac{\partial^2 F}{\partial\bar w_n\partial\tilde{\theta}}\|_{\!{}_{D_+,\Cal O_+}}
+\|\frac{\partial^2 P^1}{\partial \bar w_n\partial I}\cdot\frac{\partial^2 F}{\partial w_n\partial\tilde{\theta}}\|_{\!{}_{D_+,\Cal O_+}}\\
&+&\|\frac{\partial P^1}{\partial I}\cdot\frac{\partial^3 F}{\partial w_n\partial\bar w_n\partial\tilde{\theta}}\|_{\!{}_{D_+,\Cal O_+}}
+\|\frac{\partial F_0}{\partial I}\cdot\frac{\partial^3 P^1}{\partial w_n\partial \bar w_n\partial\tilde{\theta}}\|_{\!{}_{D_+,\Cal O_+}}\\
&+&\sum_{m\in\Z_1}\big(\|\frac{\partial^3 P^1}{\partial w_n\partial \bar w_n\partial w_m}\frac{\partial F}{\partial\bar{w}_m}\|_{\!{}_{D_+,\Cal O_+}}
+\|\frac{\partial^2 P^1}{\partial w_n\partial w_m}\frac{\partial^2 F}{\partial \bar w_n\partial\bar{w}_m}\|_{\!{}_{D_+,\Cal O_+}}\\
&+&\|\frac{\partial^2 P^1}{\partial \bar w_n\partial w_m}\frac{\partial^2 F}{\partial w_n\partial\bar{w}_m}\|_{\!{}_{D_+,\Cal O_+}}
+\|\frac{\partial^3 P^1}{\partial w_n\partial \bar w_n\partial \bar w_m}\frac{\partial F}{\partial w_m}\|_{\!{}_{D_+,\Cal O_+}}\\
&+&\|\frac{\partial^2 P^1}{\partial w_n\partial \bar w_m}\frac{\partial^2 F}{\partial \bar w_n\partial w_m}\|_{\!{}_{D_+,\Cal O_+}}
+\|\frac{\partial^2 P^1}{\partial \bar w_n\partial \bar w_m}\frac{\partial^2 F}{\partial w_n\partial w_m}\|_{\!{}_{D_+,\Cal O_+}}\big)\\
&+&\|\frac{\partial^3 P^1}{\partial \bar w_n\partial \bar w_n\partial I}\cdot\frac{\partial F}{\partial\tilde{\theta}}\|_{\!{}_{D_+,\Cal O_+}}
+\|\frac{\partial^2 P^1}{\partial \bar w_n\partial I}\cdot\frac{\partial^2 F}{\partial \bar w_n\partial\tilde{\theta}}\|_{\!{}_{D_+,\Cal O_+}}
+\|\frac{\partial P^1}{\partial I}\cdot\frac{\partial^3 F}{\partial \bar w_n\partial \bar w_n\partial\tilde{\theta}}\|_{\!{}_{D_+,\Cal O_+}}\\
&+&\|\frac{\partial F_0}{\partial I}\cdot\frac{\partial^3 P^1}{\partial \bar w_n\partial \bar w_n\partial\tilde{\theta}}\|_{\!{}_{D_+,\Cal O_+}}
+\sum_{m\in\Z_1}\big(\|\frac{\partial^3 P^1}{\partial \bar w_n\partial \bar w_n\partial w_m}\frac{\partial F}{\partial\bar{w}_m}\|_{\!{}_{D_+,\Cal O_+}}\\
\end{eqnarray*}
\begin{eqnarray*}
&+&\|\frac{\partial^2 P^1}{\partial \bar w_n\partial w_m}\frac{\partial^2 F}{\partial \bar w_n\partial\bar{w}_m}\|_{\!{}_{D_+,\Cal O_+}}
+\|\frac{\partial^3 P^1}{\partial \bar w_n\partial \bar w_n\partial \bar w_m}\frac{\partial F}{\partial w_m}\|_{\!{}_{D_+,\Cal O_+}}+\|\frac{\partial^2 P^1}{\partial \bar w_n\partial \bar w_m}\frac{\partial^2 F}{\partial \bar w_n\partial w_m}\|_{\!{}_{D_+,\Cal O_+}}\big)\Big)\\
&\leq &22\Cal B_\varrho^{\frac{1}{2}}\varepsilon^{2-\beta'}+\frac{7}{|n|}\Cal B_\varrho^{\frac{1}{2}}\varepsilon^{2-\beta'}+6\sum_{m\in Z_1}|m|^{-(a+\bar a)}e^{-2|m|\rho}\Cal B_\varrho^{\frac{1}{2}}\varepsilon^{2-\beta'}+8|n|^{-\bar a-1}e^{-2|n|\rho}\sum_{ |m|\leq E_-K_-\atop m\neq n}|m|^{a}\Cal B_\varrho^{\frac{1}{2}}\varepsilon^{2-\beta'}.
\end{eqnarray*}
According to the above estimate, we have
\begin{eqnarray*}
&&\|\lim_{n\rightarrow \infty}\frac{1}{|n|}\sum_{\upsilon=\pm}\frac{\partial^2\{P^1,F\}}{\partial w^\upsilon_n\partial w^\upsilon_n}\|_{\!{}_{D_+,\Cal O_+}}
\leq 22\Cal B_\varrho^{\frac{1}{2}}\varepsilon^{2-\beta'}+6c\Cal B_\varrho^{\frac{1}{2}}\varepsilon^{2-\beta'}
\leq \eta\varepsilon\leq \varepsilon_+,\\
&&\|\frac{1}{|n|}\sum_{\upsilon=\pm}\frac{\partial^2\{P^1,F\}}{\partial w^\upsilon_n\partial w^\upsilon_n}-\lim_{n\rightarrow \infty}\frac{1}{|n|}\sum_{\upsilon=\pm}\frac{\partial^2\{P^1,F\}}{\partial w^\upsilon_n\partial w^\upsilon_n}\|_{\!{}_{D_+,\Cal O_+}}\\
&\leq& \frac{7}{|n|}\Cal B_\varrho^{\frac{1}{2}}\varepsilon^{2-\beta'}
+8|n|^{-\bar a-1}e^{-2|n|\rho}\sum_{|m|\leq E_-K_-\atop m\neq n}|m|^{a}\Cal B_\varrho^{\frac{1}{2}}\varepsilon^{2-\beta'}
\leq  \frac{c}{|n|}\Cal B_\varrho^{\frac{1}{2}}\varepsilon^{2-\beta'}
\leq\frac{\eta\varepsilon}{|n|}\leq \frac{\varepsilon_+}{|n|}.
\end{eqnarray*}
Together with the above arguments about all the terms in $P_+$, we finally get the verification of $(A6)$ of $P_+$. In this way, associated with the special structure of $P_+$, it is obvious that the form of the normal frequency $\Omega^+_n$ satisfy
$(\ref{asymp1})$ in $(A2)$.

\section{Iteration Lemma and Convergence}

\noindent
Set $0<\beta'\leq\frac{1}{4}$ and $\kappa=\frac{4}{3}-\frac{\beta'}{3}$.
For all $\nu\ge 1$, we
define the following sequences \[
r_\nu=\frac{r_0}{2^\nu},\quad \varrho_\nu=\frac{r_\nu}{20},\quad \rho_\nu=\rho_0(1-\sum_{i=2}^{\nu+1}2^{-i}),\ \ \gamma_\nu=\frac{\gamma_0}{2}(1+2^{-\nu}),\]
\[\Cal B_\nu=\Cal B_{\varrho_\nu}=cE^4_{\nu}\varrho_\nu^{-10(\tilde{b}+\tau+1)},\ \ E_\nu=E_0(2-2^{-\nu}),\]
\beq\varepsilon_\nu=(\varepsilon_0\prod_{\mu=0}^{\nu-1}\Cal B_\mu^{\frac{1}{3\kappa^{\mu+1}}})^{\kappa^\nu},\ \ K_{\nu}=\frac{|\ln\varepsilon_{\nu}|}{\varrho_{\nu}},\eeq
\[\eta^3_\nu=\varepsilon_\nu^{1-\beta'}\Cal B_\nu,\ \ s_{\nu+1}=\eta_\nu s_\nu,\quad D_\nu=D(r_\nu,s_\nu),\]
where $c$ is a constant, and the parameters
$r_0,\varepsilon_0,s_0,\rho_0$  are defined at the beginning of the section 4.

\subsection{Iteration lemma}

\begin{Lemma}\label{Lem5.1}
Suppose that
\begin{eqnarray}
\varepsilon_0\leq (\frac{\delta_0}{80})^{\frac{1}{1-\beta'}}\prod_{\mu=0}^{\infty}\Cal B_\mu^{-\frac{1}{3\kappa^{\mu+1}}},\ \ E_0\bar{\rho}>2\varrho_0,\ \ 3200E_0^2\delta_0<\beta'\gamma_0,\ \ \delta_0\gamma_0\ll\frac{1}{32},\label{Supp}
\end{eqnarray}
and the following conditions \\
\noindent
$(1)$. $\displaystyle N_\nu=\la\bar{\omega},\bar{I}\ra+\la\tilde{\omega}_\nu(\sigma),I\ra+\sum_{n\in\Z_1}
\Omega_{n}^{\nu}(\theta,\sigma)w_n\bar w_n+\sum_{|n|\leq E_{\nu-1}K_{\nu-1}}\langle A_{|n|}^{\nu}z_{|n|},\bar z_{|n|}\rangle$ is a generalized normal form with parameters
$\sigma$
 on a closed set $\Cal O_{\nu}$ of $\R^{\tilde{b}}$;

\noindent $(2)$.  $P_\nu $ has the estimate of the vector field
\[\|X_{P_\nu}\|_{s_\nu,\bar a,\rho_{\nu},D_\nu,\Cal O_{\nu}}\le
\varepsilon_\nu.\] Then there is a subset $\Cal
O_{\nu+1}\subset\Cal O_{\nu}$,
 \[\Cal O_{\nu+1}=\Cal
O_\nu\setminus\bigcup_{k,n,m}(\Cal
R_{k}^{\nu,1}\bigcup\Cal
R_{kn}^{\nu,2}\bigcup\Cal
R_{knm}^{\nu,3}\bigcup\Cal
R_{kn}^{\nu,4}),
\]where
\begin{eqnarray*}
 &&\Cal R_k^{\nu,1}(\gamma_{\nu})=
\left\{\sigma\in \Cal O_\nu:
 \begin{array}{rcl}
 |\la k,\omega_{\nu}\ra|
<\frac{\gamma_\nu}{|k|^\tau},K_{\nu-1} <|k|\leq K_{\nu}\nonumber\\
\end{array} \right\},\\
&&\Cal R_{kn}^{\nu,2}(\gamma_0)=
\left\{\sigma\in \Cal O_\nu:
 \begin{array}{rcl}
 |\la k,\omega_{\nu}\ra\pm (\bar{\Omega}_n^{\nu}+d_n^\nu)|
<\frac{\gamma_0}{K_\nu^\tau},n\in \Z_1, |n|\leq E_\nu K_\nu\nonumber\\
\end{array} \right\},\\
&&\Cal R_{knm}^{\nu,3}(\gamma_{0})=
\left\{\sigma\in \Cal O_\nu:
 \begin{array}{rcl}
 |\la k,\omega_{\nu}\ra\pm ((\bar{\Omega}_n^{\nu}+d_n^\nu)\pm(\bar{\Omega}_m^{\nu}+d_m^\nu))|
<\frac{\gamma_0}{K_\nu^\tau},\ |n|,|m|\leq E_\nu K_\nu\nonumber\\
\end{array} \right\},\\
&&\Cal R_{kn}^{\nu,4}(\gamma_{0})=
\left\{\sigma\in \Cal O_\nu:
 \begin{array}{rcl}
 |\la k,\omega_{\nu}\ra\pm 2\bar{\Omega}_n^{\nu}|
<\frac{\gamma_0\cdot |n|}{K_\nu^\tau},n\in \Z_1, |n|> E_\nu K_\nu\nonumber\\
\end{array} \right\},
\end{eqnarray*}
and a symplectic
transformation of variables $\Phi_\nu:D_{\nu +1}
\times\Cal O_{\nu+1}\to D_{\nu}\times\Cal O_{\nu},$
satisfying
\begin{eqnarray}\label{PHI}
&&\|\Phi_\nu-id\|_{s_\nu,a,\rho,D_{\nu+1}, \Cal O_{\nu+1}},\|D\Phi_\nu-I\|_{s_\nu,a,a,\rho,D_{\nu+1}, \Cal O_{\nu+1}},\nonumber\\
&&\|D\Phi_\nu-I\|_{s_\nu,\bar a,\bar a,\rho,D_{\nu+1},\Cal O_{\nu+1}}\leq \Cal B^{\frac{1}{2}}_\nu\varepsilon_\nu^{1-\beta'},
\end{eqnarray}
such that on
$D(r_{\nu+1},s_{\nu +1})\times\Cal O_{\nu+1},$$H_{\nu+1}=H_\nu\circ\Phi_\nu$ has the form \beq
H_{\nu+1}=\la\bar\omega,\bar I\ra+\la\tilde\omega_{\nu+1},I\ra+\sum_{n\in\Z^d_1}\Omega_{n}^{\nu+1}(\theta,\sigma)w_n\bar
w_n +\sum_{|n|\leq E_{\nu}K_{\nu}}\langle A_{|n|}^{\nu+1}z_{|n|},\bar z_{|n|}\rangle+P_{\nu+1}, \label{5.4}\eeq with \beq
\label{5.5}|\omega_{\nu+1}-\omega_{\nu}|_{\Cal O_{\nu+1}}\le
\Cal B^{\frac{1}{2}}_\nu\varepsilon_{\nu}^{1-\frac{1}{5}\beta'},  \quad
 |\Omega_n^{\nu+1}-\Omega_n^{\nu}|_{-1,D_{\nu+1},\Cal O_{\nu+1}}\le
\Cal B^{\frac{1}{2}}_\nu\varepsilon_{\nu}^{1-\frac{1}{5}\beta'}.\eeq  And also $P_{\nu+1}$ satisfies the estimate
\beq \label{5.6} \|X_{P_{\nu+1}}\|_{s_{\nu+1},\bar a,\rho_{\nu+1},D_{\nu +1},\Cal
O_{\nu+1}}\le\varepsilon_{\nu+1}. \eeq
\end{Lemma}

\proof
From the above iteration formula,and by the definition of $E_0,\gamma_0$, then it is obvious that $E_{\nu}\leq 2E_0,\frac{1}{2}\gamma_0\leq\gamma_\nu\leq \gamma_0$, thus we have  \\
$$(K_{\nu}\varrho_\nu)^{2\tau+2} e^{\frac{8E^2_{\nu}\delta_0 (\gamma_0-\gamma_\nu)K_{\nu} r_\nu}{\gamma_\nu^2}}\leq(|\ln\varepsilon_\nu|)^{2\tau+2}e^{\frac{2560\delta_0E^2_{0}}{\gamma_0}|\ln \varepsilon_\nu|},$$
by $K_{\nu} r_\nu=20 K_{\nu} \sigma_\nu=20|\ln\varepsilon_\nu|$, and choosing $\delta_0$ small enough  and $0<\beta'\leq\frac{1}{4}$ satisfying the inequality defined in $(\ref{Supp})$ such that
\begin{eqnarray}
&&2560\gamma_0^{-1}E^2_0\delta_0<\frac{4}{5}\beta',\ \  e^{\frac{2560\delta_0E^2_{0}}{\gamma_0}|\ln \varepsilon_\nu|}\leq e^{{\frac{4}{5}\beta'}|\ln\varepsilon_\nu|}=\varepsilon_\nu^{-\frac{4}{5}\beta'},\label{e1}\\
&&|\ln\varepsilon_\nu|^{2\tau+2}\leq \varepsilon_\nu^{-\frac{1}{5}\beta'},\  \ \forall  \tau>0, \label{e2}
 \end{eqnarray}so we obtain
$$(K_{\nu}\varrho_\nu)^{2\tau+2} e^{\frac{8E^2_{\nu}\delta_0(\gamma_0-\gamma_\nu) K_{\nu} r_\nu}{\gamma_\nu^2}}\leq \varepsilon_\nu^{-\beta'}.
$$ \\
In view of the definition of $\eta_\nu^3=\varepsilon_\nu^{1-\beta'}\Cal B_\nu$, so if $\varepsilon_\nu^{1-\beta'}\leq \Cal B_\nu^{-1}$, we have
\begin{eqnarray}
\frac{\eta_\nu^2}{\Cal B_\nu(K_{\nu}\varrho_\nu)^{2\tau+2}}e^{-\frac{8E^2_{\nu}\delta_0(\gamma_0-\gamma_\nu) K_{\nu} r_\nu}{\gamma_\nu^2}}\geq \frac{\eta_\nu^2}{\Cal B_\nu}\varepsilon_\nu^{\beta'}=\Cal B_\nu^{-\frac{1}{3}}\varepsilon_\nu^{\frac{2+\beta'}{3}}\geq \varepsilon_\nu.\label{B}
\end{eqnarray}
To verify the inequality $\varepsilon_\nu^{1-\beta'}\leq \Cal B_\nu^{-1}$, since $\Cal B_\nu$ are increasing with $\nu$, then we have
\begin{eqnarray*}
\Cal B_\nu^{\frac{1}{1-\beta'}}=\Cal B_\nu^{\frac{1}{3(\kappa-1)}}=(\prod_{\mu=\nu}^{\infty}\Cal B_\nu^{\frac{1}{3\kappa^{\mu+1}}})^{\kappa^\nu}\leq (\prod_{\mu=\nu}^{\infty}\Cal B_\mu^{\frac{1}{3\kappa^{\mu+1}}})^{\kappa^\nu}.
\end{eqnarray*}
By the definition of $\varepsilon_\nu$ above and the smallness condition on $\varepsilon_0$ defined in $(\ref{Supp})$,
\begin{eqnarray*}
\varepsilon_\nu^{1-\beta'}\Cal B_\nu\leq (\varepsilon_0\prod_{\mu=0}^{\infty}\Cal B_\mu^{\frac{1}{3\kappa^{\mu+1}}})^{\kappa^\nu(1-\beta')}\leq (\frac{\delta_0}{80})^{\kappa^\nu}\leq 1,
\end{eqnarray*}
so the smallness condition in ($\ref{XP}$) is satisfied for any $\nu\geq 0$. In particular,noticing $\kappa\geq \frac{5}{4}$,we have
\begin{eqnarray}
\varepsilon_\nu^{1-\beta'}\Cal B_\nu\leq \frac{\delta_0}{2^{\nu+6}}.\label{B1}
\end{eqnarray}
Now there exists a coordinate transformation $\Phi_\nu:D_{\nu+1}\times\Cal O_{\nu+1}\rightarrow D_\nu\times \Cal O_\nu$ taking $H_\nu$ into $H_{\nu+1}$. Moreover, $(\ref{PHI})$ is obtained by $(\ref{XF})$,$(\ref{DXF})$,$(\ref{4.5})-(\ref{4.7})$, and $(\ref{5.5})$ is obtained by $(\ref{N+})$.
Hence, for $|k|\leq K_\nu$, $$|\la k,\omega_{\nu+1}-\omega_\nu\ra|\leq |k|\cdot|\omega_{\nu+1}-\omega_\nu|\leq K_\nu \Cal B^{\frac{1}{2}}_\nu\varepsilon_{\nu}^{1-\frac{1}{5}\beta'}\leq \Cal{B}^{\frac{1}{2}}_\nu\varepsilon_\nu^{1-\frac{1}{4}\beta'},$$
so we have
\begin{eqnarray*}
|\la k,\omega_{\nu+1}\ra|&\geq& |\la k,\omega_{\nu}\ra|-|\la k,\omega_{\nu+1}-\omega_\nu\ra|\geq \frac{\gamma_\nu}{|k|^\tau}-\Cal{B}^{\frac{1}{2}}_\nu\varepsilon_\nu^{1-\frac{1}{4}\beta'}
\geq  \frac{\gamma_{\nu+1}}{|k|^\tau},\nonumber\\
\end{eqnarray*}
this means the small divisor condition $|\la k,\omega_{\nu+1}\ra|\geq\frac{\gamma_{\nu+1}}{|k|^\tau}$ is automatically satisfied when $|k|\leq K_\nu$.\\
Moreover,we compute some estimates
\begin{eqnarray*}
|\omega_{\nu+1}|_{\Cal O_{\nu+1}}&\overset{(\ref{5.5})}{\leq}& |\omega_{\nu}|_{\Cal O_{\nu}}+|\omega_{\nu+1}-\omega_{\nu}|_{\Cal O_{\nu+1}}\leq E_\nu+\Cal B^{\frac{1}{2}}_\nu\varepsilon_\nu^{1-\frac{1}{5}\beta'}\leq E_{\nu+1},\\
\|\tilde{\Omega}_n^{\nu+1}-\tilde{\Omega}_n^\nu\|_{r_{\nu+1},2\tau+2,\Cal O}&\overset{(\ref{5.5})}{\leq}&
|n|\cdot\Cal B^{\frac{1}{2}}_\nu\varepsilon_\nu^{1-\frac{1}{5}\beta'}
\overset{(\ref{B1})}{\leq} |n|\cdot\frac{\delta_0\gamma_\nu}{2^{\nu+6}}
\leq |n|\cdot(\gamma_\nu-\gamma_{\nu+1})\delta_0,\nonumber\\
\|\tilde{\Omega}_n^{\nu+1}\|_{r_{\nu+1},2\tau+2,\Cal O}&\leq&
\|\tilde{\Omega}_n^{\nu}\|_{r_{\nu},2\tau+2,\Cal O}+\|\tilde{\Omega}_n^{\nu+1}-\tilde
{\Omega}_n^\nu\|_{r_{\nu+1},2\tau+2,\Cal O}\nonumber\\
&\leq& |n|(\gamma_0-\gamma_\nu)\delta_0+|n|(\gamma_\nu-\gamma_{\nu+1})\delta_0\nonumber\\
&\leq& |n|(\gamma_0-\gamma_{\nu+1})\delta_0,
\end{eqnarray*}
this means the assumption $(A1), (A7)$ are also satisfied after one KAM iteration;
\begin{eqnarray}
\frac{cE_\nu^4K_\nu^{4\tau+4}\delta_0}{\eta_\nu^2\gamma_0^4\varrho_\nu^{3\tilde b+3}}e^{-K_\nu\varrho_\nu}\cdot e^{\frac{24 E_\nu^2\delta_0(\gamma_0-\gamma_\nu)K_\nu r_\nu}{\gamma_\nu^2}}\varepsilon_\nu&\leq&\frac{ \Cal{B}_{\nu}}{\eta_\nu^2}|\ln\varepsilon_\nu|^{4\tau+4}e^{-|\ln\varepsilon_\nu|}\cdot e^{{\frac{12}{5}\beta'}|\ln\varepsilon_\nu|}\varepsilon_\nu\nonumber\\
&\leq&\Cal{B}_{\nu}^{\frac{1}{3}}\varepsilon_{\nu}^{\frac{4}{3}-\frac{32}{15}\beta'}\leq \eta_\nu,\label{XRF}\\
\frac{cE_\nu^2K_\nu^{2\tau+2}\delta_0}{\gamma_0^2\varrho_{\nu}^{2\tilde b+1}}e^{-K_\nu\varrho_\nu}\cdot e^{\frac{16 E_\nu^2\delta_0(\gamma_0-\gamma_\nu)K_\nu r_\nu}{\gamma_\nu^2}}&\leq& \Cal B^{\frac{1}{2}}_{\nu}\varepsilon_{\nu}^{-\frac{1}{5}\beta'+1-\frac{8}{5}\beta'}\nonumber\\
&\leq& \Cal B^{\frac{1}{2}}_{\nu}\varepsilon_{\nu}^{1-\frac{9}{5}\beta'}\leq \eta_\nu.\label{XRt}
\end{eqnarray}
Observing that $\frac{E_\nu\bar \rho}{2}\geq \frac{E_0\bar \varrho}{2}$ , it is feasible to choose $E_0$ and $\bar \rho$ satisfying $\frac{E_0\bar \rho}{2}\geq \varrho_0$ defined in $(\ref{Supp})$, then one has
\begin{eqnarray}
c\eta_\nu^{-1}e^{-\frac{E_\nu K_\nu\bar \rho}{2}}&\leq &\Cal B_{\nu}^{-\frac{1}{3}}\varepsilon_{\nu}^{-\frac{1-\beta'}{3}}e^{-K_\nu\varrho_\nu}\leq \Cal B_{\nu}^{-\frac{1}{3}}\varepsilon_{\nu}^{\frac{2+\beta'}{3}}\leq \eta_\nu,\label{XP-R1}\\
c\eta_{\nu}^{-1}e^{-K_\nu\varrho_\nu}&\leq& \Cal B_{\nu}^{-\frac{1}{3}}\varepsilon_{\nu}^{\frac{2+\beta'}{3}}\leq \eta_\nu.\label{XP-R2}
\end{eqnarray}
At last, we estimate the perturbation from $(\ref{XP+})$
\begin{eqnarray*}
&&\|X_{P_{\nu+1}}\|_{s_{\nu+1},\bar a,\rho_{\nu+1},D_{\nu+1},\Cal O_{\nu+1}}\\
&\leq&
\frac{1}{5}(\frac{cE_{\nu}^4K_{\nu}^{4\tau+4}\delta_0}{\eta_{\nu}^2\gamma_0^4\varrho_{\nu}^{3\tilde b+3}}e^{-K_{\nu}\varrho_{\nu}} e^{\frac{24 E_{\nu}^2\delta_0(\gamma_0-\gamma_\nu)K_{\nu}r_{\nu}}{\gamma_{\nu}^2}}\varepsilon_\nu
+\frac{cE_{\nu}^2K_{\nu}^{2\tau+2}\delta_0}{\gamma_0^2\varrho_{\nu}^{2\tilde b+1}}e^{-K_{\nu}\varrho_{\nu}} e^{\frac{16 E_{\nu}^2\delta_0(\gamma_0-\gamma_\nu)K_{\nu}r_{\nu}}{\gamma_{\nu}^2}}\\
&+&c\eta_{\nu}^{-1}e^{-\frac{E_{\nu}K_{\nu}\bar \rho}{2}}+c\eta_{\nu}^{-1}e^{-K_{\nu}\varrho_{\nu}}+c\eta_{\nu})\varepsilon_\nu\nonumber\\
&\leq& \frac{1}{5}(\Cal{B}_{\nu}^{\frac{1}{3}}\varepsilon_{\nu}^{\frac{4}{3}-\frac{32}{15}\beta'}+\Cal B^{\frac{1}{2}}_{\nu}\varepsilon_{\nu}^{1-\frac{9}{5}\beta'}+ 2\Cal B_{\nu}^{-\frac{1}{3}}\varepsilon_{\nu}^{\frac{2+\beta'}{3}}
+\eta_\nu)\varepsilon_\nu\nonumber\\
&\leq& \frac{1}{5}(5\eta_\nu)\varepsilon_\nu= \eta_\nu\varepsilon_\nu=\varepsilon_{\nu+1}.
\end{eqnarray*}
This completes the proof of the iteration lemma.\qed

\subsection{Convergence}

Suppose that the assumptions of Theorem \ref{KAM} are satisfied.
 Recall that
$r_0=r, s_0=s,\rho_0=\rho, N_0=N, P_0=P,E_0=E,\gamma_0=\gamma.$
Define $\delta$ in the KAM Theorem by setting
$$\delta=\delta_0\delta_r,\quad \delta_r=\frac{1}{80}(\prod_{\mu=0}^{\infty}(\Cal B_{\mu})^{-\frac{1}{3\kappa^{\mu+1}}})^{1-\beta'},
$$
where $\delta_r$ depends on $\tilde{b},\tau,r,\gamma,E$ and by the assumption
$$ \varepsilon_0:=\|X_{P_0}\|_{s_0,\bar a,\rho_0,D_0,\Cal O_0}\leq \delta^{\frac{1}{1-\beta'}}.
$$
The small divisor conditions are satisfied by setting
$$\Cal O_1=
\left\{\sigma\in \Cal O_0:
\begin{array}{rcl}
& &|\la k,\omega\ra|\geq \frac{\gamma}{|k|^\tau},\ |k|\ne 0\\
& &|\la k,\omega\ra\pm \Omega_n|\geq \frac{\gamma }{K_0^\tau},\ |n|\leq E_0K_0\\
& &|\la k,\omega\ra\pm 2\Omega_n|\geq \frac{\gamma |n|}{K_0^\tau},\\
\end{array}
\right\},$$
the assumptions of  the iteration lemma  are satisfied when
$\nu=0$ if $\varepsilon_0$ and $\gamma_0$ are sufficiently small.
Inductively, we obtain the following sequences:
\[
\Cal O_{\nu+1}\subset\Cal O_\nu,\]
\[\Psi^\nu=\Phi_0\circ\Phi_1\circ\cdots\circ\Phi_\nu:D_{\nu+1}\times\Cal
O_{\nu+1}\to D_0,\nu\ge 0,
\]
\[H\circ\Psi^\nu=H_{\nu+1}=N_{\nu+1}+P_{\nu+1}.\]
To prove the convergence of the $\Psi^\nu$ we consider the operator norms
$$\|L\|_{s,\tilde{s}}=\sup_{W\neq0}\frac{\|LW\|_{s}}{\|W\|_{\tilde{s}}}.$$
Shorten $ \|\cdot\|_{s,a,\rho}$ as $ \|\cdot\|_{s}$ and these norms satisfy $\|AB\|_{s,\tilde{s}}\leq \|A\|_{s,s}\|B\|_{\tilde{s},\tilde{s}}$ for $s\geq \tilde{s}$ as $\|W\|_{s}\leq \|W\|_{\tilde{s}}$. By the chain rule, we get
\begin{eqnarray*}
\|D\Psi^\nu\|_{s_0,s_{\nu+1},D_{\nu+1},\Cal O_{\nu+1}}&\leq& \prod_{\mu=0}^{\nu}\|D\Phi_\mu\|_{s_{\mu+1},s_{\mu+1},D_{\mu+1},\Cal O_{\mu+1}}\overset{(\ref{PHI}),(\ref{B1})}{\leq} \prod_{\mu=0}^{\infty}(1+\frac{\delta_0}{2^{\mu+6}})\leq 2,\\
\end{eqnarray*}
with the mean value theorem we obtain
\begin{eqnarray*}
\|\Psi^{\nu+1}-\Psi^{\nu}\|_{s_0,D_{\nu+2},\Cal O_{\nu+2}}&\leq& \|D\Psi^\nu\|_{s_0,s_{\nu+1},D_{\nu+1},\Cal O_{\nu+1}}\|\Phi_{\nu+1}-Id\|_{s_{\nu+2},D_{\nu+2}, \Cal O_{\nu+2}}\nonumber\\
&\leq&2\|\Phi_{\nu+1}-Id\|_{s_{\nu+2},D_{\nu+2}, \Cal O_{\nu+2}}\leq 2\Cal B^{\frac{1}{2}}_\nu\varepsilon_\nu^{1-\beta'}.
\end{eqnarray*}
For every non-negative multi-index $k=(k_1,\cdots,k_{\tilde{b}})$, by Cauchy's estimate we have
\begin{eqnarray*}
\|\partial_\theta^{k}(\Psi^{\nu+1}-\Psi^{\nu})\|^{\lambda_0}_{s_0,D_{\nu+3},\Cal O_{\nu+2}}&\leq& 2\Cal B^{\frac{1}{2}}_\nu\varepsilon_\nu^{1-\beta'}\frac{k_1!\cdots k_{\tilde{b}}!}{(\frac{r_0}{2^{\nu+2}})^{|k|}}.
\end{eqnarray*}
The right side of which super-exponentially decay with $\nu$. This shows that $\Psi^\nu$ converge uniformly on $D_*=\T^{\tilde{b}}\times \{0\}\times\{0\}\times\{0\} $ and $\Cal O_\gamma=\bigcap_{\nu\geq 0}\Cal O_\nu$ to a $C^1_W$  continuous family of smooth torus embedding
$$ \Psi:\T^{\tilde{b}}\times \Cal O_\gamma\rightarrow D(r,s).$$
Similarly, the frequencies $\omega_\nu=(\bar \omega,\tilde\omega_{\nu})$ converge uniformly on $\Cal O_\gamma$ to a $C^1_W$ continuous limit $\omega_*=(\bar \omega,\tilde\omega_{*})$, and the frequencies $\Omega_\nu$ converge uniformly on $D_*\times \Cal O_\gamma$ to a regular limit $ \Omega_*$. Moreover, we have the estimate
 \begin{eqnarray*}
&&\|X_H\circ\Psi^\nu-D\Psi^\nu\cdot X_{N_\nu}\|_{s_0,D_{\nu+1},\Cal O_{\gamma}}\\
&\leq& \|D\Psi^\nu\|_{s_0,s_{\nu+1},D_{\nu+1},\Cal O_{\gamma}}\|(\Psi^\nu )^*X_H-X_{N_\nu}\|_{s_{\nu+1},D_{\nu+1},\Cal O_{\gamma}}\nonumber \\
&\leq & c\|X_{P_\nu}\|_{s_{\nu+1},D_{\nu+1},\Cal O_{\gamma}},
\end{eqnarray*}
then $X_H\circ\Psi=D\Psi\cdot X_{N_*}$ on $ D_*$ for each $\sigma\in \Cal O_\gamma$, where $N_*$ is the generalized normal form with frequencies $\omega_*$ and $\Omega_*$. Finally  the Hamiltonian equation becomes
\begin{eqnarray*}
&&\dot{\bar{\theta}}=\bar{\omega},\quad \dot{\bar{I}}=0,\quad \dot{\tilde{\theta}}_j=\tilde{\omega}_{*j},\quad \dot{I}_j=0,\\
&&\dot{w}_n=-{\rm i}(\Omega^*_{n}w_n+a^*_{(-n)n}w_{(-n)}),\quad \dot{\bar w}_n={\rm i}(\Omega^*_{n}\bar w_n+a^*_{n(-n)}\bar w_{(-n)}),
\end{eqnarray*}
where $ \Omega^*_{n}=\bar{\Omega}^*_{n}(\sigma)+\tilde{\Omega}^*_{n}(\theta,\sigma)$. Obviously, we can obtain $ \theta=w_*t$ if we assume the initial value is zero. Then we expand $\tilde{\Omega}^*_{n}(\theta,\sigma)$ into Fourier series
\begin{eqnarray*}
\tilde{\Omega}^*_{n}(\theta,\sigma)=\sum_{k\neq 0}\tilde{\Omega}^{*k}_{n}(\sigma)e^{{\rm i} \la k,\omega_*\ra t},
\end{eqnarray*}
and let
$w_n=\tilde w_ne^{-\sum_{k\neq 0}\frac{\tilde{\Omega}^{*k}_{n}(\sigma)}{\la k,\omega_*\ra}e^{{\rm i} \la k,\omega_*\ra t}},$
then the above equation can be transformed into
\begin{eqnarray*}
&&\dot{\bar{\theta}}=\bar{\omega},\quad \dot{\bar{I}}=0,\quad \dot{\tilde{\theta}}_j=\tilde{\omega}_{*j},\quad \dot{I}_j=0,\\
&&\dot{\tilde w}_n=-{\rm i}(\bar\Omega^*_{n}\tilde w_n+a^*_{(-n)n}\tilde w_{(-n)}),\quad \dot{\bar {\tilde w}}_n={\rm i}(\bar \Omega^*_{n}\bar {\tilde w}_n+a^*_{n(-n)}\bar {\tilde w}_{(-n)}),
\end{eqnarray*}
because $\bar{\Omega}_{*n}(\sigma)$ are all real valued frequencies, $\bar a^*_{(-n)n}=a^*_{n(-n)}$, so the embedded invariant  tori are  linearly stable.
\section{Measure Estimates}

According to the iteration lemma $\ref{Lem5.1}$, we have to exclude
the following resonant set at $\nu^{\rm th}$ step of  KAM iteration
\[\Cal O_{\nu+1}=\Cal
O_{\nu}\setminus\bigcup_{|k|\leq K_\nu}\Cal R_k^\nu,\ \ \nu\geq 0,
\]
$$ \Cal R_k^\nu= \bigcup_{n,m}(\Cal
R_{k}^{\nu,1}\bigcup\Cal
R_{kn}^{\nu,2}\bigcup\Cal
R_{knm}^{\nu,3}\bigcup\Cal
R_{kn}^{\nu,4}),$$
where
\begin{eqnarray*}
&& \Cal R_k^{\nu,1}(\gamma_{\nu})=
\left\{\sigma\in \Cal O_{\nu-1}:
 \begin{array}{rcl}
 |\la k,\omega_{\nu}\ra|
<\frac{\gamma_{\nu}}{|k|^\tau},|k|\geq K_{\nu-1}
\end{array} \right\},\\
&&\Cal R_{kn}^{\nu,2}(\gamma)=
\left\{\sigma\in \Cal O_{\nu-1}:
 \begin{array}{rcl}
 |\la k,\omega_{\nu}\ra\pm (\bar{\Omega}_n^{\nu}+d_n^\nu)|
<\frac{\gamma}{K_{\nu}^\tau},n\in \Z_1, |n|\leq E_{\nu}K_{\nu}
\end{array} \right\},\\
&&\Cal R_{knm}^{\nu,3}(\gamma)=
\left\{\sigma\in \Cal O_{\nu-1}:
 \begin{array}{rcl}
 |\la k,\omega_{\nu}\ra\pm ((\bar{\Omega}_n^{\nu}+d_n^\nu)\pm (\bar{\Omega}_m^{\nu}+d_m^\nu))|
<\frac{\gamma}{K_{\nu}^\tau}, |n|,|m|\leq E_{\nu}K_{\nu}
\end{array} \right\},\\
&&\Cal R_{kn}^{\nu,4}(\gamma)=
\left\{\sigma\in \Cal O_{\nu-1}:
 \begin{array}{rcl}
 |\la k,\omega_{\nu}\ra\pm 2\bar{\Omega}_n^{\nu}|
<\frac{\gamma\cdot |n|}{K_{\nu}^\tau},n\in \Z_1, |n|> E_{\nu}K_{\nu}
\end{array} \right\}.
\end{eqnarray*}
\noindent{\bf Remark.} From the section \ref{no4.3}, one has that
at $\nu^{\rm th}$ step, small divisor condition is automatically
satisfied for $|k|\le K_{\nu-1}$ in the set $ \Cal R_k^{\nu,1}$. Hence, we only need to excise
the above resonant set $\Cal R_k^{\nu,1}$ with $|k|\geq K_{\nu-1}$.

\begin{Lemma}\label{Lem6.1} Let $\tau \geq \tilde{b}$, then the total measure we need to exclude along the KAM iteration is
$${\rm
meas}(\Cal O\setminus\Cal O_\gamma)={\rm meas}(\bigcup_{\nu\ge 0}\bigcup_{|k|\leq K_\nu}\Cal R_{k}^{\nu})<c\gamma.$$
\end{Lemma}
 \proof We firstly give the proof of the most difficult case that the measure estimate of the set $ \Cal R_{knm}^{\nu,3}$
 \[\Cal R_{knm}^{\nu,3}(\gamma)=
\left\{\sigma\in \Cal O_{\nu-1}:
 \begin{array}{rcl}
 |\la k,\omega_{\nu}\ra+ ((\bar{\Omega}_n^{\nu}+d_n^\nu)- (\bar{\Omega}_m^{\nu}+d_m^\nu))|
<\frac{\gamma}{K_{\nu}^\tau}, |n|,|m|\leq E_{\nu}K_{\nu}\nonumber\\
\end{array} \right\}.
 \]
For $\Cal R_{knm}^{\nu,3}$,
according to the assumption $(A2)$, we have
$\bar \Omega_n^\nu=|n|(1+c^\nu(\sigma))$, where $c^\nu(\sigma)$ is independent of $n$ with the estimate
$$|c^\nu(\sigma)|_{\Cal O_{\nu-1}}+|d_n^\nu(\sigma)|_{\Cal O_{\nu-1}}+|d_m^\nu(\sigma)|_{\Cal O_{\nu-1}}=O(\varepsilon_0).$$
Hence, if $|n-m|\geq C|k|$, $C$ is large enough, we have
\begin{eqnarray*}
|\la k,\omega_{\nu}\ra+ ((\bar{\Omega}_n^{\nu}+d_n^\nu)- (\bar{\Omega}_m^{\nu}+d_m^\nu))|\geq |n-m|(1-\varepsilon_0)-c'|k|\geq (\frac{C}{2}-c')|k|\geq \tilde c,
\end{eqnarray*}
in this case there is no small divisor. Hence we only need to consider when $1\leq |n-m|< C|k|$,
\begin{eqnarray*}
|\frac{\partial(\la k,\omega_{\nu}\ra+ ((\bar{\Omega}_n^{\nu}+d_n^\nu)- (\bar{\Omega}_m^{\nu}+d_m^\nu)))}{\partial\sigma}|\geq c'|k|- |n-m|\varepsilon_0\geq (c'-C\varepsilon_0)|k|\geq \tilde c'
\end{eqnarray*}
with $\varepsilon_0\ll c'$  and we get the lower bound of the partial derivative about $\la k,\omega_{\nu}\ra+ ((\bar{\Omega}_n^{\nu}+d_n^\nu)- (\bar{\Omega}_m^{\nu}+d_m^\nu))$.  Therefore, for any fixed $|k|\leq K_\nu, |n|,|m|\leq E_\nu K_\nu$, we obtain
\begin{eqnarray*}
{\rm meas}(\Cal R_{knm}^{\nu,3})\leq c\frac{\gamma_0}{K_\nu^{\tau}}.
\end{eqnarray*}
Similarly we can get the estimates of the  sets
\begin{eqnarray*}
{\rm meas}(\Cal R_{kn}^{\nu,2})\leq c\frac{\gamma_0}{K_\nu^{\tau}},\ \ {\rm meas}(\Cal R_{kn}^{\nu,4})\leq c\frac{\gamma_0}{K_\nu^{\tau}}.
\end{eqnarray*}
For the set $\Cal R_{k}^{\nu,1}$, we have $|\frac{\partial\la k,\omega\ra}{\partial\sigma}|\geq c|k|$ and for any fixed $K_{\nu-1}<|k|\leq K_\nu$, the estimate of $\Cal R_{k}^{\nu,1}$
$${\rm meas}(\Cal R_{k}^{\nu,1})\leq c\frac{\gamma_\nu}{\la k\ra^{\tau+1}}.$$
Therefore, we get
\begin{eqnarray*}
{\rm meas}(\Cal R_{k}^{\nu})&\leq& {\rm meas}(\bigcup_{n,m\leq E_\nu K_\nu}(\Cal
R_{k}^{\nu,1}\bigcup\Cal
R_{kn}^{\nu,2}\bigcup\Cal
R_{knm}^{\nu,3}\bigcup\Cal
R_{kn}^{\nu,4}))\\
&\leq& c\frac{\gamma_\nu}{|k|^{\tau+1}}
+\sum_{n,m\leq E_\nu K_\nu}c\frac{\gamma}{K_\nu^{\tau}}
\leq c\frac{\gamma_\nu}{|k|^{\tau+1}} +c\frac{\gamma}{K_\nu^{\tau+2}},\\
{\rm meas}(\Cal O\setminus\Cal O_\gamma)&=&{\rm meas}(\bigcup_{\nu\ge 0}\bigcup_{|k|\leq K_\nu}\Cal R_{k}^{\nu})\\
&\leq& \sum_{\nu\geq0}(\sum_{K_{\nu-1}<|k|\leq  K_\nu}c\frac{\gamma_\nu}{|k|^{\tau+1}}
+\sum_{|k|\leq  K_\nu}c\frac{\gamma}{K_\nu^{\tau+2}})\\
&\leq& C(\tilde b,\tau)\sum_{\nu\geq0}\frac{\gamma_0}{K_{\nu-1}^{\tau+1-\tilde b}}+c\frac{\gamma}{K_\nu^{\tau+2-\tilde b}}\\
&\leq &c\gamma,
\end{eqnarray*}
the sum of the former inequality over all $\nu$ converges if $\tau+1>\tilde b$ and we finally obtain the measure estimate.
\qed

\section{Appendix}

\sss
\begin{Lemma}\label{Lem7.1}
$$\|FG\|_{ D(r,s)  }\le \|F\|_{ D(r,s)  }\|G\|_{ D(r,s)  }.$$\end{Lemma}

\proof Since
$(FG)_{kl\alpha\beta}=\sum_{k',l',\alpha',\beta'}F_{k-k',l-l',\alpha-\alpha',\beta-\beta'}
G_{k'l'\alpha'\beta'}$, we have
\begin{eqnarray*}
\|FG\|_{ D(r,s) }&=&\sup_{{\|w\|_{a,\rho}<s}\atop{\|\bar
w\|_{a,\rho}<s}}\sum_{k,l,\alpha,\beta}|(FG)_{kl\alpha\beta}|s^{2l}|w^{\alpha}|
|\bar
w^{\beta}|e^{|k|r}\\
&\le&\sup_{{\|w\|_{a,\rho}<s}\atop{\|\bar w\|_{a,\rho}<s}}
\sum_{k,l,\alpha,\beta}\sum_{k',l',\alpha',\beta'}|F_{k-k',l-l',\alpha-\alpha',\beta-\beta'
}G_{k'l'\alpha'\beta'}|s^{2l}|w^{\alpha}| |\bar
w^{\beta}|e^{|k|r}\\
&\le&\|F\|_{ D(r,s)  }\|G\|_{ D(r,s)  }
\end{eqnarray*}
and the proof is finished.\qed

\begin{Lemma}\label{Lem7.2}
(Cauchy inequalities)
$$
\|F_{\theta}\|_{D(r-\sigma,s)}\le \frac{c}{\sigma}\|F\|_{ D(r,s)
},\ \
\|F_{I}\|_{D(r,\frac 12 s)}\le \frac {c}{s^2}\|F\|_{ D(r,s) },$$
and
$$
\|F_{w_n}\|_{D(r,\frac 12 s)}\le \frac{c}{s}|n|^{a}e^{|n|\rho}\|F\|_{ D(r,s)
},\ \
\|F_{\bar w_n}\|_{D(r,\frac 12 s)}\le \frac{c}{s}|n|^{a}e^{|n|\rho}\|F\|_{ D(r,s)
}.$$
\end{Lemma}
\begin{Lemma}\label{Lem7.3}
There exists a constant $c>0$ such that if $n\in \Z_1,\rho>0$,
$$\|F_n\|_{D(r,s)}<ce^{-|n|\rho},\ \ \|G\|_{D(r,s)}<\varepsilon,$$
then $$\|\{F_n,G\}\|_{D(r-\sigma,\frac{1}{2}s)}<c\sigma^{-1}s^{-2}\|F_n\|_{D(r,s)}\|G\|_{D(r,s)}\leq c\sigma^{-1}s^{-2}\varepsilon e^{-|n|\rho}.$$
\end{Lemma}
\proof According to Lemma 7.1 and 7.2,
$$\|\la F_{n_I},G_\theta\ra\|_{D(r-\sigma,\frac{1}{2}s)}<c\sigma^{-1}s^{-2}\|F_n\|\cdot\|G\|,$$
$$\|\la F_{n_\theta},G_I\ra\|_{D(r-\sigma,\frac{1}{2}s)}<c\sigma^{-1}s^{-2}\|F_n\|\cdot\|G\|,$$
\begin{eqnarray}
\|\sum _{m}F_{n_{w_m}}G_{\bar w_m}\|_{D(r,\frac{1}{2}s)}&\leq &\sum_{m}\|F_{n_{w_m}}\|_{D(r,\frac{1}{2}s)}\|G_{\bar w_m}\|_{D(r,\frac{1}{2}s)}\nonumber\\
&\leq & \|F_{n_{w}}\|_{D(r,\frac{1}{2}s)}\|G_{\bar w}\|_{D(r,\frac{1}{2}s)}\nonumber \\
&\leq & cs^{-2}\|F_n\|\cdot\|G\|.\nonumber
\end{eqnarray}
It follows that
$\|\{F_n,G\}\|_{D(r-\sigma,\frac{1}{2}s)}<c\sigma^{-1}s^{-2}\varepsilon e^{-|n|\rho}$
.\qed

\end{document}